\documentclass[12pt]{article}

\usepackage{amsmath,enumerate,amsbsy,amsfonts,amssymb,mathabx,amscd,graphicx,algorithm}
\usepackage[round,authoryear]{natbib}
\usepackage{authblk}
\usepackage{mathrsfs}
\usepackage{adjustbox}
\usepackage{epsfig}
\newtheorem{definition}{Definition}
\newtheorem{theorem}{Theorem}

\newtheorem{lemma}{Lemma}

\newtheorem{proposition}{Proposition}
\usepackage{multirow}

\newcommand{\beps}{\boldsymbol \epsilon}

\newcommand{\bSigma}{\boldsymbol \Sigma}

\newcommand{\bsigma}{\boldsymbol \sigma}
\newcommand{\bGamma}{\boldsymbol \Gamma}

\newcommand{\bDelta}{\boldsymbol \Delta}
\newcommand{\bXi}{\boldsymbol \Sigma}
\newcommand{\bpsi}{\boldsymbol \psi}
\newcommand{\bPsi}{\boldsymbol \Psi}
\newcommand{\bphi}{\boldsymbol \phi}
\newcommand{\bPhi}{\boldsymbol \Phi}

\newcommand{\bxi}{\boldsymbol \xi}
\newcommand{\btheta}{\boldsymbol \theta}
\newcommand{\bgamma}{\boldsymbol \gamma}
\newcommand{\bdelta}{\boldsymbol \delta}

\newcommand{\bfeta}{\boldsymbol \eta}
\newcommand{\bvar}{\boldsymbol \varepsilon}

\newcommand{\bvarepsilon}{\boldsymbol \varepsilon}
\newcommand{\bzeta}{\boldsymbol \zeta}

\newcommand{\ba}{{\mathbf a}}
\newcommand{\be}{{\mathbf e}}
\newcommand{\bbf}{{\mathbf f}}
\newcommand{\bx}{{\mathbf x}}
\newcommand{\by}{{\mathbf y}}

\newcommand{\bv}{{\mathbf v}}

\newcommand{\br}{{\mathbf r}}
\newcommand{\bw}{{\mathbf w}}
\newcommand{\bg}{{\mathbf g}}
\newcommand{\bs}{{\mathbf s}}
\newcommand{\bbb}{{\mathbf b}}

\newcommand{\bD}{{\bf D}}
\newcommand{\bA}{{\bf A}}
\newcommand{\bB}{{\bf B}}
\newcommand{\bC}{{\bf C}}
\newcommand{\bE}{{\bf E}}
\newcommand{\bI}{{\bf I}}

\newcommand{\bP}{{\bf P}}
\newcommand{\bS}{{\bf S}}
\newcommand{\bX}{{\bf X}}
\newcommand{\bY}{{\bf Y}}
\newcommand{\bZ}{{\bf Z}}
\newcommand{\bR}{{\bf R}}
\newcommand{\bU}{{\bf U}}
\newcommand{\bV}{{\bf V}}
\newcommand{\bW}{{\bf W}}
\newcommand{\bQ}{{\bf Q}}
\newcommand{\bJ}{{\bf J}}

\newcommand{\cD}{{\cal D}}

\newcommand{\cL}{{\cal L}}
\newcommand{\cM}{{\cal M}}

\newcommand{\cU}{{\cal U}}
\newcommand{\cS}{{\cal S}}

\newcommand{\mF}{f}
\newcommand{\eZ}{\mathbb{Z}}
\newcommand{\eR}{\mathbb{R}}
\newcommand{\cH}{\mathbb{H}}
\newcommand{\eS}{\mathbb{S}}

\newcommand{\tA}{\text{A}}

\newcommand{\tE}{\text{E}}
\newcommand{\tF}{\text{F}}

\newcommand{\tX}{\text{X}}
\newcommand{\tR}{\text{R}}

\newcommand{\cov}{\text{Cov}}
\newcommand{\var}{\text{Var}}

\newcommand{\LS}{\text{LS}}

\newcommand{\aic}{\text{AIC}}
\newcommand{\bic}{\text{BIC}}

\newcommand{\llangle}{\langle \langle}
\newcommand{\rrangle}{\rangle \rangle}

\newcommand{\0}{{\bf 0}}
\def\T{{ \mathrm{\scriptscriptstyle T} }}
\newtheorem{condition}{Condition}

\makeatletter
\DeclareRobustCommand\widecheck[1]{{\mathpalette\@widecheck{#1}}}
\def\@widecheck#1#2{%
	\setbox\z@\hbox{\m@th$#1#2$}%
	\setbox\tw@\hbox{\m@th$#1%
		\widehat{%
			\vrule\@width\z@\@height\ht\z@
			\vrule\@height\z@\@width\wd\z@}$}%
	\dp\tw@-\ht\z@
	\@tempdima\ht\z@ \advance\@tempdima2\ht\tw@ \divide\@tempdima\thr@@
	\setbox\tw@\hbox{%
		\raise\@tempdima\hbox{\scalebox{1}[-1]{\lower\@tempdima\box
				\tw@}}}%
	{\ooalign{\box\tw@ \cr \box\z@}}}
\makeatother
\RequirePackage[colorlinks,citecolor=blue,urlcolor=blue]{hyperref}

\begin{document}

\title{\bf \large On Consistency and Sparsity for High-Dimensional Functional Time Series with Application to Autoregressions}


\author[2]{Shaojun Guo}
\author[1]{Xinghao Qiao}
\affil[1]{Department of Statistics, London School of Economics and Political Science, U.K.}
\affil[2]{Institute of Statistics and Big Data, Renmin University of China, P.R. China}
\setcounter{Maxaffil}{0}

\renewcommand\Affilfont{\itshape\small}
\date{\vspace{-5ex}}
\maketitle
\begin{abstract}
Modelling a large collection of functional time series arises in a broad spectral of real applications. Under such a scenario, not only the number of functional variables can be diverging with, or even larger than the number of temporally dependent functional observations,  but each function itself is an infinite-dimensional object, posing a challenging task. In this paper, we propose a three-step procedure to estimate high-dimensional functional time series models. To provide theoretical guarantees for the three-step procedure, we focus on multivariate stationary processes and propose a novel functional stability measure based on their spectral properties. Such stability measure facilitates the development of some useful concentration bounds on sample (auto)covariance functions, which serve as a fundamental tool for further convergence analysis in high-dimensional settings. As functional principal component analysis (FPCA) is one of the key dimension reduction techniques in the first step, we also investigate the non-asymptotic properties of the relevant estimated terms under a FPCA framework. To illustrate with an important application, we consider vector functional autoregressive models and develop a regularization approach to estimate autoregressive coefficient functions under the sparsity constraint. Using our derived non-asymptotic results, we investigate convergence properties of the regularized estimate under high-dimensional scaling. Finally, the finite-sample performance of the proposed method is examined through both simulations and a public financial dataset.
\end{abstract}


\section{Introduction}
\label{sec.intro}
In functional data analysis, it is commonly assumed that each measured function, treated as the unit of observation, is independently sampled from some realization of an underlying stochastic process. Functional time series, on the other hand, refers to a collection of curves observed consecutively over time, where the temporal dependence across observations exhibits. 
The literature has mainly focused on univariate or bivariate functional time series, see, e.g., \cite{hormann2010,cho2013,panaretos2013,hormann2015,jirak2016,li2019} and the reference therein. Recent advances in technology have made multivariate or even high-dimensional functional time series datasets become increasingly common in many applications. Examples include cumulative intraday return trajectories \cite[]{horvath2014} and functional volatility processes \cite[]{muller2011} for a large number of stocks, daily concentration curves of particulate matter and gaseous pollutants at different locations \cite[]{li2017}, and intraday energy consumption curves for thousands of London households (available at {https://data.london.gov.uk/dataset/smartmeter-energy-use-data-in-london-households}).
These applications require understanding relationships among a relatively large collection of functional variables based on temporally dependent functional observations.

Throughout the paper, suppose we observe $p$-dimensional vector of functional time series, $\bX_t(\cdot)=\big(X_{t1}(\cdot), \dots, X_{tp}(\cdot)\big)^{\T}$ for $t=1,\dots,n,$ defined on a compact interval $\cU.$ Addressing multivariate or even high-dimensional functional time series problems poses challenges and is largely untouched in the literature. Under such a scenario, not only $p$ is large relative to $n,$ but each $X_{tj}(\cdot)$ is an infinite-dimensional object with temporal dependence across observations. A standard procedure towards the estimation of models involving high-dimensional functional data consists of three steps.
In the first step, due to the infinite-dimensional nature of functional data, some form of dimension reduction, e.g. data-driven basis expansion via {functional principal component analysis} (FPCA) or its dynamic version \cite[]{hormann2015} and pre-fixed basis expansion \cite[]{fan2015}, 
is needed to approximate each $X_{tj}(\cdot)$ by the $q_j$-dimensional truncation, which transforms the problem of modelling a $p$-dimensional vector of functional time series into that of modelling ($\sum_{j=1}^p q_j$)-dimensional vector time series. The second step involves the estimation under a high-dimensional and dependent setting, 
where some lower-dimensional structure is commonly imposed on the model parameter space. 
One large class assumes sparse function-valued parameters involved in high-dimensional functional time series models. Under functional sparsity constraints, the first step results in the estimation of block sparse vector- or matrix-valued parameters in the second step, where different regularized estimation procedures can be developed in a blockwise fashion, see, e.g., under an independent setting, \cite{fan2015,kong2016} and \cite{qiao2018a}. 
Finally, for interpretation and prediction, the third step recovers functional sparse estimates from those block sparse estimates obtained in the second step. 

The essential challenge to support such three-step estimation procedure is to provide theoretical guarantees in a range of high-dimensional settings, e.g., $\log p/n \rightarrow 0.$ Within such high-dimensional statistics framework, one main goal is to obtain some non-asymptotic results, i.e. error bounds on a given performance metric that hold with high probability for a finite sample size $n$ and provide explicit control on the dimension $p$ as well as other structural parameters. 
Compared with non-functional data, the intrinsic infinite-dimensionality of each process $X_{tj}(\cdot)$ leads to a significant rise in theoretical complexity of the problem, since one needs to develop some operator- and FPCA-based non-asymptotic results for dependent processes within an abstract Hilbert space and to propose a dependence measure to capture the effect of temporal dependence on non-asymptotic properties.
The existing theoretical work mainly for the first step has focused on studying its asymptotic properties by treating $p$ as fixed under a moment-based dependence structure \cite[]{hormann2010} or its non-asymptotic properties under either an independent setting \cite[]{koltch2017,qiao2019,av2020} or a special autoregressive structure \cite[]{Bbosq1}. These results, however, are not sufficient to evaluate the performance of the three-step procedure in a high-dimensional regime with a general dependence structure. Such a challenging task motivates us to develop some essential non-asymptotic results under the setting we consider, which fills the gap between practical implementation and theoretical justification and forms the core of our paper.

A key innovation in our paper is to propose a functional stability measure for a large class of stationary Gaussian processes, $\{\bX_t(\cdot)\},$ based on their spectral density functions. Such stability measure provides new insights into the effect of temporal dependence on theoretical properties of $\widehat\bSigma_h$'s, the estimators for autocovariance functions $\bSigma_h(u,v) = \cov\{\bX_t(u),\bX_{t+h}(v)\}$ with $h=0, \pm 1, \dots$ and $(u,v)\in \cU^2,$  and, in particular, facilitates the development of some novel concentration bounds on $\widehat\bSigma_h$ serving as a fundamental tool for further convergence analysis under high-dimensional scaling. Based on these concentration bounds, we establish non-asymptotic error bounds on relevant estimated terms under a FPCA framework so as to provide theoretical guarantees for our proposed three-step procedure. Such concentration results can also lead to convergence analysis of other possible high-dimensional functional time series models, e.g. those mentioned in Section~\ref{sec.discussion}. It is worth noting that the functional stability measure is fundamentally different from the direct extension of the stability measure \cite[]{basu2015a} to the functional domain. This is because, for truly infinite-dimensional functional objects, in contrast to the functional analog of Basu et al.'s stability measure, which just controls the maximum eigenvalue for spectral density functions of $\{\bX_t(\cdot)\},$ the functional stability measure utilizes the functional Rayleigh quotients of spectral density functions relative to $\bSigma_0$ and hence can more precisely capture the effect of small decaying eigenvalues. Such functional stability measure also leads to non-asymptotic results of normalized versions of relevant estimated terms, making the characterization of relevant tail behaviours more accurate for small eigenvalues. 
To the best of our knowledge, we are the first to propose such a dependence measure for high-dimensional functional time series and rely on it to develop some essential non-asymptotic results. 
To illustrate the proposed three-step approach and the usefulness of the derived non-asymptotic results with an important application, we consider {vector functional autoregressive} (VFAR) models, which characterize the temporal and cross-sectional inter-relationships in $\{\bX_t(\cdot)\}.$ 
One advantage of a VFAR model is that it accommodates dynamic linear interdependencies in $\{\bX_t(\cdot)\}$ into a static framework within a Hilbert space. 
Moreover, a sparse VFAR model facilitates the extraction of Granger causal networks \cite[]{basu2015b} under the functional domain. The VFAR estimation is intrinsically a very high-dimensional problem, since, in the second step of our procedure, we need to fit a {\it vector autoregressive} (VAR) model, whose dimensionality, $(\sum_{j=1}^p q_j)^2,$ grows quadratically with $\sum_{j=1}^p q_j.$ For example, estimating a VFAR model of order 1 with $p=20$ and $q_j=5$ requires estimating $100^2=10,000$ parameters. 
Under high-dimensional scaling and sparsity assumptions on the functional transition matrices, the second step
requires to estimate a block sparse VAR model. We then propose the regularized estimates of block transition matrices, on which the block sparsity constraint is enforced via a standardized group lasso penalty \cite[]{simon2012}. Using the derived non-asymptotic results, we show that the proposed three-step approach can produce consistent estimates for sparse VFAR models in high dimensions.

{\bf Related literature}. 
Our work lies in the intersection of high-dimensional statistics, functional data analysis and time series analysis, each of which corresponds to a vast literature, hence we will only review the most closely related intersectional work to ours. 
(i) For high-dimensional independent functional data, regularization methods have recently been proposed to estimate different types of functional sparse models, e.g., functional additive regression \cite[]{fan2015,kong2016,xue2020}, static functional graphical models \cite[]{li2018,qiao2018a} and its dynamic version \cite[]{qiao2019}. 
(ii) For high-dimensional time series, some essential concentration bounds were established for 
Gaussian processes \cite[]{basu2015a}, linear processes with more general noise distributions \cite[]{sun2018} and heavy tailed non-Gaussian processes \cite[]{wong2020}.
(iii) For examples of recent developments in high-dimensional VAR models, 
\cite{kock2015,basu2015a} and \cite{wong2020} studied the theoretical properties of $\ell_1$-type regularized estimates. 
\cite{basu2015b} and \cite{billio2019} considered extracting Granger causal networks from sparse VAR models.
See also \cite{han2015}, \cite{guo2016} and \cite{ghosh2019}.
(iv) For examples of research on functional autoregressive models, see \cite{Bbosq1,koko2013,aue2015}  \cite{kowal2019} 
and the reference therein.

{\bf Outline of the paper}. In Section~\ref{sec.main}, we first introduce a functional stability measure and rely on it to establish concentration bounds on the quadratic/bilinear forms of $\widehat\bSigma_h$ leading to elementwise concentration results for $\widehat\bSigma_h$. We then establish theoretical guarantees for the proposed three-step approach by 
deriving some useful non-asymptotic error bounds under a FPCA framework. In Section~\ref{sec.vfar}, we 
develop a three-step procedure to estimate the sparse VFAR model, connect with casual network modelling, present the convergence analysis of the regularized estimate and finally examine the finite-sample performance through both simulation studies and an analysis of a public financial dataset. In Section~\ref{sec.discussion}, we conclude our paper by discussing several potential extensions. 
All technical proofs are relegated to the Supplementary Material.

{\bf Notation}. We summarize here some notation to be used throughout the paper. Let $\eZ$ 
denotes the set of integers. 
For two positive sequences $\{a_n\}$ and $\{b_n\}$, 
we write $a_n \lesssim b_n$ or $b_n \gtrsim a_n$ if there exists an absolute constant $c,$ such that $a_n \leq c b_n$ for all $n.$ We write $a_n \asymp b_n$ if and only if $b_n \lesssim a_n$ and $a_n \lesssim b_n.$  We use $x \vee y = \max(x,y).$
For matrices $\bA, \bB \in {\eR}^{p_1 \times p_2},$ we let $\langle\langle \bA,\bB\rangle\rangle=\text{trace}(\bA^{\T}\bB)$ and denote the Frobenius 
norm of $\bB$ by $||\bB||_{\tF}=\big(\sum_{j,k} \text{B}_{jk}^2\big)^{1/2}.$ 
Let $L_2(\cU)$ denote a Hilbert space of square integrable functions defined on the compact set $\cU$ equipped with the inner product $\langle f,g \rangle=\int_{\cU}f(u)g(u)du$ for $f,g \in L_2(\cU)$ and the induced norm $\|\cdot\|=\langle \cdot,\cdot \rangle^{1/2}.$ We denote its $p$-fold Cartesian product by $\cH=L_2(\cU) \times  \dots \times L_2(\cU)$ and the tensor product by $\eS=L_2(\cU)\otimes L_2(\cU).$ For $\bbf=(f_1, \ldots, f_p)^{\T}$ and $\bg=(g_1, \dots,g_p)^{\T}$ in $\cH,$ we denote the inner product by
$\langle\bbf,\bg\rangle_{\cH}=\sum_{j=1}^p\langle f_j, g_j\rangle$ and the induced norm by  $\|\cdot\|_{\cH}=\langle \cdot,\cdot \rangle^{1/2}_{\cH}.$ We use $\|\bbf\|_0=\sum_{j=1}^p I(\|f_j\| \neq 0)$ 
to denote the functional version of vector $\ell_0$ norm. 
For any $K \in \eS,$ it can be viewed as the kernel function of a linear operator acting on $L_2(\cU),$ i.e. for each $f \in L_2(\cU),$
$K$ maps $f(u)$ to $K(f)(u)=\int_{\cU}K(u,v)f(v)dv.$
For notational economy, we will use $K$ to denote both the kernel function and the operator.
Moreover, we denote the operator and Hilbert-Schmidt norms by $\|K\|_{\cL}=\sup_{\|f\| \leq 1}\|K(f)\|$ and 
$\|K\|_{\cS} =\big(\int\int K(u,v)^2dudv\big)^{1/2},$ respectively. For 
$\bA=\big(\tA_{jk}\big)_{1 \leq j,k \leq p}$ with its $(j,k)$-th component $\tA_{jk} \in \eS$, 
we define the functional versions of 
Frobenius, elementwise $\ell_{\infty}$ and matrix $\ell_{\infty}$ norms by 
$\|\bA\|_{\tF} = \big(\sum_{j,k} \|\tA_{jk}\|_{\cS}^2\big)^{1/2},$ $\|\bA\|_{\max} =  \max_{j,k} \|\tA_{jk}\|_{\cS}$ and $\|\bA\|_{\infty} = \max_{j} \sum_{k} \|\tA_{jk}\|_{\cS},
$
respectively. 

\section{Main results}
\label{sec.main}
Suppose that $\{\bX_t(\cdot)\}_{t \in \eZ},$ defined on $\cU,$ is a sequence of $p$-dimensional vector of centered and covariance-stationary Gaussian processes with mean zero and $p \times p$ autocovariance functions,  $\bSigma_{h}=\big(\Sigma_{jk}^{(h)}\big)_{1\leq j,k\leq p}$ with its $(j,k)$-th component $\Sigma_{jk}^{(h)}\in \eS$ for $h \in \eZ.$
In particular when $h=0,$ one typically refers to $\Sigma_{jj}^{(0)}$ as marginal-covariance functions for $j=k,$ and cross-covariance functions for $j \neq k.$ To simplify notation, we will also use $\bSigma_{h}$ to denote the lag-$h$ autocovariance operator induced from the kernel function $\bSigma_h,$ i.e., for any given $\bPhi \in \cH,$
$$
    \bSigma_h(\bPhi)(u): = \int_{\cU} \bSigma_{h}(u,v) \bPhi(v)dv= \Big(\big\langle\bsigma_{1}^{(h)}(u,\cdot),\bPhi(\cdot)\big\rangle_{\cH}, \dots, \big\langle\bsigma_{p}^{(h)}(u,\cdot),\bPhi(\cdot)\big\rangle_{\cH}\Big)^{\T}
    \in \cH,
$$
where $\bsigma_{j}^{(h)}(u,\cdot)=\big(\Sigma_{j1}^{(h)}(u,\cdot), \dots, \Sigma_{jp}^{(h)}(u,\cdot)\big)^{\T}$ for $j=1, \dots, p.$ 
In the special case of $h=0,$ the covariance function $\bSigma_0$ is symmetric and non-negative definite, i.e. $\bSigma_0(u,v)=\bSigma_0(v,u)^{\T}$ for any $(u,v) \in \cU^2$ and $\langle\bPhi, \bSigma_0(\bPhi)\rangle_{\cH} \geq 0$ for any $\bPhi\in \cH.$  
\subsection{Functional stability measure}
\label{sec.fsm}

To introduce the functional stability measure, we first consider the {spectral density operator} of  $\{\bX_t(\cdot)\}_{t \in \eZ},$ defined from the  Fourier transform of autocovariance operators $\{\bSigma_{h}\}_{h \in \eZ},$ which encodes the second-order dynamical properties of $\{\bX_t(\cdot)\}_{t \in \eZ}.$

\begin{definition}
\label{def.sd}
We define the spectral density operator  of $\{\bX_t(\cdot)\}_{t \in \eZ}$ at frequency $\theta$ by
\begin{equation}
\label{def.sd.function}
	f_{\bX,\theta} = \frac{1}{2 \pi} \sum_{h \in \eZ} \bSigma_{h} \exp(- i h \theta), ~~\theta \in [-\pi,\pi].
\end{equation}
\end{definition}

The spectral density operator (or function) generalizes the notion of the spectral density matrix \cite[]{basu2015a} to the functional domain, and it can also be viewed as a generalization of the spectral density operator (or function) \cite[]{panaretos2013} to the multivariate setting. 
Furthermore, if $\sum_{h =0}^\infty \|\bXi_{h}\|_{\cL} < \infty,$ then ${\mF}_{\bX,\theta}$ is uniformly bounded and continuous in $\theta$ with respect to $\|\cdot\|_{\cL},$ where we denote by $\|\bXi_{h}\|_{\cL}={\text{sup}}_{\|\bPhi\|_{\cH} \leq 1, \bPhi \in \cH} \|\bXi_{h}(\bPhi)\|_{\cH}$ the operator norm of $\bXi_{h},$ and the following inversion formula holds:
\begin{equation}
\label{inv.formula}
	\bXi_{h} = \int_{-\pi}^\pi {\mF}_{\bX,\theta} \exp( i h \theta)d\theta, ~~~\text{ for all }h \in \eZ.
\end{equation}
The inversion relationships in (\ref{def.sd.function}) and (\ref{inv.formula}) indicate that spectral density operators  and autocovariance operators comprise a Fourier transform pair. Hence, to study the second-order dynamics of $\{\bX_t(\cdot)\}_{t \in \eZ},$ we can impose conditions on $\bSigma_0$ and $\{\mF_{\bX, \theta}, \theta \in [-\pi,\pi]\}$  in the following Conditions \ref{cond.cov.func} and \ref{cond.bd.fsm}, respectively, which together imply that
$\{\mF_{\bX,\theta}, \theta \in [-\pi,\pi]\}$ are trace-class operators.

\begin{condition}
	\label{cond.cov.func}
	(i) The marginal-covariance functions, $\Sigma_{jj}^{(0)}$'s, are continuous on $\cU^2$ and uniformly bounded over $j\in\{1,\ldots,p\};$
	(ii) $\lambda_0 = {\max}_{1\leq j \leq p}\int_{\cU} \Sigma_{jj}^{(0)}(u,u)du =O(1)$.
\end{condition}

\begin{condition}
\label{cond.bd.fsm}
(i) The spectral density operators $\mF_{\bX,\theta}, \theta \in [-\pi,\pi]$ exist;
(ii) The functional stability measure of $\{\bX_t(\cdot)\}_{t \in \eZ}$ 
defined as follows, is bounded, i.e.
\begin{equation}
\label{eq.bd.fsm}
	\cM(\mF_{\bX}) =  2 \pi \cdot \underset{\theta\in [-\pi, \pi], \bPhi \in \cH_0^p}{\text{esssup}} \frac{ \big\langle \bPhi, \mF_{\bX,\theta}(\bPhi) \big\rangle_{\cH}}{\big\langle\bPhi, \bXi_0(\bPhi)\big\rangle_{\cH}} < \infty,
\end{equation}
where $\cH_0 = \big\{\bPhi \in \cH :  \langle\bPhi, \bXi_0(\bPhi)\rangle_{\cH} \in (0,\infty) \big\}.$
\end{condition}

In general, we can relax Condition~\ref{cond.cov.func}(ii) by allowing $\lambda_0$ to grow at some slow rate as $p$ increases. Then our established non-asymptotic bounds, e.g., those in (\ref{bd_max_lam_phi}) and (\ref{bd_score_max}), will depend on $\lambda_0.$
We next provide several comments for Condition~\ref{cond.bd.fsm}. First, the functional stability measure is proportional to the essential supremum of the functional Rayleigh quotient of $\mF_{\bX,\btheta}$ relative to $\bXi_0$ over $\theta \in [-\pi,\pi].$ In particular, under the non-functional setting with $\bPhi \in \eR^p$ and $f_{\bX,\theta}, \bSigma_0 \in \eR^{p \times p},$ (\ref{eq.bd.fsm}) reduces to
$$
2 \pi \cdot \underset{\theta\in [-\pi, \pi], \bPhi \neq {\bf 0}}{\text{ess}\sup} \frac{\bPhi^{\T} f_{\bX,\theta}\bPhi}{\bPhi^{\T}\bSigma_0 \bPhi} < \infty,
$$
which is equivalent to the upper bound condition for the stability measure, $\widetilde\cM(f_{\bX}),$ introduced by \cite{basu2015a}, i.e.
$$\widetilde\cM(f_{\bX})=\underset{\theta\in [-\pi, \pi], \bPhi \neq {\bf 0}}{\text{ess}\sup} \frac{\bPhi^{\T} f_{\bX,\theta}\bPhi}{\bPhi^{\T}\bPhi} < \infty.$$
Second, 
if $X_{t1}(\cdot), \dots, X_{tp}(\cdot)$ are finite-dimensional objects, the upper bound conditions 
for $\cM(\mF_{\bX})$ and the functional analog of $\widetilde\cM(f_{\bX})$ are equivalent. However, for truly infinite-dimensional functional objects,
$\cM(\mF_{\bX})$ makes more sense, 
since it can more precisely capture the effect of small eigenvalues of $\mF_{\bX,\theta}$'s relative to those of $\bXi_0$. Moreover, Condition~\ref{cond.bd.fsm} is satisfied by a large class of infinite-dimensional functional data, see examples in Appendix~\ref{ap.fsm.ex} for details.
Third, it is clear that, unlike $\widetilde\cM(f_{\bX}),$ $\cM(\mF_{\bX})$ is a scale-free stability measure. In the special case of no temporal dependence, $\cM(\mF_\bX)=1.$
Fourth, since the autocovariance function characterizes a multivariate Gaussian process, it can be used to quantify the temporal and cross-sectional dependence for this class of models. In particular, the spectral density functions provide insights into the stability of the process. In our analysis of high-dimensional functional time series, we will use $\cM(\mF_{\bX})$ as a stability measure of the process of $\{\bX_t(\cdot)\}_{t \in \eZ}.$ Larger values of $\cM(\mF_{\bX})$ would correspond to a less stable process.

We next illustrate the superiority of the functional stability measure to possible competitors using VFAR models as an example. See Section~\ref{sec.vfar} for details on VFAR models.
In particular, we consider a VFAR model of order $1$, denoted by VFAR(1), as follows
\begin{equation}
\label{vfar1.ex}
\bX_t(u) = \int_{\cU} \bA(u,v) \bX_{t-1}(v)dv + \bvar_t(u), ~u\in\cU.
\end{equation}
In the special case of a symmetric $\bA$, i.e. $\bA(u,v)=\bA(v,u)^{\T},$
equation (\ref{vfar1.ex}) has a stationary solution if and only if $||\bA||_{\cL} <1.$ See Theorem~3.5 of \cite{Bbosq1} for $p=1.$ However, this restrictive condition is violated by many stable VFAR(1) models with non-symmetric $\bA.$ Moreover, it does not generalize beyond VFAR(1) models. 

\begin{figure*}[!]
 	\centering
   	\includegraphics[width=6.7cm,height=6.7cm]{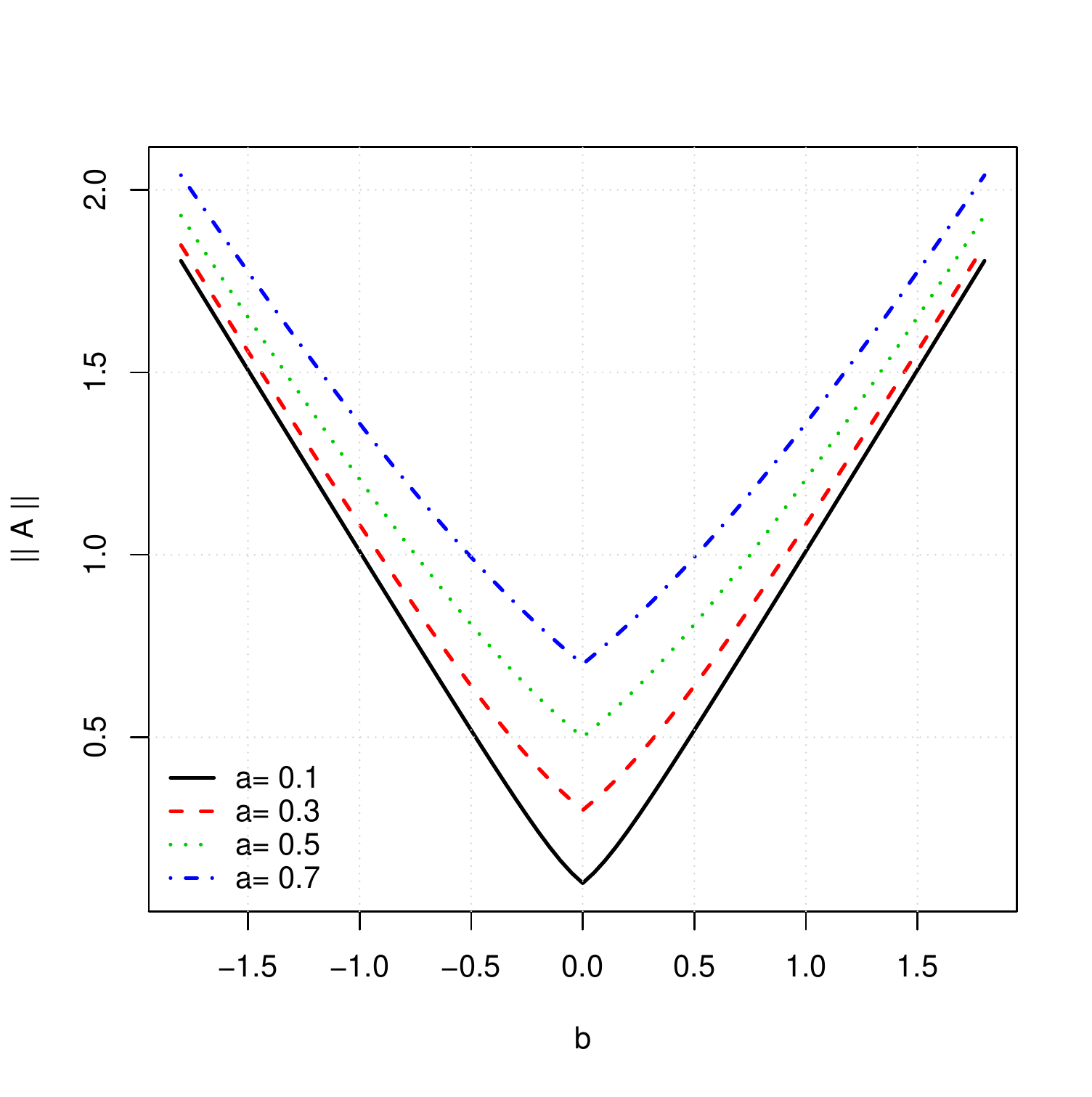}
   	\includegraphics[width=6.7cm,height=6.7cm]{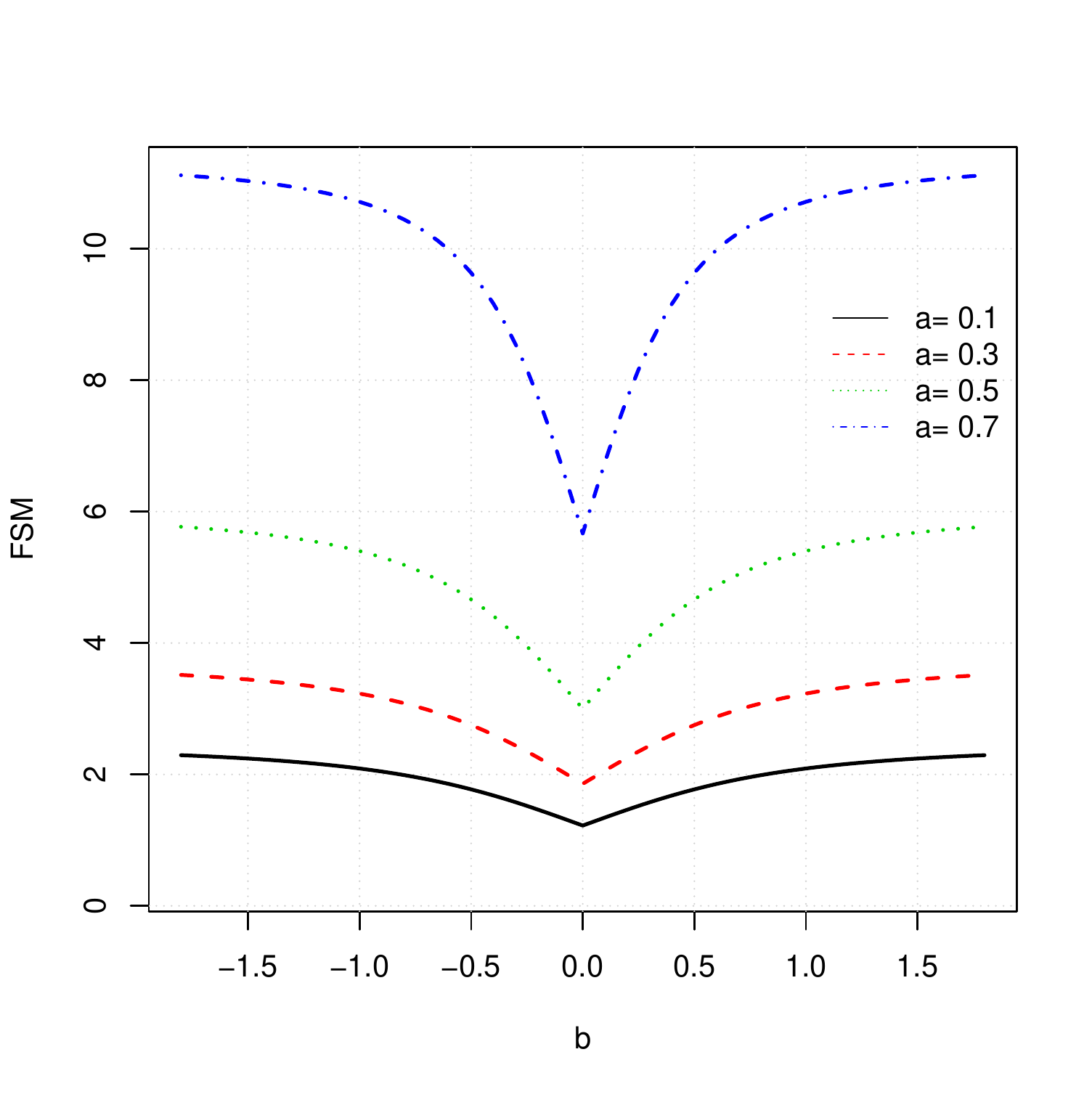}
 	\caption{\label{ill.example}{\it \small The illustrative VFAR(1) model. Left: $||\bA||_{\cL}$ as a function of $a$ and $b$, plotted against $b$ for different $a.$ Right: $\cM(\mF_{\bX})$ as a function of $a$ and $b$, plotted against $b$ for different $a.$}}
 \end{figure*}

We consider an illustrative example with
\begin{equation}
\label{vfar1.ex.A}
\bA(u,v)=\left(\begin{array}{cc}a\psi_1(u)\psi_1 (v)& b\psi_1(u)\psi_2(v)\\0 & a\psi_2(u)\psi_2(v)\end{array}\right),
\bX_t(u)=\left(\begin{array}{c}x_{t1}\psi_1(u)\\x_{t2}\psi_2(u)\end{array}\right),\bvarepsilon_t(u)=\left(\begin{array}{c}e_{t1}\psi_1(u)\\e_{t2}\psi_2(u)\end{array}\right),
\end{equation}
where
$(e_{t1},e_{t2})^{\T} \overset{\text{i.i.d.}}{\sim}N({\bf 0}, \bI_2)$
and $\|\psi_j\|=1$ for $j=1,2.$ Section~\ref{ap.ill.ex} of the Supplementary Material provides details to calculate $\rho(\bA)$ (spectral radius of $\bA$), $||\bA||_{\cL}$ and $\cM(\mF_{\bX})$ for this example. In particular, $\rho(\bA)=|a|<1$ corresponds to a stationary solution to equation~(\ref{vfar1.ex}). Figure~\ref{ill.example}
visualizes $||\bA||_{\cL}$ and $\cM(\mF_{\bX})$ for various values of $a \in (0,1).$ We observe a few apparent patterns. First, increasing $a$ results in a value for larger $||\bA||_{\cL}.$ As $|b|$ grows large enough, the condition of $||\bA||_{\cL}<1$ will be violated, but equation~(\ref{vfar1.ex}) still have a stationary solution. Second, processes with stronger temporal dependence, i.e. with larger values of $a$ or $|b|,$ have larger values of $\cM(\mF_{\bX})$ and will be considered less stable. For a high-dimensional VFAR(1) model, it is more sensible to use $\cM(\mF_{\bX})$ rather than $\|\bA\|_{\cL}$ or $\|\bA^{j}\|_{\cL}$ for some $j \geq 1$ (Theorem~5.1 of \cite[]{Bbosq1}) as a measure of stability of the process.

\begin{definition}
\label{def.sub.fsm}
For all $k$-dimensional subprocesses of $\{\bX_{t}(\cdot)\}_{t\in\eZ},$ i.e. $\big\{\big(X_{tj}(\cdot)\big):j \in J\big\}_{t \in \eZ},$ for $J \subseteq\{1,\dots,p\}$ and $|J| \leq k,$ we define the corresponding functional stability measure by
\begin{equation}
\label{eq.sub.fsm}
\cM_k(\mF_{\bX}) = 2\pi \cdot \underset{\theta\in [-\pi, \pi],\|\bPhi\|_0 \le k,\bPhi \in \cH_0}{\text{ess} \sup} \frac{\langle \bPhi, \mF_{\bX,\theta}(\bPhi)\rangle_{\cH}}{\langle \bPhi, \bXi_0(\bPhi)\rangle_{\cH}}, \ k=1,\dots,p.
\end{equation}
\end{definition}

It is obvious from definitions in  Condition~\ref{cond.bd.fsm} and (\ref{eq.sub.fsm}) that
$
\cM_1(\mF_{\bX}) \le \cM_2(\mF_{\bX}) \le \dots \le \cM_p(\mF_{\bX}) = \cM(\mF_{\bX}) <\infty
$
and $\cM_1(\mF_{\bX})=\max_{1 \leq j \leq p} \cM(\mF_{X_j}),$ which is allowed to evolve with $p$ satisfying $\cM_1(\mF_{\bX}) \lesssim M^{-(2\alpha+1)}\{n/\log(pM)\}^{1/2},$ implied from Theorems~\ref{thm_lam_phi}--\ref{thm_cov} and Condition~\ref{cond_q_n} below.

\subsection{Concentration bounds on $\widehat \bSigma_h$}
\label{sec.bd.Sigma}
Based on $n$ temporally dependent observations 
$\bX_1(\cdot),\dots, \bX_n(\cdot),$ 
we construct an empirical estimator of $\bSigma_{h}$ by
\begin{equation}
    \label{est.autocov.fn}
    \widehat \bSigma_{h}(u,v)=\frac{1}{n-h}\sum_{t=1}^{n-h} \bX_t(u)\bX_{t+h}(v)^{\T}, ~h=0,1, \dots, ~(u,v)\in \cU^2.
\end{equation}

The following theorem provides concentration bounds on $\widehat \bXi_0$ under the quadratic and bilinear forms. These concentration bounds form the core of our theoretical results, which serve as a starting point to establish further non-asymptotic error bounds presented in Sections~\ref{sec.bd.Sigma} and \ref{sec.bd.fpca}.

\begin{theorem}
	\label{thm_con_op}
	Suppose that Conditions~\ref{cond.cov.func} and \ref{cond.bd.fsm} hold. Then for any given vectors $\bPhi_1,$
	$\bPhi_2 \in \cH_0$ satisfying $\|\bPhi_1\|_0 \vee \|\bPhi_2\|_0 \le k$ with some integer $k~(1 \le k \le p),$ there exists some universal constant $c > 0$ such that for any $\eta > 0,$
	\begin{equation}
	\label{thm_Phi_1}
	P\left\{\left|\frac{\big\langle \bPhi_1, (\widehat{\bXi}_0 - \bXi_0)(\bPhi_1)\big\rangle_{\cH}} {\big\langle \bPhi_1, \bXi_0(\bPhi_1)\big\rangle_{\cH}}\right| > \cM_k(\mF_{\bX})\eta\right\} \le 2\exp\Big\{- c n \min(\eta^2, \eta)\Big\},
	\end{equation}
	and
	\begin{equation}
	\label{thm_Phi_12}
	P\left\{\left|\frac{\big\langle \bPhi_1, (\widehat{\bXi}_0 - \bXi_0)(\bPhi_2)\big\rangle_{\cH}} {\big\langle \bPhi_1, \bXi_0(\bPhi_1)\big\rangle_{\cH} +
	\big\langle \bPhi_2, \bXi_0(\bPhi_2)\big\rangle_{\cH}}\right| >   \cM_k(\mF_{\bX}) \eta\right\} \le 4\exp\Big\{- c n \min(\eta^2, \eta)\Big\}.
	\end{equation}	
\end{theorem}

The concentration inequalities in (\ref{thm_Phi_1}) and (\ref{thm_Phi_12}) suggest that the temporal dependence may affect the tail behaviors via $\cM_k(\mF_{\bX})$ in two different ways, depending on which term in the tail bounds is dominant. 
With suitable choices of $\bPhi_1$ and $\bPhi_2,$ we can derive non-asymptotic results for entries of $\widehat\bSigma_{h}$ and relevant estimated terms under a FPCA framework. For example, under the Karhunen-Lo\`eve expansion of each $X_{tj}(\cdot)$ (see details in Section~\ref{sec.bd.fpca}), choosing $\bPhi_1=(0, \dots, 0, \phi_{jl},0, \dots, 0)^{\T}$ leads to $\langle \bPhi_1, \bXi_0(\bPhi_1)\rangle_{\cH}={\lambda_{jl}}$ and $\langle \bPhi_1, \widehat\bXi_0(\bPhi_1)\rangle_{\cH}=n^{-1}\sum_{t=1}^n \langle X_{tj}, \phi_{jl} \rangle^2,$ both of which are useful terms in our further analysis. Moreover, if we choose $
\bPhi_{1} = (0, \ldots,0,\lambda_{jl}^{-1/2}\phi_{jl}, 0,\ldots,0)^{\T} $ and $\bPhi_{2} = (0, \ldots,0,\lambda_{km}^{-1/2}\phi_{km}, 0,\ldots,0)^{\T},$ an application of (\ref{thm_Phi_12}) and some calculations yield elementwise concentration bounds on $\bSigma_0$ as stated in the following theorem.


\begin{theorem}
	\label{thm_max_bound}
	Suppose that Conditions~\ref{cond.cov.func} and \ref{cond.bd.fsm} hold. Then there exists some universal constant $\tilde c > 0$ such that for any $\eta > 0$ and each $j,k =1,\ldots,p,$
	\begin{equation}
	\label{thm_Sigma_compt}
	P\left\{\big\|\widehat{\Sigma}_{jk}^{(0)} - \Sigma_{jk}^{(0)}\big\|_{\cS} >   2 \cM_1(\mF_{\bX})\lambda_0 \eta\right\} \le 4 \exp\Big\{- \tilde c n \min(\eta^2, \eta)\Big\},
	\end{equation}
	and
	\begin{equation}
	\label{thm_Sigma_max}
	P\left\{\big\|\widehat{\bSigma}_0 - \bSigma_0\big\|_{\max} >   2 \cM_1(\mF_{\bX})\lambda_0 \eta\right\} \le 4p^2 \exp\Big\{- \tilde c n \min(\eta^2, \eta)\Big\}.
	\end{equation}
	In particular, if $n \geq \rho^2\log p,$ where $\rho$ is some constant with $\rho > \sqrt{2} \tilde c^{-1/2},$ then with probability greater than $1- 4p^{2- \tilde c\rho^2}$, the estimate $\widehat{\bSigma}_0$ satisfies 
	\begin{equation}
	\label{bd_Sigma_max}
	\big\|\widehat{\bSigma}_0 - \bSigma_0\big\|_{\max} \le 2 \cM_1(\mF_{\bX})\lambda_0 \rho \sqrt{\frac{\log p}{n}}.
	\end{equation}
\end{theorem}

Under an independent setting, the diagonalwise concentration properties of $\widehat \bSigma_{0}$ were studied in \cite{koltch2017} and \cite{qiao2018a}. 
In particular, \cite{koltch2017} established the concentration bound under operator norm, i.e., with probability greater than $1-e^{-\delta}$ for all $\delta \geq 1,$ 
\begin{equation}
\label{conbd_op}
\|\widehat\Sigma_{jj}^{(0)}-\Sigma_{jj}^{(0)}\|_{\cL} \lesssim \|\Sigma_{jj}^{(0)}\|_{\cL}\left\{\sqrt{\frac{\tilde r\big(\Sigma_{jj}^{(0)}\big)}{n}} \vee \frac{\tilde r\big(\Sigma_{jj}^{(0)}\big)}{n} \vee \sqrt{\frac{\delta}{n}} \vee \frac{\delta}{n} \right\},
\end{equation}
where
$\tilde r(\Sigma_{jj}^{(0)}) \asymp \int_{\cU}\Sigma_{jj}^{(0)}(u,u)du/\|\Sigma_{jj}^{(0)}\|_{\cL}.$ 
By contrast, when $\bX_1(\cdot), \dots, \bX_n(\cdot)$ are temporally dependent, (\ref{thm_Sigma_compt}) implies that, with probability greater than $1-4e^{-\delta},$
\begin{equation}
\label{conbd_hs}
\|\widehat\Sigma_{jk}^{(0)}-\Sigma_{jk}^{(0)}\|_{\cS} \lesssim \cM_1(f_{\bX})\max_j\int_{\cU}\Sigma_{jj}^{(0)}(u,u)du \left( \sqrt{\frac{\delta}{n}} \vee \frac{\delta}{n}\right).
\end{equation}
Note (\ref{conbd_op}) and (\ref{conbd_hs}) are both essential concentration bounds of independent interests and share the the common term, $\int_{\cal U}\Sigma_{jj}^{(0)}(u,u)du,$ and the same rate via $(n^{-1}\delta)^{1/2} \vee n^{-1}\delta$ but with different multiplicative terms. However, concentration bounds in (\ref{conbd_hs}) under Hilbert--Schmidt norm with free choices of $(j,k)$ play a crutial role in deriving the convergence results for PFCA via Theorems~\ref{thm_lam_phi} and \ref{thm_cov} below in a high-dimensional and dependent setting.
In the bounds established in Theorem~\ref{thm_max_bound}--\ref{thm_cov}, 
the effects of temporal dependence are commonly captured by $\cM_1(\mF_{\bX})$ with larger values yielding a slower convergence rate. 

We next present similar concentration bounds on $\widehat{\bXi}_h$'s for $h > 0.$

\begin{proposition}
	\label{thm_sigma_k}
	Suppose that Conditions~\ref{cond.cov.func}--\ref{cond.bd.fsm} hold and $h$ is fixed. Then for any given vectors $\bPhi_1,$
	$\bPhi_2 \in \cH_0$ satisfying $\|\bPhi_1\|_0 \vee \|\bPhi_2\|_0 \le k$ with some integer $k~(1 \le k \le p),$  there exists some universal constant $c > 0$ such that for any $\eta > 0,$
	\begin{equation}
	\label{thm_sig_k1}
	P\left\{\left|\frac{\big\langle \bPhi_1, (\widehat{\bXi}_h - \bXi_h)(\bPhi_1)\big\rangle_{\cH}} {\big\langle \bPhi_1, \bXi_0(\bPhi_1)\big\rangle_{\cH}}\right| > 2 \cM_k(\mF_{\bX})\eta\right\} \le 4\exp\Big\{- c n \min(\eta^2, \eta)\Big\},
	\end{equation}
	and
	\begin{equation}
	\label{thm_sig_k2}
	P\left\{\left|\frac{\big\langle \bPhi_1, (\widehat{\bXi}_h - \bXi_h)(\bPhi_2)\big\rangle_{\cH}}
    	{\big\langle \bPhi_1, \bXi_0(\bPhi_1)\big\rangle_{\cH} +
    	\big\langle \bPhi_2, \bXi_0(\bPhi_2)\big\rangle_{\cH}}\right| >    2\cM_k(\mF_{\bX}) \eta\right\} \le 8\exp\Big\{- c n \min(\eta^2, \eta)\Big\}.
	\end{equation}	
\end{proposition}

With the same choices of $\bPhi_1$ and $\bPhi_2$ as those used in applying Theorem~\ref{thm_con_op} to prove Theorem~\ref{thm_max_bound}, the concentration bounds in (\ref{thm_sig_k1}) and (\ref{thm_sig_k2}) can lead to 
$\big\|\widehat{\bSigma}_h - \bSigma_h\big\|_{\max}=O_P\left\{\cM_1(\mF_{\bX})(\log p/n)^{1/2}\right\}.$ Moreover, these concentration results are useful to address other important high-dimensional functional time series problems, e.g., high-dimensional functional factor models and non-asymptotic analysis of dynamic FPCA, as discussed in Section~\ref{sec.discussion}.



\subsection{Rates in elementwise $\ell_{\infty}$-norm under a FPCA framework}
\label{sec.bd.fpca}
For each $j=1,\dots, p,$ we assume that $X_{tj}(\cdot)$ admits the Karhunen-Lo\`eve expansion, i.e.\\
$
    X_{tj}(\cdot)=\sum_{l=1}^{\infty}\xi_{tjl}\phi_{jl}(\cdot),
$
which forms the foundation of FPCA.
The coefficients $\xi_{tjl}=\langle X_{tj}, \phi_{jl}\rangle,$ $l\geq 1,$ namely {\it functional principal component} (FPC) scores, correspond to a sequence of random variables with $E(\xi_{tjl})=0,$ $\text{Var}(\xi_{tjl})=\lambda_{jl}$ and $\cov(\xi_{tjl},\xi_{tjl'})=0$ if $l \neq l'.$ The eigenpairs $\{(\lambda_{jl},\phi_{jl})\}_{l \geq 1}$ satisfy the eigen-decomposition $ \langle\Sigma_{jj}^{(0)}(u,\cdot), \phi_{jl}(\cdot) \rangle=\lambda_{jl}\phi_{jl}(u)$ with $\lambda_{j1} \geq \lambda_{j2}\geq \cdots.$ We say that $X_{tj}(\cdot)$ is $d_j$-dimensional if $\lambda_{jd_{j}} \neq 0$ and $\lambda_{j(d_j+1)}=0$ for some positive integer $d_j.$ If $d_j=\infty,$ all the eigenvalues are nonzero and $X_{tj}(\cdot)$ is a truly infinite-dimensional functional object.

To implement FPCA based on realizations $\{X_{1j}(\cdot), \dots, X_{nj}(\cdot)\},$ we first compute the sample estimator of $\Sigma_{jj}^{(0)}$ by $\widehat\Sigma_{jj}^{(0)}(u,v)=n^{-1}\sum_{t=1}^nX_{tj}(u)X_{tj}(v).$ Performing an eigenanalysis on $\widehat\Sigma_{jj}^{(0)},$ i.e. $\langle \widehat\Sigma_{jj}^{(0)}(u,\cdot), \widehat\phi_{jl}(\cdot)\rangle=\widehat \lambda_{jl}\widehat\phi_{jl}(u)$ for $l\geq 1,$ leads to estimated eigenpairs $(\widehat\lambda_{jl},\widehat\phi_{jl})$ and estimated FPC scores $\widehat\xi_{tjl}=\langle X_{tj}, \widehat\phi_{jl}\rangle.$ In the following, we will provide the 
non-asymptotic analysis of estimated eigenpairs. 

\subsubsection{Eigenvalues and eigenfunctions}
\label{sec.bd.eigen}
We first impose the following regularity condition.

\begin{condition}
\label{cond_eigen}
For each $j=1,\dots,p,$ all the nonzero eigenvalues of $\Sigma_{jj}^{(0)}$ are different, i.e. $\lambda_{j1} >  \lambda_{j2} > \cdots > 0,$ and there exist some positive constants $c_0$ and $\alpha>1$ such that $ \lambda_{jl}-\lambda_{j(l+1)} \geq c_0 l^{-\alpha-1}$ for $l =1,\dots,\infty.$
\end{condition}

Condition~\ref{cond_eigen} is standard in functional data analysis literature, see, e.g., \cite{hall2007} and \cite{kong2016}. 
The parameter $\alpha$ controls lower bounds for spacings between adjacent eigenvalues with larger values of $\alpha$ allowing tighter eigengaps. This condition also implies that $\lambda_{jl} \geq c_0 \alpha^{-1} l^{-\alpha}$ as $\lambda_{jl}=\sum_{k=l}^{\infty}\{\lambda_{jk}-\lambda_{j(k+1)}\}\geq c_0\sum_{k=l}^\infty k^{-\alpha-1}.$ 


In the following theorem, we present relative error bounds on $\{\widehat\lambda_{jl}\}$ and $\|\widehat\phi_{jl}-\phi_{jl}\|$ in elementwise $\ell_{\infty}$ norm, which plays an crucial rule for the further consistency analysis under high-dimensional scaling.
\begin{theorem}
\label{thm_lam_phi}
Suppose that Conditions~\ref{cond.cov.func}--\ref{cond_eigen} hold. 
Let $M$ be a positive integer possibly depending on $(n,p).$ If
$n \gtrsim M^{4\alpha+2}\cM_{1}^2(\mF_{\bX}) \log(pM),$
then there exist some positive constants $c_1$ and $c_2$ independent of $(n,p,M)$ such that, with probability greater than 
$ 1- c_1(pM)^{-c_2}$, the estimates $\{\widehat\lambda_{jl}\}$ and $\{\widehat\phi_{jl}(\cdot)\}$ satisfy
\begin{equation}
\label{bd_max_lam_phi}
\max_{1 \le j\le p,1 \le l \le M} \left\{\Big|\frac{\widehat \lambda_{jl} - \lambda_{jl}}{\lambda_{jl}}\Big| 
+ \Big\|\frac{ \widehat \phi_{jl} - \phi_{jl}}{l^{\alpha + 1}}\Big\| \right\}
\lesssim \cM_{1}(\mF_{\bX})\sqrt{\frac{\log (pM)}{n}}.
\end{equation}
\end{theorem}

We provide three remarks for the relative errors of $\{\widehat\lambda_{jl}\}$ and $\|\widehat \phi_{jl}-\phi_{jl}\|.$ First, compared with the non-asymptotic results for the absolute errors of $\{\widehat\lambda_{jl}\}$ under an independent setting \cite[]{qiao2018a}, Theorem~\ref{thm_lam_phi} does not require the upper bound condition for eigenvalues. Moreover, when $\cM_{1}(\mF_{\bX})$ remains constant with respect to $p,$ provided that $\lambda_{jl}$ converges to zero as $l$ grows to infinity, (\ref{bd_max_lam_phi}) leads to a faster rate of convergence for small eigenvalues, i.e. $|\widehat\lambda_{jl}-\lambda_{jl}|=O_P(\lambda_{jl}n^{-1/2}.$ Such rate is also sharp in the sense of \cite{jirak2016}. See also
\cite{av2020}, which, under an independent setting, established relative concentration bounds on $\{\widehat\lambda_{jl}\}$ similar to Lemma~\ref{lemma_lambda} in the Supplementary Material.
Second, when each $X_{tj}(\cdot)$ is finite-dimensional with $\lambda_{jl}=0$ for $l>d_j,$  
the estimators of zero-eigenvalues enjoy the faster rates due to the property of first order degeneracy \cite[]{bathia2010}. 
Third, error bounds on $\widehat\phi_{jl}$ under $\ell_2$ norm 
are derived using a well-known pathway bound in Lemma~4.3 of \cite{Bbosq1}. For finite-dimensional functional objects, such error bounds lead to the optimal $\sqrt{n}$-rate. 
Under an infinite-dimensional setting, our derived error bounds on $\widehat \phi_{jl}(\cdot)$'s can possibly be improved from a probabilistic perspective as long as $l$ diverges. Interestingly, \cite{jirak2016} established a sharper rate by $\Lambda_{jl}^{-1/2}\|\widehat\phi_{jl}-\phi_{jl}\| =O_P(n^{-1/2})$ with $\Lambda_{jl}=\sum_{k \neq l}^{\infty}\frac{\lambda_{jl}\lambda_{jk}}{(\lambda_{jl}-\lambda_{jk})^2}.$ In particular, when $\lambda_{jl} \asymp l^{-\alpha}$ with $\alpha > 1$, $\Lambda_{jl} \lesssim l^{2}$ as $l \to \infty.$ We believe that, with the help of techniques used in \cite{jirak2016}, the relevant bounds in Theorems~\ref{thm_lam_phi} and \ref{thm_cov} can be further sharpened. The extra complication, however, will make the theoretical justification of the VFAR estimate in Section~\ref{sec.vfar} more challenging. We leave the development of optimal results as a topic for future research.

Finally, we give two remarks on the parameter $M$. First, $M$ can be viewed as the truncated dimension of $X_{tj}(\cdot) \approx \sum_{l=1}^{M} \xi_{tjl}\phi_{jl}(\cdot)$.
In general, $M$ can depend on $j,$ say $M_j,$ then the right-side of  (\ref{bd_max_lam_phi}) becomes ${\cal M}_1(f_{\bX})\{\log(\sum_{j=1}^p M_j)/n\}^{1/2}.$
Second, the sample size lower bound in Theorem~\ref{thm_lam_phi} implies that $M\lesssim \big[n/\{\cM_{1}^2(f_{\bX})\log p\}\big]^{1/(4\alpha+2)}$, which controls the rate that $M$ can grow at most as a function of $n, p, \cM_{1}(f_{\bX})$ and $\alpha,$ with larger values of $n$ or smaller values of $p$ or $\alpha$ allowing a larger $M.$


\subsubsection{Covariance between FPC scores}
\label{sec.bd.score}
For each $j,k=1,\dots,p,$ $l,m=1, 2, \dots,$ and $h=0,1,\dots,$ let $\sigma_{jklm}^{(h)}=E(\xi_{tjl} \xi_{(t+h)km})$  and its sample estimator be $\widehat \sigma_{jklm}^{(h)}=(n-h)^{-1}\sum_{t=1}^{n-h}\widehat\xi_{tjl}\widehat\xi_{(t+h)km}.$ In the second step of the three-step procedure, our main target is to fit a sparse model based on temporally dependent estimated FPC scores, $\{\widehat\xi_{tjl}\}.$ To provide the theoretical guarantee for this step, 
we present the convergence analysis of $\{\widehat \sigma_{jklm}^{(h)}\}$ in elementwise $\ell_{\infty}$ norm as follows.

\begin{theorem}
\label{thm_cov}
	Suppose that Conditions~\ref{cond.cov.func}--\ref{cond_eigen} hold and $h$ is fixed.
	Let $M$ be a positive integer possibly depending on $(n,p).$ If
	$n \gtrsim M^{4\alpha+2}\cM_{1}^2(\mF_{\bX}) \log(pM),$
then there exist some positive constants $c_3$ and $c_4$ independent of $(n,p,M)$ such that, with probability greater than 
$ 1- c_3(pM)^{-c_4}$, the estimates $\big\{\widehat \sigma_{jklm}^{(h)}\big\}$ satisfy
\begin{equation}
\label{bd_score_max}
\underset{\underset{1 \le l,m \le M}{1\le j,k\le p}}{\max} \frac{\left|\widehat \sigma_{jklm}^{(h)} - \sigma_{jklm}^{(h)}\right|}{(l \vee m)^{\alpha+1}\lambda_{jl}^{1/2} \lambda_{km}^{1/2}}\lesssim \cM_1(\mF_{\bX}) \sqrt{\frac{\log (pM)}{n}}.
\end{equation}
\end{theorem}

We provide three comments here. First, compared with the convergence rate of absolute errors of $\{\widehat\sigma_{jklm}^{(0)}\}$ under an independent setting, i.e. $|\widehat\sigma_{jklm}^{(0)}-\sigma_{jklm}^{(0)}|=O_P\big\{(l+m)^{\alpha+1}n^{-1/2}\big\}$ \cite[]{qiao2018a}, we obtain that of scaled errors of $\{\widehat\sigma_{jklm}^{(h)}\}$ as $\lambda_{jl}^{-1/2}\lambda_{km}^{-1/2}|\widehat\sigma_{jklm}^{(h)}-\sigma_{jklm}^{(h)}|=O_P\big\{\cM_1(\mF_{\bX})(l\vee m)^{\alpha+1}n^{-1/2}\big\},$ which can more precisely characterize the effect of small eigenvalues on the convergence. 
Second, in the special case of $j=k$ and $l=m$ with $\sigma_{jjll}^{(0)}=\lambda_{jl}$ and $\widehat\sigma_{jjll}^{(0)}=\widehat\lambda_{jl},$ the scaled errors of $\{\widehat\sigma_{jjll}^{(0)}\},$ i.e. $\lambda_{jl}^{-1}|\widehat \lambda_{jl}-\lambda_{jl}|,$ would correspond to a faster convergence rate due to (\ref{bd_max_lam_phi}). Third, if we relax Condition~\ref{cond_eigen} by allowing the parameter $\alpha$ to depend on $j,$ the resulting scaled error rate becomes
$\lambda_{jl}^{-1/2}\lambda_{km}^{-1/2}|\widehat\sigma_{jklm}^{(h)}-\sigma_{jklm}^{(h)}|=O_P\big\{\cM_1(\mF_{\bX})(l^{\alpha_j+1}\vee m^{\alpha_k+1})n^{-1/2}\big\}.$

\section{Vector functional autoregressive models}
\label{sec.vfar}
Inspired from the standard VAR formulation, we propose a VFAR model of lag $L$, namely VFAR($L$), which is able to characterize linear inter-dependencies in $\{\bX_t(\cdot)\}_{t \in \eZ}$
as follows
\begin{equation}
\label{vfar1}
\bX_{t}(u) = \sum_{h = 1}^L \int_{\cU}\bA_{h}(u,v) \bX_{t-h}(v)dv  + \bvarepsilon_{t}(u), ~~t =L+1, \ldots,n,
\end{equation}
where $\bvar_t(\cdot)=\big(\varepsilon_{t1}(\cdot), \dots, \varepsilon_{tp}(\cdot)\big)^{\T}$ are
independently sampled from a $p$-dimensional vector of mean zero Gaussian processes, independent of $\bX_{t-1}(\cdot), \bX_{t-2}(\cdot), \dots,$ and 
$\bA_h=\big(A_{jk}^{(h)}\big)_{1\leq j,k\leq p}$ is the transition function at lag $h$ with $A_{jk}^{(h)} \in \eS.$
The structure of transition functions provides insights into the temporal and cross-sectional inter-relationship  amongst $p$ functional time series. To make a feasible fit to (\ref{vfar1}) in a high-dimensional regime, we assume the functional sparsity in $\bA_1, \dots, \bA_L,$ i.e. most of the components in $\big\{X_{(t-h)k}(\cdot): h=1,\dots,L, k=1,\dots,p\big\}$ are unrelated to $X_{tj}(\cdot)$ for $j=1,\dots,p.$

Due to the infinite-dimensional nature of functional data, for each $j,$ we take a standard dimension reduction approach through FPCA to approximate $X_{tj}(\cdot)$ using the leading $q_j$ principal components, i.e. $X_{tj}(\cdot)\approx \sum_{l=1}^{q_j}\xi_{tjl}\phi_{jl}(\cdot)=\bxi_{tj}^{\T}\bphi_{j}(\cdot),$ where $\bxi_{tj}=(\xi_{tj1}, \dots, \xi_{tjq_j})^{\T}$, $\bphi_{j}(\cdot)=\big(\phi_{j1}(\cdot), \dots, \phi_{jq_j}(\cdot)\big)^{\T}$ and $q_j$ is chosen data-adaptively to provide a reasonable approximation to the trajectory $X_{tj}(\cdot).$

Once the FPCA has been performed for each $X_{tj}(\cdot),$ we let $\bV_{j}^{(h)} \in {\eR}^{(n-L)\times q_j}$ with its row vectors given by $\bxi_{(L+1-h)j}, \dots, \bxi_{(n-h)j}$ and
$\bPsi_{jk}^{(h)}=\int_{\cU}\int_{\cU}\bphi_k(v)A_{jk}^{(h)}(u,v)\bphi_j(u)^{\T}dudv \in {\eR}^{q_k \times q_j}.$ Then further derivations in Section~\ref{ap.vfar.deriv} of the Supplementary Material lead to the matrix representation of (\ref{vfar1}) as
\begin{equation}
\label{vfar3}
\bV_{j}^{(0)} = \sum_{h = 1}^L \sum_{k=1}^p
\bV_{k}^{(h)}\bPsi_{jk}^{(h)} + \bR_j + \bE_j, ~~j=1,\dots,p,
\end{equation}
where $\bR_{j}$ and $\bE_j$ are $(n-L) \times q_j$ error matrices whose row vectors are formed by the truncation and random errors, respectively. Hence, we can rely on the block sparsity pattern in $\{\bPsi_{jk}^{(h)}:h=1,\dots,L,j,k=1,\dots,p\}$ to recover the functional sparsity structure in $\{A_{jk}^{(h)}:h=1,\dots,L,j,k=1,\dots,p\}.$ It is also worth noting that (\ref{vfar3}) can be viewed as a $\big(\sum_{j=1}^p q_j\big)$-dimensional VAR(L) model with the error vector consisting of both the truncation and random errors.

\subsection{Estimation procedure}
\label{sec.vfar.est}
The estimation procedure proceeds in the following three steps.

{\bf Step~1}. We perform FPCA based on observed curves, $X_{1j}(\cdot), \dots, X_{nj}(\cdot)$ and thus obtain estimated eigenfunctions $\widehat\bphi_{j}(\cdot) = \big(\widehat\phi_{jl}(\cdot), \dots, \widehat\phi_{jq_j}(\cdot)\big)^{\T}$ and FPC scores $\widehat\bxi_{tj} = \big(\widehat\xi_{tj1},\dots, \widehat\xi_{tjq_j}\big)^{\T}$ for each $j.$ See Section~\ref{ap.select.tune} of the Supplementary Material for the selection of $q_j$'s in practice.

{\bf Step 2}. Motivated from the matrix representation of a VFAR(L) model in (\ref{vfar3}), we propose a penalized {\it least squares} (LS) approach, which minimizes the following optimization criterion over $\{\bPsi_{jk}^{(h)}: h=1,\dots, L, k=1,\dots,p\}:$
\begin{equation}
    \label{vfar.crit}
    \frac{1}{2}\left\|\widehat\bV_{j}^{(0)}-\sum_{h=1}^L \sum_{k=1}^p \widehat\bV_{k}^{(h)}\bPsi_{jk}^{(h)}\right\|_F^2 + \gamma_{nj} \sum_{h=1}^L \sum_{k=1}^p  \left\|\widehat\bV_{k}^{(h)}\bPsi_{jk}^{(h)}\right\|_F,
\end{equation}
where $\widehat\bV_{j}^{(h)},$ the estimate of $\bV_{j}^{(h)},$ is a $(n-L) \times q_j$ matrix with its $i$-th row vector given by $\widehat\bxi_{(L+i-h)j}$ for $i=1, \dots, (n-L),$ and $\gamma_{nj} \geq 0$ is a regularization parameter. The $\ell_1/\ell_2$ type of standardized group lasso penalty \cite[]{simon2012} in (\ref{vfar.crit}) forces the elements of $\bPsi_{jk}^{(h)}$ to either all be zero or non-zero. Potentially, one could modify (\ref{vfar.crit}) by adding an unstandardized group lasso penalty
\cite[]{yuan2006} in the form of $\gamma_{nj} \sum_{h=1}^L \sum_{k=1}^p  \left\|\bPsi_{jk}^{(h)}\right\|_F$ to produce the block sparsity in $\{\bPsi_{jk}^{(h)}\}.$ However, orthonormalization within each group would correspond to the uniformly
most powerful invariant test for inclusion of a group, hence we use a standardized group lasso penalty here. 
In Section~\ref{ap.alg.vfar} of the Supplementary Material, we develop a block version of {fast iterative shrinkage-thresholding algorithm} (FISTA), which mirrors recent gradient-based techniques \cite[]{beck2009, dono2015}, to solve the optimization problem in (\ref{vfar.crit}) with the solution given by $\{\widehat \bPsi_{jk}^{(h)}\}.$ The proposed block FISTA algorithm is easy to implement and converges very fast, thus is suitable for solving large-scale optimization problems.

{\bf Step~3}. Finally, we recover the functional sparse estimates of elements in $\{A_{jk}^{(h)}\}$ by the block sparse estimates in $\{\widehat \bPsi_{jk}^{(h)}\}$ via 
\begin{equation}
\label{est.A}
\widehat A_{jk}^{(h)}(u,v)=\widehat\bphi_k(v)^{\T}\widehat\bPsi_{jk}^{(h)}\widehat\bphi_j(u), ~~h=1,\dots,L, ~ j,k=1,\dots,p.
\end{equation}

\subsection{Functional network Granger causality}
\label{sec.ngc}
In this section, we extend the definition of {network Granger causality} (NGC) under a VAR framework \cite[]{Bluk2005} to the functional domain and then use the extended definition under our proposed VFAR framework to understand the causal relationship in $\{\bX_t(\cdot)\}_{t \in \eZ}$

In an analogy to the NGC formulation, a {functional NGC} (FNGC) model consists of $p$ nodes, one for each functional variable, and a number of edges with directions connecting a subset of nodes. Specifically, functional times series of $\{X_{tk}(\cdot)\}_{t\in\eZ}$ is defined to be Granger causal for that of $\{X_{tj}(\cdot)\}_{t\in\eZ}$ or equivalently there is an edge from node $k$ to node $j,$ if $A_{jk}^{(h)}(u,v) \neq 0$ for some $(u,v)\in \cU^2$ or $h\in \{1,\dots, L\}.$ Then our proposed FNGC model can be represented by a directed graph $G=(V,E)$ with vertex set $V=\{1, \dots, p\}$ and edge set
$$E=\left\{(k,j): A_{jk}^{(h)}(u,v) \neq 0 \text{ for some } (u,v) \in \cU^2 \text{ or } h \in \{1,\dots, L\}, (j,k) \in V^2\right\}.$$ 
It is worth noting that, at lag $h,$ $\|A_{jk}^{(h)}\|_{\cS}$ can be viewed as explaining the global Granger-type casual impact of $X_k(\cdot)$ on $X_j(\cdot),$ while $A_{jk}^{(h)}(u,v)$ itself accounts for the local Granger-type casual impact of $X_k(v)$ on $X_j(u).$ To explore the FNGC structure and the direction of influence from one node to the other, we need to develop an approach to estimate $E,$ i.e. identifying the locations of non-zero entries in $\widehat \bA_1,\dots, \widehat\bA_L$ under the Hilbert-Schmidt norm, the details of which are presented in Section~\ref{sec.vfar.est}. 

\subsection{Theoretical properties}
\label{sec.vfar.theory}
According to Section~\ref{ap.vfar1.rep} of the Supplementary Material, all VFAR($L$) models in (\ref{vfar1}) can be reformulated as a VFAR(1) model. Without loss of generality, we consider a VFAR(1) model in the form of 
$$
\bX_t(u)=\int_{\cU}\bA(u,v)\bX_{t-1}(v)dv + \bvarepsilon_t(u), ~~ t=2, \dots,n, ~u \in \cU.
$$

To simplify our notation in this section, we focus on the setting where $q_j$'s are the same across $j=1,\dots,p.$ However, our theoretical results extend naturally to the more general setting. In our empirical studies, we select different $q_j$'s, see Section~\ref{ap.select.tune} of the Supplementary Material for details. Let $\widehat\bZ=(\widehat\bV_{1}^{(1)}, \dots, \widehat\bV_{p}^{(1)})^{
\T} \in {\eR}^{(n-1)\times pq},$ $
\bPsi_{j} = \big((\bpsi_{j1}^{(1)})^\T, \dots, (\bpsi_{jp}^{(1)})^\T\big)^\T \in \eR^{pq \times q}.
$ and
$\widehat\bD=\text{diag}\big(\widehat\bD_1,\dots, \widehat\bD_p\big) \in {\eR}^{pq\times pq},$ where $\widehat \bD_k=\big\{(n-1)^{-1}(\widehat\bV_{k}^{(1)})^{\T}\widehat\bV_{k}^{(1)}\big\}^{1/2} \in {\eR}^{q \times q}$ for $k=1,\dots,p.$ Then minimizing (\ref{vfar.crit}) over $\bPsi_j\in {\eR}^{pq\times q}$ is equivalent to minimizing the following criterion over $\bB_{j} \in {\eR}^{pq\times q},$
\begin{eqnarray}
\label{vfar.crit.theory}
-\langle\langle \widehat \bY_j, \bB_{j}\rangle\rangle +
\frac{1}{2} \langle\langle \bB_{j}, \widehat \bGamma\bB_{j}\rangle\rangle  + \gamma_{nj}  \|\bB_{j}\|_{1}^{(q)},
\end{eqnarray}
where $\widehat \bY_j = (n-1)^{-1}\widehat\bD^{-1} \widehat{\bZ}^\T\widehat{\bV}_{j}^{(0)},$  $\widehat \bGamma= (n-1)^{-1}\widehat\bD^{-1}\widehat{\bZ}^\T\widehat{\bZ} \widehat\bD^{-1}.$ Let $\widehat\bB_{j}$ be the minimizer of (\ref{vfar.crit.theory}), then $\widehat\bPsi_{j}=\widehat\bD^{-1}\widehat\bB_{j}$ with its $k$-th row block given by $\widehat\bPsi_{jk}$ and $\widehat\bA =(\widehat A_{jk})$ with its $(j,k)$-th entry, $\widehat A_{jk}(u,v)=\widehat \bphi_{k}(v)^{\T}\widehat \bPsi_{jk}\widehat\bphi_{j}(u)$ for $j,k=1,\dots,p$ and $(u,v)\in \cU^2.$

Before imposing the condition on the entries of $\bA=(A_{jk}),$ we begin with some notation. For the $j$-th row of $\bA,$ we denote the set of non-zero functions by $S_j=\big\{k \in \{1,\dots,p\}: \|A_{jk}\|_{\cS} \neq 0 \big\}$ and its cardinality by $s_j=|S_j|$ for $j=1,\dots,p.$ We also denote the maximum degree or row-wise cardinality by $s=\max_j s_j$ (possibly depends on $n$ and $p$), corresponding to the maximum number of non-zero functions in any row of $\bA.$ For a sparse VFAR($L$) model in (\ref{vfar1}) with $s_{jh}=\sum_{k=1}^p I(\|A^{(h)}_{jk}\|_{\cal S} \neq 0)$ for $j=1, \dots, p$ and $h=1, \dots, L,$ according to Appendix~\ref{ap.vfar1.rep}, the equivalent sparse VFAR(1) model in (\ref{vfar.eqv}) has the maximum degree $s=\max_j (\sum_{h} s_{jh}) \lesssim \max_{j,h} s_{jh}$ assuming a fixed lag order $L.$ 

\begin{condition}
\label{cond.fvar.bias}	
For each $j=1,\dots,p$ and $k\in S_j,$ $\tA_{jk}(u,v) = \sum_{l,m=1}^\infty a_{jklm} \phi_{jl}(u) \phi_{km}(v)
$
and there exist some positive constants $\beta > \alpha/2+1$ and $\mu_{jk}$ such that $|a_{jklm}| \le \mu_{jk}(l+m)^{- \beta - 1/2}$ for $l,m \geq 1.$ 
\end{condition}

For each $(j,k)$, the basis with respect to which coefficients $\{a_{jklm}\}_{l,m \geq 1}$ are defined is determined by $\{\phi_{jl}(\cdot)\}_{l\geq 1}$ and $\{\phi_{km}(\cdot)\}_{m\geq 1}$ 
The parameter $\beta$ in Condition~\ref{cond.fvar.bias} determines the decay rate of the upper bounds for coefficients $\{a_{jklm}\}_{l,m \geq 1}$ and hence characterizes the degree of smoothness in  $\{A_{jk}\},$ with larger values of $\beta$ yielding smoother functions.
See also \cite{hall2007} and 
\cite{kong2016} for similar smoothness conditions in functional linear models.

We next establish the consistency of the VFAR estimate based on the sufficient conditions in Conditions~\ref{cond.fvar.RE}--\ref{cond.fvar.max.error} below. 
To be specific, we first establish an upper bound on $\|\widehat{\bA} - \bA\|_\infty$ in Theorem~\ref{thm.vfar} below. 
Using the convergence results in Section~\ref{sec.bd.fpca}, we then show that all VFAR models in (\ref{vfar1}) satisfy Conditions~\ref{cond.fvar.RE}--\ref{cond.fvar.max.error} with high probability through Propositions~\ref{res.re}--\ref{prop_Error_max} below. As a consequence, the error bound in Theorem~\ref{thm.vfar} holds with high probability.

Before stating these conditions, we give some notation. For a block matrix $\bB = (\bB_{jk}) \in \eR^{p_1q \times p_2 q}$ with its $(j,k)$-th block $\bB_{jk} \in {\eR}^{q \times q},$  we define its $q$-block versions of Frobenius norm, elementwise $\ell_{\infty}$ norm and matrix $\ell_1$ norm by $\|\bB\|_{\tF} = \big(\sum_{j,k}\|\bB_{jk}\|_{\tF}^2\big)^{1/2},$  $\|\bB\|_{\max}^{(q)} = \max_{j,k} \|\bB_{jk}\|_{\tF}$ and
$
\|\bB\|_{1}^{(q)} = {\max}_{k}\sum_{j} \|\bB_{jk}\|_{\tF},
$ respectively.

\begin{condition}
\label{cond.fvar.RE}
The symmetric matrix $\widehat \bGamma \in {\eR}^{pq \times pq}$ satisfies the restricted eigenvalue condition with tolerance $\tau_1>0$ and curvature $\tau_2>0$ if
\begin{equation}
\label{cond_RE}
\btheta^\T \widehat \bGamma \btheta \geq  \tau_2 \|\btheta\|^2 -\tau_1 \|\btheta\|_1^2 \quad \forall \btheta \in {\eR}^{pq}.
\end{equation}
\end{condition}

Condition~\ref{cond.fvar.RE} comes from a class of conditions commonly referred to as restricted eigenvalue (RE) conditions in the lasso literature \cite[]{bickel2009,loh2012}. Intuitively speaking, Condition~\ref{cond.fvar.RE} implies that $\btheta^\T \widehat \bGamma \btheta/\|\btheta\|^2$ is strictly positive as long as $\|\btheta\|_1$ is not large relative to $\|\btheta\|.$ Denote the estimation error by $\bDelta_j=\widehat\bB_j-\bB_j,$ this condition ensures that $\langle\langle\bDelta_j, \widehat\bGamma\bDelta_j\rangle\rangle \ge \tau_2 \|\bDelta_j\|_{\tF}^2/2$ if $\tau_2\ge 32\tau_1 q^2 s$, as stated in Theorem~\ref{thm.vfar} below.

\begin{condition}
\label{cond.fvar.eigen}	
	There exist some positive constants $C_{\lambda}$ and $C_{\phi}$ independent of $(n,p,q)$ such that
    \begin{equation}
    \label{fvar.eigen}	
    \begin{split}
    \max_{1 \le j \le p, 1 \le l \le q} \left|\frac{\widehat{\lambda}_{jl}^{-1/2} - \lambda_{jl}^{-1/2}}{\lambda_{jl}^{-1/2}}\right| \le C_{\lambda}  \cM_1(\mF_{\bX})\sqrt{\frac{\log (pq)}{n}},\\
    \max_{1 \le j \le p, 1 \le l \le q}\|\widehat \phi_{jl} - \phi_{jl}\| \le C_{\phi} \cM_1(\mF_{\bX})q^{\alpha + 1} \sqrt{\frac{\log (pq)}{n}}.
    \end{split}
    \end{equation}
\end{condition}

\begin{condition}
	\label{cond.fvar.max.error}
	There exits some positive constant $C_{E}$ independent of $(n,p,q)$ such that
	\begin{equation}
	\label{fvar.max.error}
	\Big\|\widehat \bY_j- \widehat \bGamma\bB_j\Big\|_{\max}^{(q)} \leq C_E \cM_1(\mF_{\bX}) s_j \Big\{q^{\alpha + 2}\sqrt{\frac{\log (pq)}{n}} + q^{-\beta +1} \Big\},~j=1,\ldots,p.
	\end{equation}
	\end{condition}
	
Condition~\ref{cond.fvar.eigen} and \ref{cond.fvar.max.error} are two deviation conditions, which guarantee the good behaviours of relevant estimated terms by controlling their deviation bounds. Specifically, Condition~\ref{cond.fvar.max.error} ensures that $\widehat\bY_j$ and $\widehat\bGamma$ are nicely concentrated around their population versions. See also similar deviation conditions in the lasso literature \cite[]{loh2012,basu2015a}.

We are now ready to present the theorem on the convergence rate of the VFAR estimate. 

\begin{theorem}
\label{thm.vfar}
    Suppose that Conditions \ref{cond.cov.func}--\ref{cond.fvar.max.error} hold with $\tau_2\ge 32\tau_1 q^2 s.$ Then, for any regularization parameter, $\gamma_{nj} \ge 2C_E \cM_1(\mF_{\bX}) s_j\big\{q^{\alpha + 2}(\log(pq)/n)^{1/2} + q^{-\beta +1} \big\},$ $\gamma_n=\max_{j}\gamma_{nj}$ and $q^{\alpha/2}s\gamma_n\to 0$ as $n, p, q \to \infty,$ any minimizer $\widehat\bB_j$ of (\ref{vfar.crit.theory}) satisfies
$$
\|\widehat \bB_{j} - \bB_{j}\|_F \le \frac{24 s_j^{1/2}\gamma_{nj}}{\tau_2},
~~\|\widehat \bB_{j} - \bB_{j}\|_{1}^{(q)} \le \frac{96s_j\gamma_{nj}}{\tau_2} ~ \text{ for }j =1,\ldots,p,
$$
and the estimated transition function, $\widehat \bA,$ satisfies
\begin{equation}
\label{err.A}
\|\widehat \bA - \bA\|_\infty \leq \frac{96\alpha^{1/2}q^{\alpha/2}s\gamma_{n}}{c_0^{1/2}\tau_2}\Big\{1+o(1)\Big\}.
\end{equation}
\end{theorem}

The convergence rate of $\widehat\bA$ under functional matrix $\ell_{\infty}$ norm is governed by dimensionality parameters ($n, p, s$) and 
internal parameters ($\cM_1(\mF_{\bX}),$ $q,$ $\tau_1,$ $\tau_2,$ $\alpha, \beta$).
We provide three remarks for the error bound in (\ref{err.A}). 
First, it is easy to see that larger values of $\alpha$ (tighter eigengaps) or $\cM_1(\mF_{\bX})$ (less stable process of $\{\bX_t(\cdot)\}$) or $s$ (denser structure in $\bA$) yield a slower convergence rate, while enlarging $\beta$ or $\tau_2$ will increase the entrywise smoothness in $\bA$ or the curvature of the RE condition, respectively, thus resulting in a faster rate. Second, the convergence rate consists of two terms corresponding to the variance-bias tradeoff as commonly considered in nonparametric statistics. Specifically, the variance is of the order $O_P\big[\cM_1(\mF_{\bX})s^2q^{(3\alpha+4)/2}\{\log(pq)/n\}^{1/2}\big]$ and the bias term is bounded by $O\big\{\cM_1(\mF_{\bX})s^2q^{(\alpha-2\beta+2)/2}\big\}.$ To balance both terms, we can choose an optimal $q$ satisfying $\log(pq)q^{2\alpha+2\beta+2}\asymp n$, which leads to $q \asymp \{n /\log (p \vee n)\}^{1/(2\alpha + 2 \beta + 2)}.$ Third, when each $X_{tj}(\cdot)$ is finite-dimensional, although the truncation step is no longer required, the FPC scores still need to be estimated. The resulting convergence rate becomes $O_P\big\{\cM_1(\mF_{\bX})s^2(\log p/n)^{1/2}\big\},$ which is slightly different from that of the high-dimensional VAR estimate in \cite{basu2015a}.

Finally, we turn to 
verify that, when $\bX_1(\cdot),\dots,\bX_n(\cdot)$ are drawn from model (\ref{vfar1}), Conditions~\ref{cond.fvar.RE}--\ref{cond.fvar.max.error} are satisfied with high probability as stated in the following  
Propositions~\ref{res.re}--\ref{prop_Error_max}. Before presenting these propositions, we list two regularity conditions.


\begin{condition}
\label{cond_min_bound}
For $\bSigma_0=\big(\Sigma_{jk}^{(0)}\big)_{1\leq j,k\leq p},$ we denote by $\bD_0 = \text{diag}(\Sigma_{11}^{(0)}, \dots, \Sigma_{pp}^{(0)})$ the diagonal function. 
The infimum $\underline{\mu}$ of the functional Rayleigh quotient of $\bXi_0$ relative to $\bD_0,$ defined as follows, is bounded below by zero, i.e.
\begin{eqnarray}
\underline{\mu} = \inf_{\bPhi \in \widebar{\cH}_0} \frac{\big \langle\bPhi, \bXi_0(\bPhi) \big \rangle_{\cH}}
{\big \langle\bPhi, \bD_0(\bPhi) \big \rangle_{\cH}}  > 0,
\end{eqnarray}
where $\widebar{\cH}_0 = \{\bPhi \in \cH: \big \langle\bPhi, \bD_0(\bPhi) \big \rangle_{\cH} \in (0, \infty)\}.$
\end{condition}

\begin{condition}
\label{cond_q_n}
The sample size $n$ satisfies the bound $n \gtrsim q^{4\alpha+2}\cM_{1}^2(\mF_{\bX}) \log(pq).$
\end{condition}

In Condition~\ref{cond_min_bound}, the lower bound on $\underline{\mu},$ chosen as the curvature $\tau_2$ in the proof of Proposition~\ref{res.re}, can be understood as requiring the minimum eigenvalue of the correlation function for $\bX_t(\cdot)$ to be bounded below by zero. Specially, if $X_{tj}(\cdot)$ is $d_j$-dimensional for $j=1,\dots,p$ 
it is easy to show that $\underline{\mu}$ reduces to the minimum eigenvalue of the correlation matrix for the ($\sum_jd_j$)-dimensional vector,  
$\bxi_t=(\xi_{t11}, \dots, \xi_{t1d_1}, \dots, \xi_{tp1}, \dots, \xi_{tpd_p})^{\T}.$ 
Condition~\ref{cond_q_n} is required here due to its presence in Theorems~\ref{thm_lam_phi} and \ref{thm_cov}.


\begin{proposition} (Verify Condition~\ref{cond.fvar.RE})
	\label{res.re}
	Suppose that Conditions \ref{cond.bd.fsm}--\ref{cond_eigen} and \ref{cond_min_bound}--\ref{cond_q_n} hold. Then there exist three positive constants $C_{\Gamma}$, $c_5$ and $c_6$ independent of $(n,p,q)$ such that
	\begin{eqnarray*}
		\btheta^T \widehat \bGamma \btheta \ge  \underline{\mu} \big\|\btheta\big\|^2_2 - C_\Gamma \cM_1(\mF_{\bX}) q^{\alpha +1}\sqrt{\frac{\log (pq)}{n}}\big\|\btheta\big\|_1^2
	\end{eqnarray*}
with probability greater than $1 - c_5(pq)^{-c_6}$.
\end{proposition}
\begin{proposition} (Verify Condition~\ref{cond.fvar.eigen})
\label{prop_Error_eigen}
Suppose that Conditions~\ref{cond.bd.fsm}--\ref{cond_eigen} and \ref{cond_q_n} hold. Then there exist four positive constants $C_{\phi}$, $C_{\lambda}$ $c_5$ and $c_6$ independent of $(n,p,q)$ such that 
(\ref{fvar.eigen}) holds with probability greater than $1 -c_5(pq)^{-c_6}.$
\end{proposition}

\begin{proposition} (Verify Condition~\ref{cond.fvar.max.error})
\label{prop_Error_max}
    Suppose that Conditions \ref{cond.bd.fsm}--\ref{cond.fvar.bias} and \ref{cond_q_n} hold. Then there exist three positive constants $C_E$, $c_5$ and $c_6$ independent of $(n,p,q)$ such that 
    (\ref{fvar.max.error}) holds
	with probability greater than $1 -c_5(pq)^{-c_6}.$
\end{proposition}

Propositions \ref{res.re}--\ref{prop_Error_max} can be proved by applying the convergence results in Theorems~\ref{thm_lam_phi} and \ref{thm_cov}. With suitable choices of common constants $c_5$ and $c_6$ in Propositions \ref{res.re}--\ref{prop_Error_max}, we can show that the joint probability for the three events corresponding to the non-asymptotic upper bounds in
(\ref{bd_max_lam_phi}) and (\ref{bd_score_max}) is greater than $1-c_5(pq)^{-c_6}.$ Consequently, 
with probability greater than $1-c_5(pq)^{-c_6},$ the estimate $\widehat\bA$ satisfies the error bound in (\ref{err.A}).

\subsection{Simulation studies}
\label{sec.sim}
In this section, we conduct a number of simulations to compare the finite-sample performance of our proposed method to potential competitors.

In each simulated scenario, we generate functional variables by $X_{tj}(u)=\bs(u)^{\T}\btheta_{tj}$ for $j=1,\dots,p$ and $u \in \cU=[0,1],$ where $\bs(\cdot)$ is a $5$-dimensional Fourier basis function and each $\btheta_{t}=(\btheta_{t1}^{\T},\dots,\btheta_{tp}^{\T} )^{\T} \in \eR^{5p}$ is generated from a stationary VAR(1) process, $\btheta_t=\bB\btheta_{t-1}+\bfeta_{t},$ with block transition matrix $\bB \in \eR^{5p\times 5p},$ whose $(j,k)$-th block is given by $\bB_{jk}$ for $ j,k=1,\dots,p$ and innovations $\bfeta_{t}$'s being independently sampled from $N({\bf 0}, \bI_{5p}).$ To mimic a real data setting, we generate observed values, $W_{tjs},$ with measurement errors, $W_{tjs}=X_{tj}(u_s)+e_{tjs},$ from $T=50$ equally spaced time points, $0=u_1,\dots,u_T=1$ with errors $e_{tjs}$'s being randomly sampled from $N(0,0.5^2).$ In our simulations, we generate $n=100$ or $n=200$ observations of $p=40$ or $p=80$ functional variables, and we aim to show that, although our method is developed for fully observed functional time series, it still works well even for the dense design with measurement error. It is worth noting that, as discussed in Section~\ref{sec.intro}, the VFAR estimation is naturally a very high-dimensional problem. For example, to fit a VFAR(1) model under our most ``low-dimensional" setting with $p=40$ and $n=200,$ we need to estimate $40^2\times 5^2 = 40,000$ parameters based on only 200 observations.

According to Section~\ref{ap.vfar1.sim} of the Supplementary Material, $\bX_t(\cdot)$ follows from a VFAR(1) model in (\ref{vfar1.ex}), 
where $\varepsilon_{tj}(u)=\bs(u)^{\T}\bfeta_{tj}$ and autocoefficient functions satisfy $A_{jk}(u,v)=\bs(u)^{\T}\bB_{jk}\bs(v)$ for $j,k=1,\dots,p, (u,v) \in \cU^2.$ Hence,
the functional sparsity structure in $\bA$ can be characterized by the block sparsity pattern in $\bB.$ In the following, we consider two different scenarios to generate $\bB.$
\begin{itemize}
    \item[(i)] {\bf Block sparse}. We generate a block sparse $\bB$ without any special structure. Specifically, we generate  $\bB_{jk}=w_{jk}\bC_{jk}$ for $j,k=1\dots,p,$ where entries in $\bC_{jk}$ are randomly sampled from $N(0,1)$ and $w_{jk}$'s are generated from $\{0,1\}$ under the constraint of $\sum_{k=1}^p w_{jk}=5$ for each $j,$ such that the same row-wise cardinality for $\bB$ can be produced in a blockwise fashion. To guarantee the stationarity of $\{\bX_t(\cdot)\},$ we rescale $\bB$ by $\iota\bB/\rho(\bB),$ where $\iota$ is generated from Unif[0.5,1].

    \item[(ii)] {\bf Block banded}. We generate a block banded $\bB,$ with entries in $\bB_{jk}$ being randomly sampled from $N(0,1)$ if $|j-k|\leq 2,$ and being zero at other locations. $\bB$ is then rescaled as described in (i).
\end{itemize}

For each $j=1, \dots,p,$ we perform regularized FPCA \cite[]{Bramsay1} on observations $\{W_{tjs}\}$ to obtain smoothed estimates of $\phi_{jl}(\cdot)$'s, and use $5$-fold cross-validation to choose $q_j$ and the smoothing parameter, the details of which are presented in Sections~\ref{ap.select.tune} and \ref{ap.rfpca} of the Supplementary Material. Typically $q_j=$ 4, 5 or 6 are selected in our simulations. 
To choose the regularization parameters $\gamma_{nj}$'s, there exist a number possible methods such as AIC/BIC and cross-validation. While the third one is computationally intensive and remains largely unexplored in the literature for serially dependent functional observations, we take an approach motivated by the information criterion for sparse additive models \cite[]{voorman2014}. For each $j,$ our proposed information criterion is
\begin{equation}
\label{ic_vfar}
\text{IC}_j(\gamma_{nj}) = n \text{log} \Big\{\big\|\widehat\bV_{j}^{(0)} - \sum_{h = 1 }^L\sum_{k = 1 }^p \widehat\bV_{k}^{(h)}\widehat \bPhi_{jk}^{(h)}(\gamma_{nj})  \big\|_F^2\Big\} + \kappa_n \text{df}_j(\gamma_{nj}),		
\end{equation}
where $\kappa_n=2$ and $\log n$ correspond to AIC and BIC, respectively, and $\text{df}_j(\gamma_{nj})$ is the effective degrees of freedom used in fitting (\ref{vfar.crit}) with
\begin{equation}
\label{vfar.df}
\text{df}_j(\gamma_{nj}) = \sum_{h = 1 }^L \sum_{k = 1 }^p \Bigg\{ I\Big((h,k): ||\widehat \bPhi_{jk}^{(h)}(\gamma_{nj})||_F \neq 0\Big) +  (q_{j} q_{k} -1) \frac{||\widehat\bV_{k}^{(h)}\widehat \bPhi_{jk}^{(h)}(\gamma_{nj})||_F^2}{||\widehat\bV_{k}^{(h)}\widehat \bPhi_{jk}^{(h)}(\gamma_{nj})||_F^2 + \gamma_{nj}} \Bigg\}.
\end{equation}

We compare our proposed $\ell_1/\ell_2$-penalized LS estimate using all selected principal components, namely $\ell_1/\ell_2$-$\text{LS}_{\text{a}},$ to its two competitors. One method, $\ell_1/\ell_2$-$\text{LS}_{2},$ relies on minimizing $\ell_1/\ell_2$-penalized LS based on the first two estimated principal components, which capture partial curve information. The other approach, $\ell_1$-$\text{LS}_{1},$ projects the functional data into a standard format by computing the first estimated FPC score and then implements an $\ell_1$ regularization approach \cite[]{basu2015a} for the VAR estimation on this data. We examine the sample performance of three approaches, $\ell_1/\ell_2$-$\text{LS}_{\text{a}},$ $\ell_1/\ell_2$-$\text{LS}_{2},$ and $\ell_1$-$\text{LS}_{1}$ in terms of model selection consistency and estimation accuracy.
\begin{itemize}
    \item {\bf Model selection}. We plot the true positive rates against false positive rates, respectively defined as
    $\frac{\#\{(j,k):  ||\widehat A_{jk}^{(\bgamma_n)}||_{\cS}\neq0 \text{ and } ||A_{jk}||_{\cS}\neq 0\}}{\#\{(j,k): ||A_{jk}||_{\cS}\neq 0\}}$ and
    $\frac{\#\{(j,k):  ||\widehat A_{jk}^{(\bgamma_n)}||_{\cS}\neq0 \text{ and } ||A_{jk}||_{\cS}= 0\}}{\#\{(j,k): ||A_{jk}||_{\cS}= 0\}},$ over a sequence of $\bgamma_n=(\gamma_{n1}, \dots, \gamma_{np})$ values to produce a ROC curve. We compute the {area under the ROC curve} (AUROC) with values closer to one indicating better performance in recovering the functional sparsity structure in $\bA.$

    \item {\bf Estimation error}. We calculate the relative estimation accuracy for $\widehat\bA$ by $\|\widehat{\bA} - \bA\|_{\tF}/\|\bA\|_{\tF},$ where $\widehat\bA$ is the regularized estimate based on the optimal regularization parameters selected by minimizing AICs or BICs in (\ref{ic_vfar}). 
\end{itemize}

\begin{table}[t]
	\caption{\label{auroc.table} The mean and standard error (in parentheses) of AUROCs over 100 simulation runs. The best values are in bold font.}
	\begin{center}
		\resizebox{5.2in}{!}{
			\begin{tabular}{cccccccc}		
				\hline
				& &  \multicolumn{3}{c}{Model~(i)}& \multicolumn{3}{c}{Model~(ii)} 
				\\\hline
				$n$ & $p$ & $\ell_1/\ell_2$-$\LS_a$  & $\ell_1/\ell_2$-$\LS_2$ & $\ell_1$-$\LS_1$ & $\ell_1/\ell_2$-$\LS_a$ & $\ell_1/\ell_2$-$\LS_2$ & $\ell_1$-$\LS_1$ \\
				
				\multirow{2}{*}{100} & 40 & {\bf 0.840(0.018)} & 0.690(0.019) & 0.591(0.023) & {\bf 0.872}(0.016) & 0.719(0.022) & 0.609(0.024)\\
				 & 80 & {\bf 0.829(0.015)} & 0.682(0.017) & 0.585(0.015) & {\bf 0.869}(0.014) & 0.714(0.017) & 0.600(0.017) \\
				\hline
				\multirow{2}{*}{200} & 40 & {\bf 0.951(0.011)} & 0.764(0.020) & 0.616(0.021) & {\bf 0.971}(0.006) & 0.795(0.018) & 0.639(0.023) \\
				 & 80 & {\bf 0.948(0.010)} & 0.770(0.017) & 0.626(0.015) & {\bf 0.969}(0.005) & 0.799(0.014) & 0.644(0.015) \\
				\hline
			\end{tabular}
		}	
	\end{center}
\end{table}

To investigate the support recovery consistency, we report the average AUROCs of three comparison methods under both model settings in Table~\ref{auroc.table}. 
In all simulations, we observe that $\ell_1/\ell_2$-$\text{LS}_{\text{a}}$ with most of curve information being captured, provides highly significant improvements over its two competitors and $\ell_1$-$\text{LS}_{1}$ gives the worst results. 
To evaluate the estimation accuracy, Table~\ref{err.table} presents numerical results of relative errors of different regularized estimates. We also report the performance of the LS estimate in the oracle case, where we know locations of non-zero entries of $\bA$ in advance. Several conclusions can be drawn from Table~\ref{err.table}. First, in all scenarios, the proposed BIC-based $\ell_1/\ell_2$-$\text{LS}_{\text{a}}$ method provides the highest estimation accuracy among all the comparison methods. Second, the performance of AIC-based methods severally deteriorate in comparison with their BIC-based counterparts. Given the high-dimensional, dependent and functional natural of the model structure, computing the effective degrees of freedom in (\ref{vfar.df}) leads to a very challenging task and requires further investigation. In practice, with some prior knowledge about the targeted graph density under a FNGC framework, one can choose the values of $\gamma_{nj}$ that result in the graph with a desired sparsity level. Third, $\LS_{\text{oracle}}$ estimates give much worse results than BIC-based regularized estimates. This is not surprising, since even in the ``large $n,$ small $p$" scenario, e.g. $n=100, p=40$ for Model~(i), implementing LS requires estimating $5\times5^2=125$ parameters based on only 100 observations, which intrinsically results in a high-dimensional estimation problem.

\begin{table}[t]
	\caption{\label{err.table} The mean and standard error (in parentheses) of relative estimation errors of $\widehat\bA$ over 100 simulation runs. The best values are in bold font.}
	\begin{center}
		\resizebox{5.3in}{!}{
			\begin{tabular}{ccccccc}
				\hline
				Model &$(n,p)$ & $\bgamma_{n}$ & $\ell_1/\ell_2$-$\LS_{\text{a}}$  & $\ell_1/\ell_2$-$\LS_2$ & $\ell_1$-$\LS_1$ & $\LS_{\text{oracle}}$\\
				\hline
				\multirow{8}{*}{(i)}&\multirow{2}{*}{(100,40)} & $\aic$& 1.783(0.044) & 1.455(0.032) & 1.047(0.008) & \multirow{2}{*}{1.483(0.044)} \\
				&& $\bic$& {\bf 0.971(0.007)} & 0.997(0.002) & 1.002(0.002) & \\
				&\multirow{2}{*}{(100,80)} & $\aic$& 1.546(0.037) & 1.576(0.030) & 1.124(0.013) & \multirow{2}{*}{1.530(0.050)} \\
				&                         & $\bic$& {\bf 0.991(0.003)} & 0.999(0.001) & 1.002(0.001) &  \\
				&\multirow{2}{*}{(200,40)} & $\aic$& 1.419(0.035) & 1.159(0.016) & 1.020(0.004) & \multirow{2}{*}{0.990(0.029)}  \\
				&                          & $\bic$& {\bf 0.850(0.012)} & 0.989(0.003) & 0.999(0.002) & \\
				&\multirow{2}{*}{(200,80)} & $\aic$& 1.544(0.045) & 1.350(0.022) & 1.032(0.004) & \multirow{2}{*}{1.016(0.033)}  \\
				&                         & $\bic$& {\bf 0.915(0.013)} & 0.994(0.002) & 1.000(0.001) & \\ \hline
				\multirow{8}{*}{(ii)}&\multirow{2}{*}{(100,40)} & $\aic$ & 1.679(0.040) & 1.400(0.026) & 1.039(0.008) & \multirow{2}{*}{1.363(0.039)}  \\
				&                         & $\bic$ & {\bf 0.957(0.009)} & 0.995(0.002) & 1.000(0.002) & \\
				&\multirow{2}{*}{(100,80)} & $\aic$ & 1.435(0.028) & 1.489(0.020) & 1.105(0.012) & \multirow{2}{*}{1.383(0.036)}    \\
				&                        & $\bic$ & {\bf 0.983(0.004)} & 0.998(0.001) & 1.001(0.001) &  \\
				&\multirow{2}{*}{(200,40)} & $\aic$ & 1.329(0.028) & 1.128(0.013) & 1.014(0.004) & \multirow{2}{*}{0.909(0.024)}   \\
				&                         & $\bic$ & {\bf 0.824(0.009)} & 0.985(0.003) & 0.998(0.002) & \\
				&\multirow{2}{*}{(200,80)} & $\aic$ & 1.428(0.027) & 1.296(0.013) & 1.024(0.004) & \multirow{2}{*}{0.926(0.019)} \\
				&                         & $\bic$ & {\bf 0.881(0.010)} & 0.990(0.002) & 0.998(0.001) &  \\
				\hline
			\end{tabular}
		}	
	\end{center}
\end{table}

\subsection{Real data analysis}
\label{sec.real}
In this section, we apply our proposed method to a public financial dataset, which was downloaded from \url{https://wrds-web.wharton.upenn.edu/wrds}. The dataset consists of high-frequency observations of prices for a collection of S\&P100 stocks from $n=251$ trading days in year 2017. We removed several stocks for which the data are not available during the observational period. See Table~\ref{name.tb} in the Supplementary Material for tickers, company names and classified sectors of the inclusive $p=98$ stocks.
We then obtained data at a sampling frequency one-minute per datum such that the impact of microstructure noise is reduced \cite[]{zhang2005}. The daily trading period (9:30-16:00) is thus converted to $\cU=[0,T]$ with $T=390$ minutes. 
Let $P_{tj}(u_k)$ ($t=1, \dots,n, j=1, \dots, p, k=1, \dots, T$) be the price of the $j$-th stock at intraday time $u_k$ on the $t$-th trading day. We denote the {cumulative intraday return} (CIDR) trajectory, in percentage, by $r_{tj}(u_k)=100[\log\{P_{tj}(u_k)\}-\log\{P_{tj}(u_1)\}]$ \cite[]{horvath2014}. This transformation not only guarantees the shape of CIDR curves nearly the same as original daily price curves, but also makes the assumption of stationarity for curves more plausible. According to the definition, CIDR curves always start from zero enforcing level stationary and the logarithm helps reduce potential scale inflation. See~\cite{horvath2014} for an empirical study on testing the stationary of CIDR curves.

Our main target is to construct a directed graph under a FNGC framework and to display the Granger-type casual relationships among intraday returns of different stocks. To achieve this goal, we first center all the series of $\{r_{tj}(\cdot)\}$ about their empirical means and then apply the three-step approach to fit a sparse VFAR(1) model on the demeaned curves. In the first step, we implement regularized FPCA to smooth these curves and therefore the effects of microstructure noise are further reduced to be almost negligible.
See Sections~\ref{ap.select.tune}--\ref{ap.rfpca} and \ref{ap.alg.vfar} of the Supplementary Material for details on regularized FPCA with the selection of relevant tuning parameters and block FISTA algorithm, respectively. 

\begin{figure*}[ht]
	\centering
	\includegraphics[width=6.7cm,height=6.7cm]{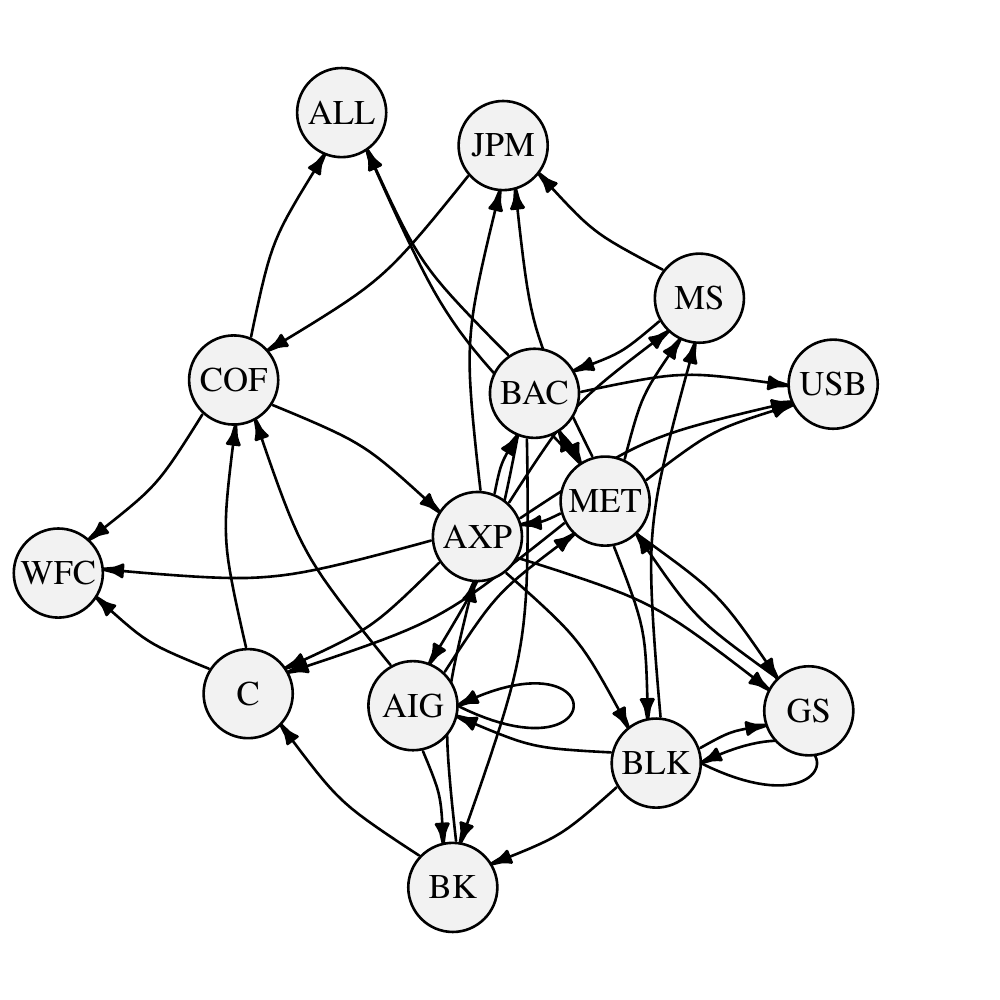}
	\includegraphics[width=6.7cm,height=6.7cm]{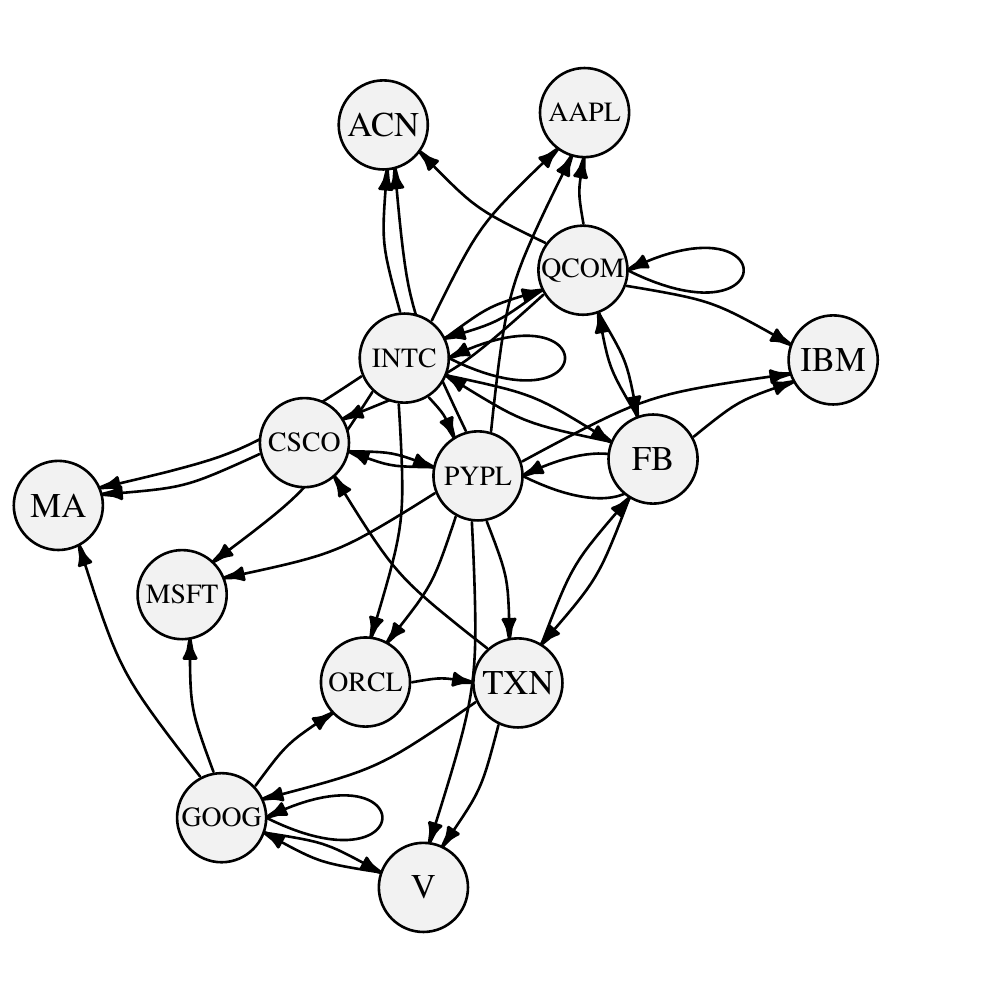}
	\vspace{-0.35cm}
	\caption{\label{Net.plot14}{\it \small Left and right graphs plot the directed networks with indegree=3 for $p=14$ stocks in the financial and IT sectors, respectively.}}
\end{figure*}

\begin{figure*}[t]
	\centering
	\includegraphics[width=9cm,height=7.2cm]{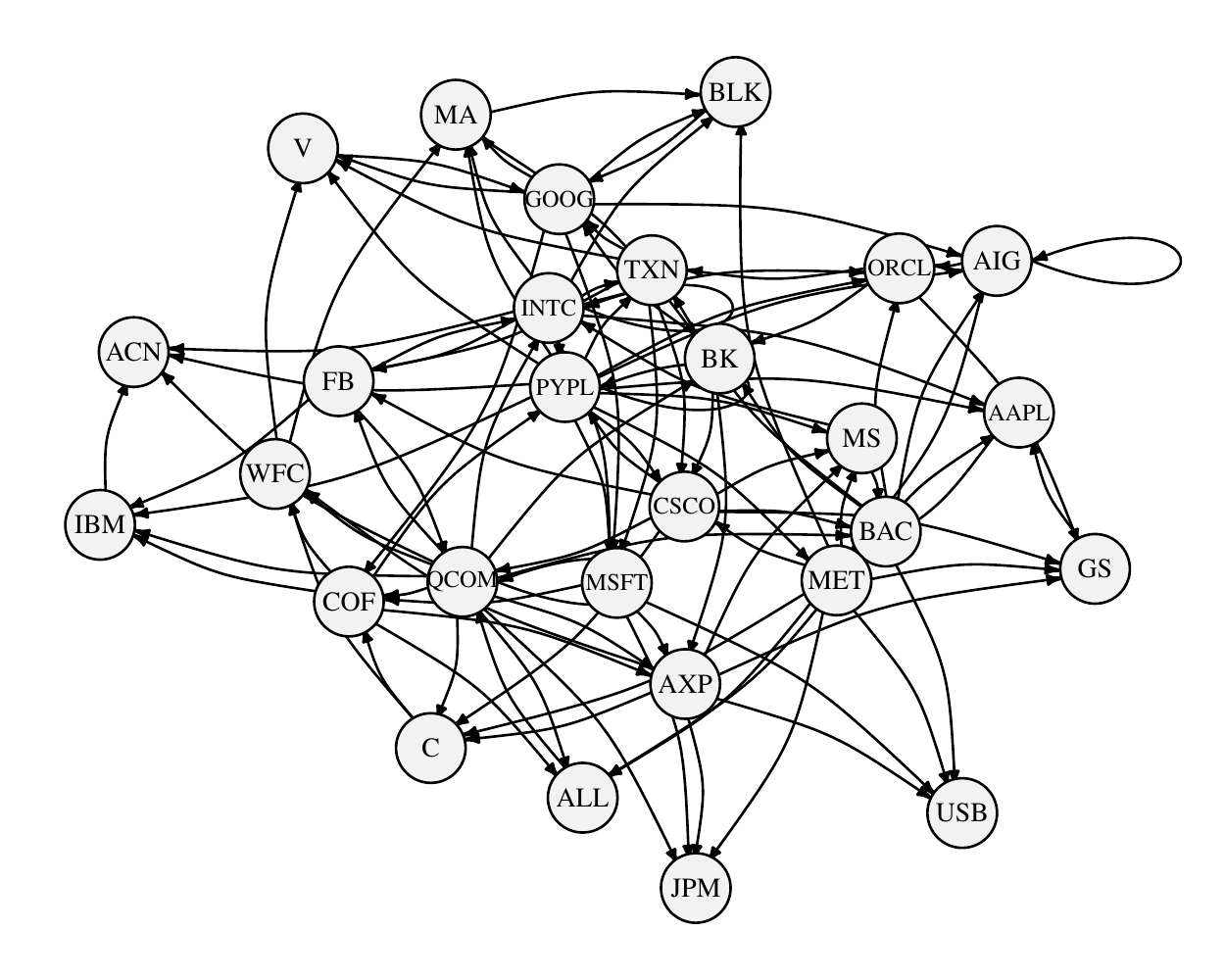}
	\vspace{-0.4cm}
	\caption{\label{Net.plot28}{\it \small The directed graph with indegree=4 for $p=28$ stocks in the financial and IT sectors.}}
\end{figure*}

To better visualize and interpret the network, we focus on $p=14$ stocks in the financial and information technology (IT) sections, respectively and set the row-wise sparsity to 3/14, i.e. each node receives connections from 3 (indegree) out of 14 nodes.  A more systematic method for determining the network sparsity level, e.g. via a significance testing, needs to be developed. Figure~\ref{Net.plot14} displays two directed networks based on the identified sparsity structures in $\widehat\bA$ (estimated transition function) for stocks in two sectors.
It suggests that ``MET" (Metlife) and ``INTC" (Intel) together with ``PYPL" (PayPal), placed in the center of each network, provide the lowest levels of column-wise sparsity in the financial and IT sectors, respectively, thus resulting in 
highest Granger-type causal impacts on all the stocks in terms of their CIDR curves.
Moreover, we consider $p=28$ stocks in both financial and IT sectors. Setting the row-wise sparsity to 1/7, we plot a larger directed graph in Figure~\ref{Net.plot28}. We observe that more IT companies, e.g. PayPal, Qualcomm and Intel, have relatively higher causal impacts. Interestingly, PayPal, as a leading financial technology (FinTech) company belonging to both financial and IT sectors leads to the the highest causal influence on others. See also Section~\ref{supp_sec_real} of the Supplementary Material for additional empirical analysis.

\section{Discussion}
\label{sec.discussion}

We identify several important directions for future research. The first topic considers the functional extension of high-dimensional factor models \cite[]{bai2002,lam2012}, where the observations, $\bX_t(\cdot)$'s,
can be decomposed as the sum of two unobservable and mutually orthogonal components
\begin{equation}
    \label{ffm}
    \bX_t(u) = \bW_t(u) + \beps_t(u), ~t=1,\dots,n, ~u \in \cU.
\end{equation}
Here $\bW_t(\cdot)=\bB \bbf_t(\cdot)$ are the common components driven by $r$ (much smaller than $p$) functional factors $\bbf_t(\cdot)=\big(f_{t1}(\cdot), \dots, f_{tr}(\cdot)\big)^{\T},$
$\bB \in {\eR}^{p\times r}$ is the factor loading matrix 
and $\beps_t(\cdot)=\big(\epsilon_{t1}(\cdot), \dots, \epsilon_{tp}(\cdot)\big)^{\T}$ are
idiosyncratic components. For each $h \in \eZ$ and $\{\bX_t(\cdot)\}_{t=1}^{n},$ denote  $\bSigma^{X}_{h}(u,v) = \cov\{\bX_t(u), \bX_{t+h}(v)\}$ and its sample estimator by $\widehat{\bSigma}^{X}_{h}(u,v).$ To estimate such functional factor model~(\ref{ffm}), we discuss two different approaches. (i) The first one is based on the following integrated covariance decomposition,
\begin{equation}
\label{decomp_cov}
 \int\int \bSigma^{X}_{0}(u,v) dudv=\bB \Big\{\int\int \bSigma^{f}_{0}(u,v) dudv\Big\}\bB^{\T} + \int\int \bSigma^{\epsilon}_{0}(u,v) dudv.   
\end{equation}
Intuitively, by imposing some eigenvalue conditions on two terms on the right-hand side of (\ref{decomp_cov}) similar to those in \cite{fan2013}, the above decomposition is asymptotically identified as $p \rightarrow \infty$ and hence $\bB$ can be recovered by performing an eigenanalysis of $\int\int \widehat\bSigma^{X}_{0}(u,v) dudv.$ 
(ii) If $\{\beps_t(\cdot)\}$ follows a white noise process, then, inspired from \cite{lam2012} and the fact that $ \bSigma^{X}_{h}(u,v)=\bB\bSigma^{f}_{h}(u,v)\bB^{\T}$ for $h\geq 1,$ an autocovariance-based procedure can be developed to estimate model~(\ref{ffm}). To theoretically support both estimation procedures under high-dimensional settings, the main challenge is to investigate convergence properties of $\widehat \bSigma_{h}^{X}-\bSigma_{h}^{X}$ for $h=0,1,\dots,$ and hence our concentration results in Theorem~\ref{thm_con_op} and Proposition~\ref{thm_sigma_k} can be applied.

Second, in the dimension reduction step of the three-step procedure, one can also perform dynamic FPCA \cite[]{hormann2015} based on $\{X_{tj}(\cdot)\}_{t=1}^n$ for each $j.$ Such dimension reduction technique provides an optimal truncated approximation for functional time series, but is computationally intensive as it relies on the eigenanalysis of spectral density functions $f_{\bX,\theta}$ with sample estimators given by $\widehat f_{\bX,\theta}=(2 \pi)^{-1}\sum_h w_H(h)\widehat\bSigma_{h} \exp(- i h \theta),$ where $w_{H}(\cdot) = w(\cdot/H)$ is some appropriate weight function  with $H$ (the lag window size). 
To provide theoretical guarantees for relevant estimated terms under a dynamic FPCA framework similar to those in Theorems~\ref{thm_lam_phi} and \ref{thm_cov}, our established non-asymptotic error bounds on
$\widehat \bSigma_{h}$ for $h \in \eZ$ become applicable. 

Third, within the proposed functional stability measure framework, we believe our established concentration results can be extended beyond Gaussian functional time series to accommodate linear processes with functional sub-Gaussian errors. 
It is also interesting to develop suitable concentration results for heavy-tailed non-Gaussian functional time series under the functional generalization of the $\beta$-mixing condition in \cite{wong2020}. 

Fourth, it is of great interest to develop new inference tools for high-dimensional functional time series and to apply these techniques to quantify deviations of autocoefficient functions in sparse VFAR models. Fifth, our analysis is based on the estimation where smoothness parameters are assumed to be known. It is interesting to develop adaptive estimation procedures that do not require the knowledge of the parameter space and automatically adjust to the smoothness properties. However, this would pose complicated challenges under the high-dimensional, functional and dependent setting we consider.

These topics are beyond the scope of the current paper and will be pursued elsewhere.



 \section*{Appendix}
 \setcounter{equation}{0}
 \renewcommand{\theequation}{A.\arabic{equation}}

\appendix
\section{Examples satisfying Condition~\ref{cond.bd.fsm}}
\label{ap.fsm.ex}
It is clear that Condition~\ref{cond.bd.fsm} holds for finite-dimensional functional data. In the following, we give five illustrative examples for infinite-dimensional functional data, where the upper bounds of their functional stability measures can be easily controlled. As long as the denominator of $\cM(\mF_{\bX})$ is arbitrarily small, the numerator can also be arbitrarily small, Condition~\ref{cond.bd.fsm} in this sense still holds for a large class of infinite-dimensional functional data including, but not limited to the five examples below.



{\it Example~1}: When $\bX_t$ and $\bX_s$ are independent for any $t \neq s$ and each $X_{jt}(\cdot)$ is infinite dimensional for $j=1,\dots, p,$ then $\cM(\mF_{\bX})=1<\infty.$

{\it Example~2}: For convenience, we consider an univariate functional time series example. Let
$$
X_{t}(u) = \sum_{l =1}^\infty \lambda_{l}^{1/2} a_{tl} \phi_{l}(u), ~u \in [0,1], ~\lambda_l =  l^{-\alpha}, ~t=1,\dots,n,
$$
where $\alpha>1,$ $\{\phi_l(\cdot)\}_{l=1}^\infty$ are orthonormal basis functions, 
$$
\left (\begin{array}{cc}
a_{t(2l-1)} \\
a_{t(2l)}
\end{array}
\right ) =    \bA
\left (\begin{array}{cc}
a_{(t-1)(2l-1)} \\
a_{(t-1)(2l)}
\end{array}
\right ) +
\left (\begin{array}{cc}
e_{t(2l-1)} \\
e_{t(2l)}
\end{array}
\right ),
$$
$$
\bA = \frac{1}{4} \left (\begin{array}{cc}
2& 1 \\
1 & 2
\end{array}
\right), ~~
\left (\begin{array}{cc}
e_{t(2l-1)} \\
e_{t(2l)}
\end{array}
\right ) \stackrel{i.i.d.}{\sim} N\left \{\left (\begin{array}{cc}
0 \\
0
\end{array}
\right ),\left (\begin{array}{cc}
\frac{11}{16}& -\frac{1}{4} \\
-\frac{1}{4} & \frac{11}{16}
\end{array}
\right) \right\},
$$
and $\big(e_{t(2l-1)},
e_{t(2l)}\big)^{\T}$'s are independent for all $t$ and $l.$ Here the covariance structure of $\big(e_{t(2l-1)},
e_{t(2l)}\big)^{\T}$ can guarantee that the covariance matrix of $\big(\xi_{t(2l-1)},
\xi_{t(2l)}\big)^{\T}$ is an identity matrix, which means that $\xi_{(t-1)(2l-1)}$ and $\xi_{(t-1)(2l)}$ are independent. 
In this setting, some specific calculations yield that 
$
\Sigma_{0}(u,v) = E\{X_{t}(u)X_{t}(v)\} = \sum_{l=1}^\infty \lambda_{l}\phi_{l}(u)\phi_{l}(v)
$
and for $h = 1,2,\ldots,$
$$
	 \Sigma_{h}(u,v) = E\{X_{t}(u)X_{(t+h)}(v)\} 
	= \sum_{l=1}^\infty\Big(\lambda_{2l-1}^{1/2}\phi_{2l-1}(u),\lambda_{2l}^{1/2}\phi_{2l}(u)\Big)
	\bA^h
	\left (\lambda_{2l-1}^{1/2}\phi_{2l-1}(v),
		\lambda_{2l}^{1/2}\phi_{2l}(v)\right)^{\T}.
$$

Let
$
\boldsymbol{\Omega}(\theta) = \bI_{2}+\sum_{h=1}^\infty \big\{\exp(-ih\theta)+\exp(ih\theta)\big\}\bA^h.
$
We can calculate the spectral density function of $X_t(\cdot)$ at frequency $\theta \in [-\pi, \pi]$ by
$$
f_{X,\theta}(u,v) 
= \frac{1}{2\pi}\sum_{l=1}^\infty\Big(\lambda_{2l-1}^{1/2}\phi_{2l-1}(u),\lambda_{2l}^{1/2}\phi_{2l}(u)\Big)
\boldsymbol{\Omega}(\theta)
\left (\lambda_{2l-1}^{1/2}\phi_{2l-1}(v),\lambda_{2l}^{1/2}\phi_{2l}(v)\right)^{\T}.
$$
For each $\phi(u) = \sum_{l=1}^\infty x_{l} \phi_l(u)$ with $\boldsymbol{x} = (x_1,x_2,\ldots)^{\T}
\in \mathbb{R}^{\infty}$ and  $\sum_{l=1}^\infty x_{l}^2 < \infty$, we have
\begin{eqnarray*}
	&&\frac{2\pi\int_{[0,1]^2}  \phi(u) \phi(v)f_{X,\theta}(u,v)du dv}{\int_{[0,1]^2} \phi(u) \phi(v)\Sigma_0(u,v)du dv}\\
	&=& \frac{\sum_{l=1}^\infty \Big(x_{2l-1}\lambda_{2l-1}^{1/2},x_{2l}\lambda_{2l}^{1/2}\Big)\boldsymbol{\Omega}(\theta)
		\big(x_{2l-1}\lambda_{2l-1}^{1/2},
			x_{2l}\lambda_{2l}^{1/2}\big)^{\T}}{\sum_{l=1}^\infty \Big(x_{2l-1}^2\lambda_{2l-1}+x_{2l}^2\lambda_{2l}\Big)}.
\end{eqnarray*}
Since $(a+b)/(c+d) \le \max(a/c,b/d)$ for any positive real numbers $a,b,c,d,$ we have 
\begin{eqnarray*}
	\cM(\mF_{X})= \underset{\theta \in [-\pi,\pi],\boldsymbol{x} \in {\mathbb R}^\infty}{\mbox{esssup}}
	\sup_{l \ge 1}\frac{\Big(x_{2l-1}\lambda_{2l-1}^{1/2},x_{2l}\lambda_{2l}^{1/2}\Big)\boldsymbol{\Omega}(\theta)
		\left (
			x_{2l-1} \lambda_{2l-1}^{1/2},x_{2l} \lambda_{2l}^{1/2}\right)^{\T}}{x_{2l-1}^2\lambda_{2l-1}+x_{2l}^2\lambda_{2l}}.
\end{eqnarray*}
It is easy to see that $\boldsymbol{\Omega}(\theta)$ is positive definite and its maximum eigenvalue, $\lambda_{\max}(\boldsymbol \Omega),$ exists and is uniformly bounded over $\theta \in [-\pi,\pi].$ As a result, it follows that
$$
\cM(\mF_{X}) \le \lambda_{\max}(\boldsymbol{\Omega}) < \infty.
$$
It is also worth noting that the dimension of $\bA$ can be larger and we only require that all eigenvalues of $\bA$ are in $[-\delta',\delta']$ for some $\delta' \in (0,1)$ to guarantee the stationarity of $\{a_{tm}\}$. 

{\it Example~3}: It can be shown that for any $p>1$, if $\{X_{t1},t\in \mathbb{Z}\},\ldots, $  $\{X_{tp},t\in \mathbb{Z}\}$ are independent and $\sup_{1\le j\le p}\mathcal{M}(\mF_{X_{j}})<\infty$, then the functional stability measure of $\bX_t = (X_{t1},\ldots,X_{tp})^{\T}$ is
$$
\mathcal{M}(\mF_{\bX}) \le \sup_{1\le j\le p}\mathcal{M}(\mathcal{F}_{X_{j}})<\infty.
$$

{\it Example~4}: Consider a general case $\boldsymbol{\mathbf{Y}}_t(u) = \bA \bX_t(u)$ with $\bA \in {\mathbb R}^{p \times p}$ and $\mathcal{M}(\mF_{\bX}) <\infty$. We can easily obtain that
$$
\mathcal{M}(\mF_{\boldsymbol{\mathbf Y}}) \le \mathcal{M}(\mF_{\bX}) < \infty,
$$ 
which implies that linear transformation of $\bX_t$ does not increase the functional stability measure.	It is also worth noting that components of $\boldsymbol{\mathbf{Y}}_t$ can be dependent in this example.

{\it Example~5}: Consider a more general scenario, $\boldsymbol{\mathbf{Y}}_t(u) = \bA \bX_t(u) + \boldsymbol{\mathbf {\xi}}_t(u)$, where $\bA \in {\mathbb R}^{p \times r},$ $\bX_t(u)$ is a $r$-dimensional vector of Gaussian processes and $\bX_t(u),$  $\boldsymbol{\mathbf {\xi}}_s(u)$ are independent for all $t$ and $s$. When $r$ is fixed, it implies that $\boldsymbol{\mathbf{Y}}_t(u)$ can be expressed under a factor model structure. Note that $(a + b)/(c +d) \le \max(a/c,b/d)$ for all positive real numbers $a,b,c$ and $d.$
Hence if $\max\{\mathcal{M}(\mF_{\bX}),\mathcal{M}(\mF_{\boldsymbol{\xi}}) \} <\infty,$ then 
$$\mathcal{M}(\mF_{\boldsymbol{\mathbf Y}}) \le \max\{\mathcal{M}(\mF_{\bX}),\mathcal{M}(\mF_{\boldsymbol{\xi}}) \} <\infty.
$$


\linespread{1.2}\selectfont
\bibliography{paperbib}
\bibliographystyle{dcu}

\newpage
\setcounter{page}{1}
\setcounter{section}{1}
\setcounter{equation}{0}
\renewcommand{\theequation}{B.\arabic{equation}}

\begin{center}
	{\noindent \bf \large Supplementary Material to ``On Consistency and Sparsity for High-Dimensional Functional Time Series with Application to Autoregressions"}\\
\end{center}
\begin{center}
	{\noindent Shaojun Guo and Xinghao Qiao}
\end{center}

\linespread{1.4}\selectfont

This supplementary material contains proofs of main theorems in Appendix~\ref{asec.pf}, additional technical proofs in Appendix~\ref{ap.tech.proof}, 
derivations of functional stability measure for the illustrative VFAR(1) example in Appendix~\ref{ap.ill.ex}, some derivations for VFAR models in Appendix~\ref{ap.deriv}, details of the algorithms to fit sparse VFAR models in Appendix~\ref{ap.alg} and additional empirical results in Appendix~\ref{ap.emp}.

\section{Proofs of main theorems}
\label{asec.pf}

\subsection{Proof of Theorem \ref{thm_con_op}}
\label{asec.pf.thm1}
(i) Define $\bY = \big(\langle \bPhi_1,\bX_1\rangle_{\cH}, \ldots, \langle \bPhi_1,\bX_n\rangle_{\cH}\big)^\T$, then $\bY \sim N(\0, \bQ)$, where $Q_{rs} = \big\langle \bPhi_1, \bXi_{r-s}(\bPhi_1)\big\rangle_{\cH}$ for $r,s =1,\ldots,n$.
Note $\big\langle\bPhi_1, \widehat{\bXi}_0(\bPhi_1)\big\rangle_{\cH} = {n^{-1}} \bZ^{\T} \bQ \bZ$ with $\bZ \sim N(\0,\bI_{n})$ and $\big\langle\bPhi_1, {\bXi}_0(\bPhi_1)\big\rangle_{\cH} = \tE\big(n^{-1}\bZ^{\T} \bQ \bZ\big).$
By the Hanson-Wright inequality in 
{\color{blue} Rudelson and Vershynin (2013)}, 
$$
P\left\{\Big|\big\langle \bPhi_1, \big(\widehat{\bXi}_0- \bXi_0\big)(\bPhi_1)\big\rangle_{\cH}\Big|>\epsilon\right\} \le 2 \exp\left\{- c \min\left( \frac{n^2 \epsilon^2}{\|\bQ\|_F^2}, \frac{n\epsilon}{\|\bQ\|}\right)\right\}
$$
for some constant $c> 0.$ By $\|\bQ\|_F^2/n \le \|\bQ\|^2$ and letting $\epsilon = \eta\|\bQ\|$, we obtain that
\begin{equation}
\label{coneq_Phi_1}
P\left\{\Big|\big\langle \bPhi_1, \big(\widehat{\bXi}_0- \bXi_0\big)(\bPhi_1)\big\rangle_{\cH}\Big|> \eta\|\bQ\|\right\} \le 2 \exp\left\{- c n\min\left( \eta^2,\eta\right)\right\}
\end{equation}
for some universal constant $c> 0$.

Next we derive an upper bound on the operator norm $\|\bQ\|.$ Specifically, for any $\bw = (w_1,\ldots,w_n)^{\T} \in \eR^n$ with $ ||\bw||=1,$ define $G_{\bw}(\theta) = \sum_{r =1 }^n w_r \exp(-i r\theta)$ and its conjugate by $G^*_{\bw}(\theta).$ Then we obtain that
\begin{eqnarray*}
	\bw^{\T}\bQ\bw &=& \sum_{r = 1}^n \sum_{s =1 }^n w_r w_s \big\langle \bPhi_1, \bXi_{r-s}(\bPhi_1)\big\rangle_{\cH}\\
	& =& \sum_{r = 1}^n \sum_{s =1 }^n w_r w_s \int_{-\pi}^{\pi}\big\langle\bPhi_1, \mF_{\bX,\theta}(\bPhi_1)\big\rangle_{\cH} \exp \{i(r-s)\theta\} d\theta \\
	&=& \int_{-\pi}^{\pi}\big\langle\bPhi_1, \mF_{\bX,\theta}(\bPhi_1)
	\big\rangle_{\cH} G_{\bw}(\theta) G^{*}_{\bw}(\theta)d\theta,
\end{eqnarray*}
where the second line follows from the inversion formula (\ref{inv.formula}).
For a fixed $\bPhi\in \cH$, denote
$
\cM (\mF_{\bX},{\bPhi}) = 2\pi \cdot {\text{ess} \sup}_{\theta\in [-\pi, \pi] }\big|\langle \bPhi,\mF_{\bX,\theta}(\bPhi) \rangle_{\cH}\big|.
$
Since $\big\langle\bPhi_1, \mF_{\bX,\theta}(\bPhi_1)\big\rangle_{\cH}$ is Hermitian and $\int_{-\pi}^{\pi}G_{\bw}(\theta) G^{*}_{\bw}(\theta)d\theta = 2\pi$, we have
$
\|\bQ\| \le \cM(\mF_{\bX},\bPhi_1).
$
Then it follows from (\ref{eq.sub.fsm}) that
$$
\|\bQ\| \le \cM(\mF_{\bX},\bPhi_1) \le \cM_k(\mF_{\bX}) \big\langle \bPhi_1, \bXi_0(\bPhi_1)\big\rangle_{\cH}.
$$
This result, together with (\ref{coneq_Phi_1}) implies (\ref{thm_Phi_1}).

(ii) Note that
$$
	4\big\langle \bPhi_1, (\widehat{\bXi}_0 - \bXi_0)(\bPhi_2)\big\rangle_{\cH} \leq 
	\big\langle \widetilde{\bPhi}_1, 
	(\widehat{\bXi}_0 -
	\bXi_0)(\widetilde{\bPhi}_1)\big\rangle_{\cH}
	 - \big\langle \widetilde{\bPhi}_2, (\widehat{\bXi}_0 - \bXi_0)(\widetilde{\bPhi}_2)\big\rangle_{\cH}.
$$
where $\widetilde{\bPhi}_1 = \bPhi_1 + \bPhi_2,$ $\widetilde{\bPhi}_2 = \bPhi_1 - \bPhi_2$
and
$\cM(\mF_{\bX},\widetilde{\bPhi}_i) \le 2\{\cM(\mF_{\bX},\bPhi_1)+ \cM(\mF_{\bX}, \bPhi_2) \}$ for $i=1, 2.$
Combing these with results in (i) leads to
\begin{eqnarray*}
&&P\Big[\Big|\big\langle \bPhi_1, (\widehat{\bXi}_0 - \bXi_0)(\bPhi_2)\big\rangle_{\cH}\Big|> \{\cM(\mF_{\bX},\bPhi_1)+ \cM(\mF_{\bX}, \bPhi_2)\}\eta\Big]\\
&\leq& \sum_{i=1}^2P\Big[\Big|\big\langle \widetilde\bPhi_i, (\widehat{\bXi}_0 - \bXi_0)(\widetilde\bPhi_i)\big\rangle_{\cH}\Big|> \cM(\mF_{\bX},\widetilde\bPhi_i)\eta\Big]\le 4 \exp\Big\{- c n\min\left( \eta^2,\eta\right)\Big\}
\end{eqnarray*}
for some universal constant $c> 0.$ This, together with, $\cM(\mF_{\bX},\bPhi_i) \le \cM_k(\mF_{\bX}) \big\langle \bPhi_i, \bXi_0(\bPhi_i)\big\rangle_{\cH}$ for $i=1, 2,$ implies (\ref{thm_Phi_12}), which completes the proof.  $\square$

\subsection{Proof of Theorem~\ref{thm_max_bound}}
\label{asec.pf.thm2}

First, we derive the concentration bound on $\|\widehat \Sigma_{jk}^{(0)} - \Sigma_{jk}^{(0)}\|_{\cS}$ for each $j$ and $k$.
Let 
$ \Delta_{jklm} = (\lambda_{jl} \lambda_{km})^{-1/2}\big\langle  \phi_{jl},(\widehat \Sigma_{jk}^{(0)} - \Sigma_{jk}^{(0)})(\phi_{km}) \big \rangle$
 for  $j, k= 1,\ldots,p,$ and $l,m = 1,\ldots, \infty.$ 
Then we have that
$
\big\|\widehat{\Sigma}_{jk}^{(0)} - \Sigma_{jk}^{(0)} \big\|_{\mathcal{S}}^2
= \sum^\infty_{l,m =1} \lambda_{jl}\lambda_{km} \Delta_{jklm}^2.
$
By Jensen's inequality, we have that
\begin{equation}
\tE\Big\{ \big\|\widehat{\Sigma}_{jk}^{(0)} - \Sigma_{jk}^{(0)} \big\|_{\mathcal{S}}^{2q}\Big\}
\le  \Big( \sum_{l,m= 1}^{\infty} \lambda_{jl} \lambda_{km}\Big)^ {q-1} \sum_{l,m = 1}^{\infty}
\lambda_{jl} \lambda_{km} \tE\big|\Delta_{jklm}  \big|^{2q}\le  \lambda_0^{2q} \sup_{l,m} \tE\big|\Delta_{jklm}  \big|^{2q}.
\label{Sigma.bd}
\end{equation}

For any given $(j, k,l,m),$ let $$
\bPhi_{1} = (0, \ldots,0,\lambda_{jl}^{-1/2}\phi_{jl}, 0,\ldots,0)^{\T} \text{ and } \bPhi_{2} = (0, \ldots,0,\lambda_{km}^{-1/2}\phi_{km}, 0,\ldots,0)^{\T}.
$$
By the definition of $\Delta_{jklm}$ and orthonormality of $\{\phi_{jl}(\cdot)\}$ and $\{\phi_{km}(\cdot)\}$ for each $j,k=1,\dots,p,$ we have
$
\Delta_{jklm}  = \langle\bPhi_{1}, (\widehat{\bXi}_0 - \bXi_0)(\bPhi_{2})\big\rangle_{\cH},
$
$\langle \bPhi_{1},\bXi_0(\bPhi_{1})\rangle_{\cH}=\langle \bPhi_{2},\bXi_0(\bPhi_{2})\rangle_{\cH}=1.$ Applying (\ref{thm_Phi_12}) in Theorem~\ref{thm_con_op}, we can obtain that
	\begin{equation}
	\label{con_Delta}
	P\Big\{ \big|\Delta_{jklm}  \big| > 2 \cM_1(\mF_{\bX}) \eta
	\Big\} \le 4 \exp\Big\{- c n \min(\eta^2, \eta)\Big\},
	\end{equation}
for $j, k= 1,\ldots,p,$ $l = 1,\ldots,d_j$ and $m = 1,\ldots, d_k.$ It then follows from Lemma~\ref{lemma.moment} in Appendix~\ref{ap.tech.proof} in the Supplementary Material that for each integer $q \ge 1$,
$$
\left\{2 \cM_1(\mF_{\bX})  \right\}^{-2q} \tE\big|\Delta_{jklm}  \big|^{2q} \le q! 4(4c^{-1}n^{-1})^q + 4(2q)!(4c^{-1}n^{-1})^{2q}.
$$
This together with (\ref{Sigma.bd}) implies that
\begin{equation}
\label{Sigma_moment}
	\left(2 \cM_1(\mF_{\bX}) \lambda_0 \right)^{-2q} \tE\Big\{\big\|\widehat{\Sigma}_{jk}^{(0)} - \Sigma_{jk}^{(0)} \big\|_{\mathcal{S}}^{2q}\Big\}  \le  q! 4(4c^{-1}n^{-1})^q + (2q)!4(4c^{-1}n^{-1})^{2q}.
\end{equation}
Finally, it follows from Lemma~\ref{lemma.moment} that there exists some universal constant $\tilde c > 0 $ such that
\begin{eqnarray*}
	P\Big \{\big\|\widehat{\Sigma}_{jk}^{(0)} - \Sigma_{jk}^{(0)} \big\|_{\mathcal{S}} \ge 2\cM_1(\mF_{\bX}) \lambda_0 \eta   \Big \}\le 4 \exp\Big\{- \tilde c n \min(\eta^2, \eta)\Big\}.
\end{eqnarray*}

Using the definition of
$\|\widehat \bSigma_0 - \bSigma_0\|_{\max} = {\max}_{1 \le j,k \le p} \|\widehat \Sigma_{jk}^{(0)} - \Sigma_{jk}^{(0)}\|_{\cS}$ and applying the union bound of probability, we obtain that
\begin{eqnarray*}
	P\Big \{ \|\widehat \bSigma_0 - \bSigma_0\|_{\max} \ge 2 \cM_1(\mF_{\bX}) \lambda_0 \eta   \Big \}\le 4 p^2\cdot\exp\Big\{ - \tilde c n \min(\eta^2, \eta)\Big\}.
\end{eqnarray*}

Let $\eta = \rho \sqrt{\log p /n} \leq 1$ and $ \rho^2 \tilde c > 2$, which can be achieved for sufficiently large $n$.  We obtain that
\begin{eqnarray*}
	P\left \{ \|\widehat \bSigma_0 - \bSigma_0\|_{\max} \ge 2 \cM_1(\mF_{\bX}) \lambda_0 \rho \sqrt{\frac{\log p}{n}}   \right \}\le 4 p^{2 - \tilde c\rho^2}. 
\end{eqnarray*}
The proof is complete. $\square$



\subsection{Proof of Theorem~\ref{thm_lam_phi}}
\label{ap_pf_lam_phi}
To simplify our notation, for each $j,k=1,\dots,p,$ we will denote $\Sigma_{jk}^{(0)}$ and $\widehat\Sigma_{jk}^{(0)}$ by $\Sigma_{jk}$ and $\widehat\Sigma_{jk}, $ respectively, in our subsequent proofs. Let $\delta_{jl} = {\min}_{1 \le k \le l}\{\lambda_{jk} - \lambda_{j(k+1)}\}$ and $\widehat \Delta_{jk} = \widehat \Sigma_{jk} - \Sigma_{jk}$ for $j,k=1, \dots,p$ and $l=1,2 \dots.$ It follows from (4.43) and Lemma~4.3 of \cite{Bbosq1} that
\begin{equation}
\label{eigen.bd}
\underset{l \ge 1}{\sup}~ |\widehat \lambda_{jl} - \lambda_{jl}| \le \|\widehat \Delta_{jj}\|_{\cS}~~\mbox{and}~~\underset{l \ge 1}{\sup} ~\delta_{jl}\|\widehat \phi_{jl} - \phi_{jl}\| \le 2\sqrt{2} \|\widehat \Delta_{jj}\|_{\cS}.
\end{equation}
Moreover, we can express $\widehat \lambda_{jl} - \lambda_{jl}$ and $\widehat \phi_{jl} - \phi_{jl},$ as stated in Lemma~\ref{lemma.lambda.phi} in 
Appendix~\ref{ap.tech.proof}. The proof of Theorem~\ref{thm_lam_phi} relies on the concentration inequalities for eigenvalues and eigenvectors as stated in the following Lemmas~\ref{lemma_lambda} and \ref{lemma_phi}, whose proofs are provided in Appendix~\ref{ap.tech.proof}. 
\begin{lemma}
    \label{lemma_lambda}
Suppose that Conditions~\ref{cond.cov.func}--\ref{cond_eigen} hold. Then there exists some universal constant $\tilde c_1>0$ such that for each $j=1,\dots,p, l=1, \dots,\infty,$ and any $\eta >0$,
\begin{equation}
\label{thm_con_lambda}
P\left\{\left |\frac{\widehat \lambda_{jl} - \lambda_{jl}}{\lambda_{jl}}\right|  >  \cM_1(\mF_{\bX}) \eta  + \rho_1 l^{2\alpha + 1} \cM_{1}^2(\mF_{\bX}) \eta^2\right\} \le 4 \exp\Big\{- \tilde c_1 n \min(\eta^2,\eta)\Big\},
\end{equation}
where $\rho_1 = 16\sqrt{2} c_0^{-2} \alpha \lambda_0^2.$ 
\end{lemma}

\begin{lemma}
    \label{lemma_phi}
Suppose that Conditions~\ref{cond.cov.func}--\ref{cond_eigen} hold. Then there exists some universal constant $\tilde c>0$ such that for each $j=1,\dots,p, l=1, \dots,\infty,$ and any $\eta > 0,$
	\begin{equation}
	\label{thm_con_phi}
	P\left\{\Big\|\widehat \phi_{jl} - \phi_{jl}\Big\| > 4 \sqrt{2}  \cM_1(\mF_{\bX}) \lambda_0 c_0^{-1} l^{\alpha + 1}  \eta\right\}
	\le  4 \exp\Big\{ - \tilde c n \min(\eta^2, \eta)\Big\}.
	\end{equation}
	
\end{lemma}
{\it Proof of Theorem~\ref{thm_lam_phi}}.
Applying the union bound of probability in (\ref{thm_con_lambda}), we obtain that
$$
P\left\{\underset{1 \leq j\leq p, 1\leq l \leq M}{\max}\left |\frac{\widehat \lambda_{jl} - \lambda_{jl}}{\lambda_{jl}}\right|  >  \cM_1(\mF_{\bX}) \eta  + \rho_1 l^{2\alpha + 1} \cM_{1}^2(\mF_{\bX}) \eta^2\right\} \le 4pM \exp\Big\{- \tilde c_1 n \min(\eta^2,\eta)\Big\}.
$$
Let $\eta = \tilde\rho_2 \sqrt{\log (pM)/n} \leq 1$ and $1+\rho_1 M^{2\alpha+1}\cM_{1}(\mF_{\bX})\eta\leq \tilde\rho_1$, 
which can be achieved for sufficiently large $n \gtrsim M^{4\alpha+2}\cM_{1}^2(\mF_{\bX})\log(pM).$
We obtain that
\begin{equation}
\label{max_coneq_lam}
   P\left\{\underset{1 \leq j\leq p, 1\leq l \leq M}{\max}\left |\frac{\widehat \lambda_{jl} - \lambda_{jl}}{\lambda_{jl}}\right|  >  \tilde\rho_1\tilde\rho_2\cM_1(\mF_{\bX}) \sqrt{\frac{\log (pM)}{n}} \right\} \le 4(pM)^{1-\tilde c_1\widetilde\rho_2^2}. 
\end{equation}
Finally, letting $\eta = \tilde\rho_3 \sqrt{\frac{\log (pM)}{n}}<1$ and following the same developments, we obtain that 
\begin{equation}
\label{max_coneq_phi}
	P\left\{\underset{1 \leq j\leq p, 1\leq l \leq M}{\max}\left(\frac{\big\|\widehat \phi_{jl} - \phi_{jl}\big\|}{l^{\alpha+1}}\right) > 4 \sqrt{2}\lambda_0 c_0^{-1} \widetilde\rho_3 \cM_1(\mF_{\bX})  \sqrt{\frac{\log (pM)}{n}} \right\}\le 4(pM)^{1-\tilde c\widetilde\rho_3^2}.
\end{equation}
It follows from (\ref{max_coneq_lam}) and (\ref{max_coneq_phi}) that, for sufficiently large $n \gtrsim M^{4\alpha+2}\cM_{1}^2(\mF_{\bX})\log(pM)$ and suitable choices of constants $c_1, c_2>0$, (\ref{bd_max_lam_phi}) holds. 
$\square$

\subsection{Proof of Theorem~\ref{thm_cov}}
\label{ap.pf.thm_cov}
The proof of Theorem~\ref{thm_cov} is based on the following Lemma~\ref{lemma_score}.
\begin{lemma}
    \label{lemma_score}
	Suppose that Conditions~\ref{cond.cov.func}--\ref{cond_eigen} hold. Then there exist some positive constants $\rho_4, \rho_5, \tilde c_3$ and $\tilde c_4$ such that
	
	(i) for each $j = 1,\ldots,p,$ $l=1,\dots, \infty,$ and any $\eta >0,$
		\begin{equation}
		\label{con_bd_score1}
		P\left\{\left |\frac{\widehat \sigma_{jjll}^{(0)} - \sigma_{jjll}^{(0)}}{\lambda_{jl}}\right|  >  \cM_1(\mF_{\bX}) \eta  + \rho_1 \cM_{1}^2(\mF_{\bX}) l^{2\alpha + 1}\eta^2\right\} \le 4 \exp\Big\{- \tilde c_3n \min(\eta^2,\eta)\Big\};
		\end{equation}
		
		(ii) for each $j,k = 1,\ldots,p,$ $l,m=1,\dots,\infty,$ but $j \neq k$ or $l \neq m,$ a fixed $h$ and any $\eta >0$,
		\begin{equation}
	    \label{con_bd_score2}
		\begin{split}
		&P\left\{ \left| \frac{\widehat \sigma_{jklm}^{(h)} - \sigma_{jklm}^{(h)}}{\lambda_{jl}^{1/2} \lambda_{km}^{1/2}}\right| \ge \rho_4 \cM_1(\mF_{\bX}) (l \vee m)^{\alpha + 1}  \eta
		+ \rho_5 \cM_1^2(\mF_{\bX})  (l \vee m)^{3\alpha + 2}\eta^2\right\} \\
		&\le \tilde c_4 \exp\big\{ - \tilde c_3 n\min(\eta^2, \eta)\big\} + \tilde c_4\exp\big\{ - \tilde c_3  \cM_1^{-2}(\mF_{\bX}) n (l \vee m)^{-2(\alpha + 1)}\big\}.
		\end{split}
		\end{equation}
\end{lemma}
{\it Proof of Lemma~\ref{lemma_score}.} For the special case of $(j,k,l,m)$ with $j=k,$ provided that $\widehat \sigma_{jjlm}^{(0)} = \widehat \lambda_{jl}I(l = m)$ and $\sigma_{jjlm}^{(0)} = \lambda_{jl}I(l = m),$ 
(\ref{con_bd_score1}) follows directly from Lemma~\ref{lemma_lambda}. 

For general cases of $(j,k,l,m)$ with $j \neq k,$
$
\widehat{\sigma}_{jklm}^{(h)} = (n-h)^{-1}\sum_{t = 1}^{n-h} \widehat{\xi}_{tjl} \widehat{\xi}_{(t+h)km}
$
and  $\sigma_{jklm}^{(h)}= E(\xi_{tjk} \xi_{(t+h)lm}).$
Let $\widehat r_{jl} = \widehat \phi_{jl} - \phi_{jl},$ then $\widehat{\sigma}_{h,jklm} - \sigma_{h,jklm}$ can be decomposed as
\begin{eqnarray*}
\label{thm_cov_1}
\widehat{\sigma}_{jklm}^{(h)} - \sigma_{jklm}^{(h)} 
& = & \big\langle \widehat r_{jl}, \widehat \Sigma_{jk}^{(h)} (\widehat r_{km}) \big \rangle +
 \Big(\big\langle \widehat r_{jl}, \widehat \Delta_{jk}^{(h)} (\phi_{km}) \big \rangle +
 \big\langle \phi_{jl}, \widehat \Delta_{jk}^{(h)} (\widehat r_{km}) \big \rangle \Big)\\
 &&  + \Big(\big\langle \widehat r_{jl}, \Sigma_{jk}^{(h)} (\phi_{km})\big \rangle
 + \big\langle \phi_{jl}, \Sigma_{jk}^{(h)} (\widehat r_{km}) \big \rangle \Big) +
\big\langle \phi_{jl}, \widehat \Delta_{jk}^{(h)}(\phi_{km}) \big \rangle \\
	& = & I_1 + I_2 + I_3 + I_4.
\end{eqnarray*}

For a fixed $h > 0,$ let $\Omega_{jk}^{(h)} = \Big\{ \|\widehat \Delta_{jk}^{(h)}\|_{\cS} \le \lambda_0\Big\}$ and
$
\widetilde\Omega_{jk}^{(h)} = \Big \{ \|\widehat \Delta_{jk}^{(h)}\|_{\cS} \le 4 \cM_1(\mF_{\bX}) \lambda_0 \eta\Big\}.
$
It follows from the same developments as in the proof of (\ref{thm_Sigma_compt}) in Theorem~\ref{thm_max_bound} and Proposition~\ref{thm_sigma_k} that there exists some universal constant $\tilde c>0$ such that for any $\eta>0,$ a fixed $h \neq 0$ and each $j,k=1, \dots, p,$
	\begin{equation}
	\label{thm_Sigma_compt_h}
	P\left\{\big\|\widehat{\Sigma}_{jk}^{(h)} - \Sigma_{jk}^{(h)}\big\|_{\cS} >   4 \cM_1(f_{\bX})\lambda_0 \eta\right\} \le 8 \exp\Big\{- \tilde c n \min(\eta^2, \eta)\Big\}.
	\end{equation}

On the event $\Omega_{jk}^{(h)} \cap \widetilde \Omega_{jj}^{(0)} \cap \widetilde \Omega_{kk}^{(0)} \cap \widetilde \Omega_{jk}^{(h)},$ it follows from Condition~\ref{cond_eigen} with $\lambda_{jl} \geq c_0 \alpha^{-1} l^{-\alpha}$, (\ref{eigen.bd}), Lemma~\ref{lemma.sub.ieq} in Appendix~\ref{ap.tech.proof} that
\begin{equation}
\label{thm_cov_2}
\begin{split}
\left |\frac{I_1}{\lambda_{jl}^{1/2}\lambda_{km}^{1/2}} \right|
& 
\lesssim (lm)^{\alpha/2}\|\widehat r_{jl}\| \big(\|\widehat \Delta_{jk}^{(h)}\| + \|\Sigma_{jk}^{(h)}\|_{\cS}\big)\|\widehat r_{km}\|\\
&
\lesssim \cM_{1}^2(\mF_{\bX})(l\vee m)^{3\alpha + 2} \eta^2,
\end{split}
\end{equation}
\begin{equation}
\label{thm_cov_3}
\begin{split}
\left|\frac{I_{2}}{\lambda_{jl}^{1/2} \lambda_{km}^{1/2}}\right|
&
\lesssim(lm)^{\alpha/2}\|\widehat \Delta_{jk}^{(h)}\|_{\cS}\Big( l^{\alpha + 1}\|\widehat \Delta_{jj}^{(0)}\|_{\cS} + m^{\alpha + 1}\|\widehat \Delta_{kk}^{(0)}\|_{\cS}\Big) \\
&
\lesssim \cM_1^2(\mF_{\bX})(l \vee m)^{2\alpha + 1}  \eta^2.
\end{split}
\end{equation}
For the term $I_4$, it follows from (\ref{thm_sig_k2}) in Proposition~\ref{thm_sigma_k} and the fact $\lambda_{jl}+\lambda_{jm} \geq 2\lambda_{jl}^{1/2}\lambda_{jm}^{1/2}$ that 
\begin{equation}
\label{thm_cov_4}
P\left\{ \left|\frac{I_4}{\lambda_{jl}^{1/2}\lambda_{km}^{1/2}}\right| \ge 4 \cM_1(\mF_{\bX}) \lambda_0 \eta\right \} \le 8 \exp\Big\{ - c n \min(\eta^2, \eta)\Big\}.
\end{equation}
Finally, we consider the term $I_{3}$.
By Lemma~\ref{lemma.sub.ieq}, we have that $\|\Sigma_{jk}^{(h)}(\phi_{km})\| \le \lambda_{km}^{1/2} \lambda_0^{1/2}$
and $\|\Sigma_{jk}^{(h)}(\phi_{jl})\| \le\lambda_{jl}^{1/2} \lambda_0^{1/2}.$
These results together with (\ref{thm_con_phi_g}) in Lemma~\ref{lemma_phi_g} in Appendix~\ref{ap.tech.proof} and Condition~\ref{cond_eigen} with $\lambda_{jl} \geq c_0 \alpha^{-1} l^{-\alpha}$ imply that
\begin{equation}
\label{thm_cov_5}
\left|\frac{I_3}{\lambda_{jl}^{1/2} \lambda_{km}^{1/2}} \right| 
\lesssim \cM_1(\mF_{\bX}) (l \vee m)^{\alpha + 1}  \eta
+  \cM_1^2(\mF_{\bX})  (l \vee m)^{(5\alpha + 4)/2} \eta^2
\end{equation}
holds with probability greater than $1-16 \exp\big\{ - \tilde c_2n \min(\eta^2, \eta)\big\} -
8\exp\big\{ - \tilde c_2  \cM_1^{-2}(\mF_{\bX}) n (l \vee m)^{-2(\alpha + 1)}\big\},
$
with some positive constant $\tilde c_2.$

Combining (\ref{thm_Sigma_compt_h})--(\ref{thm_cov_5}) and by Theorem~\ref{thm_max_bound}, we obtain that there exist
four positive constants $\rho_4$, $\rho_5,$ $\tilde c_3,$ $\tilde c_4$ such that (\ref{con_bd_score2}) in Lemma~\ref{lemma_score} holds. For the case of $h=0,$ we follow the same developments as above by applying Theorems~\ref{thm_con_op}--\ref{thm_max_bound} and hence  (\ref{con_bd_score2}) follows with the different choice of relevant positive constants. The proof of Lemma~\ref{lemma_score} is complete. $\square$

{\it Proof of Theorem~\ref{thm_cov}}.
Let $\eta=\tilde \rho_4 \sqrt{\log(pM)/n}<1$ and
$\rho_4+\rho_5 M^{2\alpha+1} \cM_1(\mF_{\bX}) \eta \leq \tilde \rho_5,$
which can be achieved for sufficiently large $n \gtrsim M^{4\alpha+2}\cM_{1}^2(\mF_{\bX})\log(pM).$
Following the similar techniques as used in the proof of (\ref{bd_max_lam_phi}) in Theorem~\ref{thm_lam_phi}, we can obtain (\ref{bd_score_max}), which completes the proof. $\square$

\subsection{Proof of Theorem~\ref{thm.vfar}}
Since $\widehat \bB_{j} \in \eR^{pq \times q}$ is the minimizer of (\ref{vfar.crit.theory}), we have
$$
-\langle\langle \widehat \bY_j, \widehat \bB_{j}\rangle\rangle +
\frac{1}{2} \langle\langle \widehat \bB_{j}, \widehat \bGamma \widehat\bB_{j}\rangle\rangle
 + {\gamma_{nj}}  \|\widehat\bB_{j}\|_{1}^{(q)}
 \le -\langle\langle \widehat \bY_j, \bB_{j}\rangle\rangle +
 \frac{1}{2} \langle\langle \bB_{j}, \widehat \bGamma \widehat\bB_{j}\rangle\rangle
 + {\gamma_{nj}}  \|\bB_{j}\|_{1}^{(q)}.
$$
Letting $\bDelta_j = \widehat \bB_{j} - \bB_{j}$ and  $S_j^{c}$ be the complement of $S_j$ in the set $\{1, \dots,p\},$ we have
\begin{eqnarray*}
	\frac{1}{2}\langle\langle \bDelta_j, \widehat \bGamma \bDelta_j \rrangle
	&\leq& \llangle \bDelta_j, \widehat \bY_j - \widehat \bGamma \bB_{j}\rrangle +
	{\gamma_{nj}}  \Big( \|\bB_{j}\|_{1}^{(q)} - \|\bB_{1j} + \bDelta_j\|_{1}^{(q)}\Big)\\
	&\leq & \llangle \bDelta_j, \widehat \bY_j - \widehat \bGamma \bB_{j}\rrangle +
	{\gamma_{nj}} \Big( \|\bB_{jS_j}\|_{1}^{(q)}- \|\bB_{jS_j} + \bDelta_{jS_j}\|_{1}^{(q)} - \|\bDelta_{jS_j^c}\|_{1}^{(q)}\Big)\\
	&\leq& \llangle \bDelta_j, \widehat \bY_j - \widehat \bGamma \bB_j\rrangle +
	{\gamma_{nj}} \Big( \|\bDelta_{jS_j}\|_{1}^{(q)} - \|\bDelta_{jS_j^c}\|_{1}^{(q)}\Big)
\end{eqnarray*}
By Lemma~\ref{lemma.norm.ieq} in Appendix~\ref{ap.tech.proof}, Condition~\ref{cond.fvar.max.error} and the choice of $\gamma_{nj},$ we have
$$
\big|\llangle \bDelta_j, \widehat \bY_j - \widehat \bGamma \bB_j\rrangle \big|
\le \|\widehat \bY_j - \widehat \bGamma \bB_j\|_{\max}^{(q)} \|\bDelta_j\|_{1}^{(q)} \leq \frac{\gamma_{nj}}{2} \big( \|\bDelta_{jS_j}\|_{1}^{(q)} + \|\bDelta_{jS_j^c}\|_{1}^{(q)}\big).
$$
Combing the above two results, we have
\begin{eqnarray*}
	0 \leq \frac{1}{2}\langle\langle \bDelta_j, \widehat \bGamma \bDelta_j \rrangle \le
	\frac{3\gamma_{nj}}{2} \|\bDelta_{jS_j}\|_{1}^{(q)}   - \frac{\gamma_{nj}}{2} \|\bDelta_{jS_j^c}\|_{1}^{(q)},
\end{eqnarray*}
which implies $\|\bDelta_{jS_j^c}\|_{1}^{(q)} \le 3\|\bDelta_{jS_j}\|_{1}^{(q)}$ and therefore $\|\bDelta_j\|_{1}^{(q)} \le 4\|\bDelta_{jS_j}\|_{1}^{(q)} \le 4 \sqrt{s_j} \|\bDelta_{j}\|_{F}.$ 
This result together with Condition~\ref{cond.fvar.RE} and $\tau_2\ge  32\tau_1 q^2 s_j$ implies that \begin{equation}
\langle\langle \bDelta_j, \widehat \bGamma \bDelta_j \rrangle
\ge \tau_2 \|\bDelta_j\|_F^2 - \tau_1 q^2\big\{\|\bDelta_j\|_{1}^{(q)}\big\}^2
\ge \big(\tau_2 - 16\tau_1 q^2 s_j\big)\|\bDelta_j\|_F^2 \ge \frac{\tau_2}{2}\|\bDelta_j\|_F^2.
\end{equation}
Therefore,
$$
\frac{\tau_2}{4}\|\bDelta_j\|_F^2 \le \frac{3}{2} \gamma_{nj} \|\bDelta_j\|_{1}^{(q)} \le 6 \gamma_{nj}s_j^{1/2} \|\bDelta_{j}\|_{F},
$$
which implies that
\begin{equation}
\label{err.B}
\|\bDelta_j\|_F \le \frac{24s_j^{1/2}\gamma_{nj}}{\tau_2} \text{ and  } \|\bDelta_j\|_{1}^{(q)} \le \frac{96 s_j\gamma_{nj}}{\tau_2},
\end{equation}
as is claimed in Theorem~\ref{thm.vfar}.

Next we prove the upper bound on $\widehat\bA - \bA.$

For $k\in S_j,$ it follows from  $\bPsi_{jk}=\int\int\bphi_k(v)A_{jk}(u,v)\bpsi_j(u)^{\T}dudv,$ Condition~\ref{cond.fvar.bias} with $A_{jk}(u,v)=\bphi_k(v)^{\T}\ba_{jk}\bphi_j(u) + (\sum_{l,m=1}^{\infty}-\sum_{l,m=1}^{q})a_{jklm}\phi_{jl}(u)\phi_{km}(v)$ and orthonormality of $\{\phi_{jl}(\cdot)\}_{l\geq 1}$ and $\{\phi_{km}(\cdot)\}_{m\geq 1}$ that  $\|\bPsi_{jk}\|_F = ||\ba_{jk}||_F=\big\{\sum_{l,m=1}^q \mu_{jk}^2 (l+m)^{-2\beta-1}\big\}^{1/2}\leq \big\{\mu_{jk}^2\int_{1}^q\int_{1}^q (x+y)^{-2\beta-1}dxdy\big\}^{1/2}=O(\mu_{jk}).$ For $k \in S_j^c,$ we have $\bPsi_{jk}={\bf 0}.$ Hence
\begin{equation}
    \label{bd.Psi}
    \|\bPsi_{j}\|_1^{(q)} = \sum_{k=1}^p ||\bPsi_{jk}||_F =O\big(\sum_{k \in S_j}\mu_{jk}\big) =O(s_j).
\end{equation}

Observe that
$
\widehat \bPsi_j-\bPsi_j=\widehat\bD^{-1}\widehat\bB_j-\bD^{-1}\bB_j=(\widehat \bD^{-1}-\bD^{-1})\bB_j + \bD^{-1}(\widehat\bB_j-\bB_j) + (\widehat \bD^{-1}-\bD^{-1})(\widehat\bB_j-\bB_j).
$
It follows from the diagonal structure of $\widehat\bD^{-1}$ and $\bD^{-1}$ that
\begin{equation}
\label{err.Psi0}
\begin{split}
\|\widehat \bPsi_{j} - \bPsi_j\|_1^{(q)}\le & \|(\widehat \bD^{-1} - \bD^{-1})\|_{\max} \|\bB_j\|_1^{(q)} + \|\bD^{-1}\|_{\max}\|\widehat \bB_j - \bB_j\|_{1}^{(q)}\\
&+\|(\widehat \bD^{-1} - \bD^{-1})\|_{\max} \|\widehat \bB_j - \bB_j\|_1^{(q)}.
\end{split}
\end{equation}
By Conditions~\ref{cond_eigen}, \ref{cond.fvar.eigen} and the fact $\widehat\bD_k=\text{diag}\big(\widehat\lambda_{k1}^{1/2}, \dots, \widehat\lambda_{kq}^{1/2}\big),$ 	$\bD_k=\text{diag}\big(\lambda_{k1}^{1/2}, \dots, \lambda_{kq}^{1/2}\big),$
we have $\|(\widehat \bD^{-1} - \bD^{-1})\|_{\max} \leq \alpha^{1/2}c_0^{-1/2}q^{\alpha/2}C_{\lambda}\cM(\mF_{\bX}) \sqrt{\frac{\log (pq)}{n}}$ and $\|\bD^{-1}\|_{\max} \le \alpha^{1/2}c_0^{-1/2} q^{\alpha/2}.$ By Condition~\ref{cond.cov.func} and (\ref{bd.Psi}), we have $||\bB||_1^{(q)}\leq ||\bD||_{\max}^{(q)}\|\bPsi_{j}\|_1^{(q)}= O(\lambda_0^{1/2}s_j).$ These results together with (\ref{err.B}) implies that
\begin{equation}
    \label{err.Psi}
    \|\widehat \bPsi_{j} - \bPsi_j\|_1^{(q)} \leq \frac{96\alpha^{1/2}q^{\alpha/2}s_j\gamma_{nj}}{c_0^{1/2}\tau_2}\Big\{1+o(1)\Big\},
\end{equation}
where the constant comes from the second term in (\ref{err.Psi0}), since the first and third terms are of smaller orders relative to the second term.

For each $j,k=1,\dots,p,$ note that
\begin{eqnarray*}	
\widehat \tA_{jk}(u,v) - \tA_{jk}(u,v)	&=& \widehat\bphi_k(v)^\T \widehat\bPsi_{jk} \widehat\bphi_j(u) -\bphi_k(v)^\T \bPsi_{jk} \bphi_j(u) + \tR_{jk}(u,v)\\
&=&\widehat \bphi_k(v)^{\T} \widehat \bPsi_{jk} \left\{\widehat \bphi_j(u) -\bphi_j(u) \right\} +
 \left\{\widehat \bphi_k(v)-\bphi_k(v)\right\}^{\T} \widehat \bPsi_{jk}\bphi_j(u)  \\
	& & +~ \bphi_k(v)^{\T} (\widehat \bPsi_{jk} -\bPsi_{jk})\bphi_j(u) + \tR_{jk}(u,v),
\end{eqnarray*}

We bound the first three terms. By Lemma~\ref{lemma.err.A} in Appendix~\ref{ap.tech.proof}, we have
\begin{eqnarray}
\label{error.Psi}
&&\left\|\widehat \bphi_k(v)^{\T} \widehat \bPsi_{jk} \left\{\widehat \bphi_j(u) -\bphi_j(u) \right\}\right\|_{\cS}\leq
 q^{1/2} \max_{1 \le l \le q}\|\widehat \phi_{jl} - \phi_{jl}\| \|\widehat \bPsi_{jk} \|_F, \nonumber \\
&&\left\|\left\{\widehat \bphi_k(v)-\bphi_k(v)\right\}^{\T} \widehat \bPsi_{jk}\bphi_j(u) \right\|_{\cS} \leq
q^{1/2} \max_{1 \le m \le q}\|\widehat \phi_{km} - \phi_{km}\| \|\widehat \bPsi_{jk} \|_F,\\
&&\left \| \bphi_k(v)^{\T} (\widehat \bPsi_{jk} - \bPsi_{jk})\bphi_j(u)\right\|_{\cS}= \|\widehat \bPsi_{jk} - \bPsi_{jk}\|_F. \nonumber
\end{eqnarray}

We then bound the fourth term. By $R_{jk}(u,v)=(\sum_{l,m=1}^{q}-\sum_{l,m=1}^{\infty})a_{jklm}\phi_{jl}(u)\phi_{km}(v),$ we have
\begin{eqnarray*}
||R_{jk}||_{\cS}^2 &=& O(1) \Big\|\sum_{l=q+1}^{\infty}\sum_{m=1}^{\infty}a_{jklm}\phi_{jl}(u)\phi_{km}(v)\Big\|_{\cS}^2\\
&=&O(1) \sum_{l=q+1}^{\infty}\sum_{m=1}^{\infty}a_{jklm}^2 \leq O(1)\mu_{jk}^2\sum_{l=q+1}^{\infty}\sum_{m=1}^{\infty}(l+m)^{-2\beta-1} = O(\mu_{jk}^2q^{-2\beta+1}).
\end{eqnarray*}
This together with Condition~\ref{cond.fvar.bias} implies that
\begin{equation}
    \label{err.R}
    \underset{1\leq j \leq p}{\max}\sum_{k=1}^p||R_{jk}||_{\cS}\leq O\big(q^{-\beta+1/2}\underset{1\leq j \leq p}{\max}\sum_{k \in S_j}\mu_{jk}\big)=O\big(sq^{-\beta+1/2}\big).
\end{equation}

It follows from (\ref{bd.Psi}), (\ref{err.Psi}), (\ref{error.Psi}), (\ref{err.R}) and the fact $\|\widehat \bPsi_{j} \|_1^{(q)} \le \|\widehat \bPsi_{j} - \bPsi_j\|_1^{(q)} + \| \bPsi_{j} \|_1^{(q)} =O(s_j)$ that
\begin{eqnarray*}
\|\widehat \bA - \bA\|_\infty &\le& 2q^{1/2} \max_{\underset{1 \le l \le q}{1 \le j \le p}}\|\widehat \phi_{jl} - \phi_{jl}\|\max_{1 \le j \le p}\|\widehat \bPsi_{j} \|_1^{(q)}
+ \max_{1 \le j \le p}\|\widehat \bPsi_{j}  - \bPsi_j\|_1^{(q)} + ||\bR||_{\infty}.\\
&\leq &  \frac{96\alpha^{1/2}q^{\alpha/2}s\gamma_{n}}{c_0^{1/2}\tau_2}\Big\{1+o(1)\Big\},
\end{eqnarray*}
where the constant comes from $\max_{j}\|\widehat \bPsi_{j}  - \bPsi_j\|_1^{(q)},$ since other terms are of smaller orders of this term. The proof is complete. $\square$

\renewcommand{\theequation}{C.\arabic{equation}}
\section{Additional technical proofs}
\label{ap.tech.proof}

\subsection{Proof of Proposition 1}
\label{pf.prop1}
Let $\bY_{1,t} = \bX_t + \bX_{t+h}$ , $\bSigma_{\bY_{1}, \ell} (u,v) = \cov\{\bY_{1,t}(u), \bY_{1,(t+ \ell)}(v)\}, \ell \in \eZ, (u,v) \in \cU^2.$ 
Define the spectral density operator of $\bY_{1,t}$ by
$$
\mF_{\bY_1,\theta} = \frac{1}{2 \pi} \sum_{\ell = -\infty}^{\infty} \bXi_{\bY_1,\ell}\exp(- i \ell \theta), ~\theta \in [-\pi, \pi]. 
$$
Then we can obtain that
$
\mF_{\bY_1,\theta} = \{2 + \exp(-ih\theta) + \exp( ih\theta)\}\mF_{\bX,\theta}.
$
Similarly, by letting ${\bY}_{2,t}(u) = \bX_t(u) - \bX_{t+h}(u)$,
$\bSigma_{\bY_2, \ell} (u,v) = \cov\{\bY_{2,t}(u), \bY_{2,(t+\ell)}(v)\},\ell \in \eZ (u,v), \in \cU^2,$ 
and $\mF_{\bY_2,\theta}$ be the spectral density operator of $\bY_{2},$ $\theta \in [-\pi, \pi],$ we have
$
\mF_{\bY_2,\theta} = \{2 - \exp(-ih\theta) - \exp( ih\theta)\}\mF_{\bX,\theta}. 
$ Note that
$$
4\big\langle \bPhi_1, (\widehat{\bXi}_h - \bXi_h)(\bPhi_1)\big\rangle_{\cH} = \big\langle \bPhi_1, (\widehat{\bXi}_{\bY_1,0} - \bXi_{\bY_1,0})(\bPhi_1)\big\rangle_{\cH} - \big\langle \bPhi_1, (\widehat{\bXi}_{{\bY}_2,0} - \bXi_{{\bY}_2,0})(\bPhi_1)\big\rangle_{\cH}
$$
and
$
\cM(\mF_{{\bY}_i},\bPhi_1)\le 4 \cM(\mF_{\bX},\bPhi_1)
$ for $i=1,2.$
Combing these with results in the proof of (\ref{thm_Phi_1}) leads to
\begin{eqnarray*}
	&&P\Big[\Big|\big\langle \bPhi_1, (\widehat{\bXi}_{h} - \bXi_{h})(\bPhi_1)\big\rangle_{\cH}\Big|> 2\cM(\mF_{\bX},\bPhi_1)\eta\Big]\\
	&\leq& \sum_{i=1}^2 P\Big[\Big|\big\langle \bPhi_1, (\widehat{\bXi}_{\bY_i,0} - \bXi_{\bY_i,0})(\bPhi_1)\big\rangle_{\cH}\Big|> \cM(\mF_{\bY_i},\bPhi_1)\eta\Big] \leq 4 \exp\Big\{- c n\min\left( \eta^2,\eta\right)\Big\},
\end{eqnarray*}
for some constant $c> 0.$ This result, together with, $\cM(\mF_{\bX},\bPhi_1) \le \cM_k(\mF_{\bX}) \big\langle \bPhi_1, \bXi_0(\bPhi_1)\big\rangle_{\cH}$ implies (\ref{thm_sig_k1}).

Note that
$$
	4\big\langle \bPhi_1, (\widehat{\bXi}_h - \bXi_h)(\bPhi_2)\big\rangle_{\cH}\leq
	\big\langle \widetilde{\bPhi}_1, (\widehat{\bXi}_h - \bXi_h)(\widetilde{\bPhi}_1)\big\rangle_{\cH}
	-\big\langle \widetilde{\bPhi}_2, (\widehat{\bXi}_h - \bXi_h)({\widetilde\bPhi}_2)\big\rangle_{\cH},
$$
where $\widetilde{\bPhi}_1 = \bPhi_1 + \bPhi_2,$ $\widetilde{\bPhi}_2 = \bPhi_1 - \bPhi_2$ and
$\cM(\mF_{\bX},\widetilde{\bPhi}_i) \le 2\{\cM(\mF_{\bX},\bPhi_1)+ \cM(\mF_{\bX}, \bPhi_2) \}$ for $i=1, 2.$
Combing these with results and the proof of (\ref{thm_sig_k1}) leads to
\begin{eqnarray*}
	&&P\Big[\Big|\big\langle \bPhi_1, (\widehat{\bXi}_h - \bXi_h)(\bPhi_2)\big\rangle_{\cH}\Big|> 2\{\cM(\mF_{\bX},\bPhi_1)+ \cM(\mF_{\bX}, \bPhi_2)\}\eta\Big]\\
	&\leq& \sum_{i=1}^2P\Big[\Big|\big\langle \widetilde\bPhi_i, (\widehat{\bXi}_h - \bXi_h)(\widetilde\bPhi_i)\big\rangle_{\cH}\Big|> 2\cM(\mF_{\bX},\widetilde\bPhi_i)\eta\Big]\le 8 \exp\Big\{- c n\min\left( \eta^2,\eta\right)\Big\}
\end{eqnarray*}
for some constant $c> 0.$
This, together with, $\cM(\mF_{\bX},\bPhi_i) \le \cM_k(\mF_{\bX}) \big\langle \bPhi_i, \bXi_0(\bPhi_i)\big\rangle_{\cH}$ for $i=1, 2,$ implies (\ref{thm_sig_k2}), which completes the proof.  $\square$

\subsection{Proof of Proposition \ref{res.re}}
It is  easy to see that  $\btheta^T \widehat \bGamma \btheta  = \btheta^T \bGamma \btheta + \btheta^T (\widehat \bGamma - \bGamma) \btheta.$ Hence we have
\begin{eqnarray*}
	\btheta^T \widehat \bGamma \btheta	 \ge \btheta^T \bGamma \btheta -  \|\widehat \bGamma - \bGamma\|_{\max} \|\btheta\|_1^2.
\end{eqnarray*}
By Condition \ref{cond_min_bound}, $\lambda_{\min}(\bGamma) \ge \underline{\mu},$ where $\lambda_{\min}(\bGamma)$ denotes the minimum eigenvalue of $\bGamma.$ Together with Lemma~\ref{lemma_Gamma_max} in Section~\ref{sec.lemma_Gamma_max}, this proposition follows. $\square$
\subsection{Proof of Proposition~\ref{prop_Error_eigen}}
Note that on the event $\big\{|\widehat \lambda_{jl}-\lambda_{jl}| \le 2^{-1} \lambda_{jl}\big\},$ we have $\widehat\lambda_{jl}\geq \lambda_{jl}/2,$ $\widehat \lambda_{jl}^{-1/2} \leq \sqrt{2}\lambda_{jl}^{-1/2}$ and
$
|\widehat \lambda_{jl}^{-1/2} - \lambda_{jl}^{-1/2}| \le \frac{\widehat\lambda_{jl}^{-1}|\widehat\lambda_{jl}-\lambda_{jl}|\lambda_{jl}^{-1}}{\widehat\lambda_{jl}^{-1/2}+\lambda_{jl}^{-1/2}}\le 2\lambda_{jl}^{-3/2}|\widehat \lambda_{jl} - \lambda_{jl}|,
$ which implies that $\left|\frac{\widehat\lambda_{jl}^{-1/2}-\lambda_{jl}^{-1/2}}{\lambda_{jl}^{-1/2}}\right| \leq 2\left|\frac{\widehat\lambda_{jl}-\lambda_{jl}}{\lambda_{jl}}\right|.$ Then it follows from Theorem~\ref{thm_lam_phi} that there exist positive constants $C_{\lambda}, C_{\phi},$ $c_5$ and $c_6$ such that the first and second deviation bounds in (\ref{fvar.eigen}) respectively hold with probability greater than $1 -c_5(pq)^{-c_6}.$ The proof is complete. $\square$

\subsection{Proof of Proposition~\ref{prop_Error_max}}
Notice that
\begin{eqnarray}
\label{term1}
\widehat \bY_j - \widehat \bGamma \bB_j
& = & \Big\{(n-1)^{-1}\widehat \bD^{-1} \widehat \bZ^\T \widehat \bV_j - (n-1)^{-1}\bD^{-1} \tE(\bZ^\T \bV_j) \Big\}\\
& & 	+ (n-1)^{-1}\bD^{-1} \tE\Big\{\bZ^\T (\bV_j - \bZ\bD^{-1} \bB_j)\Big\} \nonumber - \big(\widehat \bGamma - \bGamma\big) \bB_j. 
\end{eqnarray}


First, we show the deviation bounds of $\widehat \bD^{-1} (n-1)^{-1}\widehat \bZ^\T \widehat \bV_j - \bD^{-1} \tE((n-1)^{-1}\bZ^\T \bV_j).$  We  decompose this term as
$\widehat \bD^{-1} \Big\{(n-1)^{-1}\widehat \bZ^\T \widehat \bV_j- \tE((n-1)^{-1}\bZ^\T \bV_j) \Big\}+ (\widehat \bD^{-1} - \bD^{-1}) \tE((n-1)^{-1}\bZ^\T \bV_j).$
It follows from 
Theorem~\ref{thm_cov} that there exists positive constants $C_1^*,$ $c_5$ and $c_6$ that
\begin{equation}
\label{term2-1}
\sup_{j,k}\left\|\bD_k^{-1} \Big\{(n-1)^{-1}\widehat \bZ_k^\T \widehat \bV_j- \tE\big((n-1)^{-1}\bZ_k^\T \bV_j\big) \Big\}\bD_j^{-1}\right\|_{\max} \le C^*_1 \cM_1(\mF_{\bX}) q^{\alpha + 1} \sqrt{\frac{\log (pq)}{n}},
\end{equation}
with probability greater than $1 - c_5 (pq)^{-c_6}.$ Note that $\widehat \bD_k = \mbox{diag}(\widehat \lambda_{k1}^{1/2}, \ldots, \widehat \lambda_{kq}^{1/2})$ and $\bD_k = \mbox{diag}(\lambda_{k1}^{1/2}, \ldots, \lambda_{kq}^{1/2}),$ it follows from Proposition~\ref{prop_Error_eigen} that there exists positive constant $C_2^*,$ such that
\begin{eqnarray}
\label{term2-2}
\left\|\left(\widehat \bD^{-1} - \bD^{-1}\right)\bD\right\|_{\max} \le C_2^* \cM_1(\mF_{\bX}) \sqrt{\frac{\log (pq)}{n}},
\end{eqnarray}
with probability great than $1 - c_5 (pq)^{-c_6}.$
By Condition~\ref{cond.cov.func},
we have $\max_j||\bD_j||_F \leq \lambda_0^{1/2}$ and $\|\bD^{-1}\tE\big((n-1)^{-1}\bZ^\T \bV_j\big)\|_{\max}^{(q)} \le q^{1/2}\|\bD^{-1}\tE\big((n-1)^{-1}\bZ^\T \bV_j\big)\bD_j^{-1}\|_{\max}||\bD_j||_F=O(q^{1/2}),$ where the fact that, for $q \times q$ matrix $\bA$ and a diagonal matrix $\bB,$ $\|\bA \bB\|_{F} \le q^{1/2} \|\bA\|_{\max} \|\bB\|_F,$ is used.
These results together with (\ref{term2-1}) and (\ref{term2-2}) imply that there exists $C_3^*$
\begin{eqnarray}
\label{term2}
\left\|\widehat \bD^{-1} (n-1)^{-1}\widehat \bZ^\T \widehat \bV_j- \bD^{-1} \tE\big((n-1)^{-1}\bZ^\T \bV_j\big)\right\|_{\max}^{(q)} \le C_3^* \cM_1(\mF_{\bX}) q^{\alpha + 3/2} \sqrt{\frac{\log (pq)}{n}}
\end{eqnarray}

Second, consider the bias term $
(n-1)^{-1}\bD^{-1} \tE\{\bZ^\T (\bV_j - \bZ\bD^{-1} \bB_j)\}.$ By Appendix~\ref{ap.vfar.deriv}, $\bR_j$ is a $(n-1) \times q$ matrix  whose row vectors are formed by $\{\br_{tj}, t=2, \dots, n\}$ with $\br_{tj}=(r_{tj1}, \dots, r_{tjq})^{\T}$ and $r_{tjl}=\sum_{k=1}^p \sum_{m=q+1}^{\infty} \langle\phi_{jl}, \langle A_{jk},\phi_{km}\rangle \rangle  \xi_{(t-1)km}$ for $l=1,\dots,q.$ It follows Conditions~\ref{cond.cov.func}, \ref{cond.fvar.bias} and similar arguments in deriving~(\ref{err.R}) and (\ref{term2}) that
there exists some positive constant $C_4^*$ such that
\begin{eqnarray}
\label{term3}
&&\big\|(n-1)^{-1}\bD^{-1} \tE\{\bZ^\T (\bV_j - \bZ\bD^{-1} \bB_j)\}\big\|_{\max}^{(q)}\nonumber\\
&&\leq q^{1/2}\big\|(n-1)^{-1}\bD^{-1} \tE(\bZ^\T\bR_j)\widetilde \bD^{-1}\big\|_{\max} ||\widetilde\bD||_F
\le C_4^* s_j q^{-\beta + 1 },
\end{eqnarray}
where $\widetilde \bD = \lambda_0 \bI_q.$

Third, it follows from Lemma~\ref{lemma_Gamma_max} in Appendix~\ref{sec.lemma_Gamma_max}, Lemma~15 in the Supplementary Material of {\color{blue} Qiao et al. (2019)} and $||\bB||_1^{(q)}= O\big(\lambda_0^{1/2}s_j\big)$ that there exist some positive constants $C_5^*$ such that
\begin{eqnarray}
\label{term4}
\left \|\big(\widehat \bGamma - \bGamma\big)\bB_j \right\|_{\max}^{(q)} \leq \left \|\widehat \bGamma - \bGamma \right\|_{\max}^{(q)}||\bB_j||_1^{(q)}\le \cM_1(\mF_{\bX}) s_j q^{\alpha + 2}\sqrt{\frac{\log (pq)}{n}},
\end{eqnarray}
with probability great than $1 - c_5 (pq)^{-c_6}.$

Combing results in (\ref{term1}), (\ref{term2}), (\ref{term3}) and  (\ref{term4}) implies that there exist positive constants $C_E, c_4$ and $c_5$ such that
$$ 
||\widehat \bY_j- \widehat \bGamma\bB_j||_{\max}^{(q)} \leq C_E \cM_1(\mF_{\bX}) s_j \Big \{q^{\alpha +2} \sqrt{\frac{\log (pq)}{n}}  + q^{-\beta + 1}\Big \}, ~ j =1,\ldots, p,
$$
with probability greater than $1 -c_5(pq)^{-c_6}.$ The proof is complete. $\square$

\subsection{Proposition~\ref{thm_F_bound} and its proof}
\begin{proposition}
\label{thm_F_bound}
	Suppose that Conditions~\ref{cond.cov.func} and \ref{cond.bd.fsm} hold. Then
	 there exists some universal constant $\tilde c > 0$ such that
	for any $\eta > 0$
	\begin{equation}
	\label{thm_Sigma_F}
	P\left\{\big\|\widehat{\bSigma}_0 - \bSigma_0\big\|_{F} >   2 \cM_1(\mF_{\bX}) \lambda_0 \eta \right\} \le \frac{p^2}{\eta^2 n} \big (16 \tilde c^{-1} + 128 \tilde c^{-2}n^{-1} \big).
	\end{equation}
	In particular, if the sample size $n$ satisfies the bound $n>128(\tilde\rho^2\tilde c^{2}-16\tilde c)^{-1},$  where $\tilde\rho$ is some positive constant with $\tilde\rho>4\tilde c^{-1/2},$ then with probability greater than $1 - \tilde\rho^{-2} \big (16 \tilde c^{-1} + 128 \tilde c^{-2}n^{-1} \big),$ the estimate $\widehat\bSigma_0$ satisfies the bound
	\begin{equation}
	\label{bd_Sigma_F}
	\big\|\widehat{\bSigma}_0 - \bSigma_0\big\|_{F} \le  2 \cM_1(\mF_{\bX}) \lambda_0 \tilde\rho\sqrt{\frac{p^2}{n}}.
	\end{equation}
\end{proposition}
{\it Proof}. It follows from the definition of $
\|\widehat \bSigma_0 - \bSigma_0\|_{F}^2 = \sum_{j, k =1}^p \|\widehat \Sigma_{jk}^{(0)} - \Sigma_{jk}^{(0)}\|_{\cS}^2,$
Chebyshev's inequality and (\ref{Sigma_moment}) with $q=1$ that for any
$\eta > 0,$
\begin{eqnarray*}
	P\left\{\|\widehat \bSigma_0 - \bSigma_0\|_{F}>  2\cM_1(\mF_{\bX}) \lambda_0\eta \right\} & \le&\frac  {1} {\left(2 \cM_1(\mF_{\bX}) \lambda_0 \right)^{2} \eta^2} \sum_{j, k =1}^p \tE\|\widehat \Sigma_{jk}^{(0)} - \Sigma_{jk}^{(0)}\|_{\cS}^2 \\
	&\le &  \frac{p^2}{\eta^2} (16 \tilde c^{-1}n^{-1} + 128 \tilde c^{-2}n^{-2}) \\
	&=&  \frac{p^2}{\eta^2 n} (16 \tilde c^{-1} + 128 \tilde c^{-2}n^{-1}).
\end{eqnarray*}
By letting $\eta = \tilde\rho \sqrt{p^2/n}$ with $\rho > 0$, we have that
\begin{eqnarray*}
	P\left\{\|\widehat \bSigma_0 - \bSigma_0\|_{F}>  2 \cM_1(\mF_{\bX}) \lambda_0\tilde\rho\sqrt{\frac{p^2}{n}} \right\}
	\le \tilde\rho^{-2} (16 \tilde c^{-1} + 128 \tilde c^{-2}n^{-1}).
\end{eqnarray*}
The proof is complete. $\square$

\subsection{Proof of Lemma~\ref{lemma_lambda}}
By Lemma~\ref{lemma.lambda.phi}, we obtain that
\begin{equation}
\label{lambda_expand}
\frac{\widehat \lambda_{jl} - \lambda_{jl}}{\lambda_{jl}} = \frac{\big\langle \phi_{jl}, \widehat \Delta_{jj}(\phi_{jl})\big \rangle}{\lambda_{jl}}  + \frac{R_{jl}}{\lambda_{jl}}, ~j =1,\dots,p, l = 1,\ldots,L.
\end{equation}
Note that $\lambda_{jl} = \big\langle \phi_{jl}, \Sigma_{jj}(\phi_{jl})\big \rangle$. It follows from (\ref{thm_Phi_1}) in Theorem~\ref{thm_con_op} that for any $\eta > 0$,
\begin{equation}
\label{con.bd.lambda.ratio}
    P\left\{\left|\frac{\big\langle \phi_{jl}, \widehat \Delta_{jj}(\phi_{jl}) \big \rangle}{\lambda_{jl}}\right| > \cM_1(\mF_{\bX}) \eta\right\} \le 2 \exp\Big\{- c n \min(\eta^2,\eta)\Big\}.
\end{equation}
We next turn to the term $\big|R_{jl}/\lambda_{jl}\big|$. By (\ref{eigen.bd}), Lemma~\ref{lemma.lambda.phi} and Condition~\ref{cond_eigen} with $\delta_{jl} \ge c_0 l^{-\alpha -1}$ and $\lambda_{jl} \ge c_0 \alpha^{-1} l^{-\alpha},$ we have
$$
\left|\frac{R_{jl}}{\lambda_{jl}}\right| \le 4\sqrt{2}c_0^{-2} \alpha l^{2\alpha+1} \|\widehat \Delta_{jj}\|_{\cS}^2.
$$
It then follows from (\ref{thm_Sigma_compt}) in Theorem~\ref{thm_max_bound} that there exists some constant $\tilde c> 0 $ such that for any $\eta > 0$
\begin{equation}
\label{con.bd.R}
P\left\{ \left|\frac{R_{jl}}{\lambda_{jl}}\right| > 4\sqrt{2}c_0^{-2} \alpha l^{2\alpha+1} \{2 \cM_1(\mF_{\bX}) \lambda_0\eta\}^2\right\} \le 4 \exp\Big\{- \tilde c n \min(\eta^2,\eta)\Big\}.
\end{equation}
Let $\tilde c_1 = \min(c, \tilde c)$. It follows from $\rho_1 = 16\sqrt{2} c_0^{-2} \alpha \lambda_0^2$, (\ref{lambda_expand}), (\ref{con.bd.lambda.ratio}) and (\ref{con.bd.R}) that (\ref{thm_con_lambda}) in Lemma~\ref{lemma_lambda} holds, which completes our proof. $\square$

\subsection{Proof of Lemma~\ref{lemma_phi}}
It follows from (\ref{eigen.bd}), Condition~\ref{cond_eigen} with $\delta_{jl} \ge c_0 l^{-\alpha -1}$ and (\ref{thm_Sigma_compt}) in Theorem~\ref{thm_max_bound} that there exists some universal constant $\tilde c$ such that for any $\eta>0,$ (\ref{thm_con_phi}) holds.$\square$

\subsection{Lemma~\ref{lemma_phi_g} and its proof}
\label{ap.lemma_phi_g}
\begin{lemma}
    \label{lemma_phi_g}
Suppose that Conditions~\ref{cond.cov.func}--\ref{cond_eigen} hold. Then there exists some universal constant $\tilde c_2>0$ such that for each $j=1,\dots,p, l=1, \dots,d_j,$ any given function $g \in \cH$ and $\eta > 0,$
	\begin{equation}
	\label{thm_con_phi_g}
	\begin{split}
	&P\left\{\left|\big \langle \widehat \phi_{jl} - \phi_{jl}, g\big \rangle \right| \ge
	\rho_2\|g^{-jl}\|_{\lambda}  \cM_1(\mF_{\bX}) \lambda_{jl}^{1/2}l^{\alpha + 1}  \eta
	+ \rho_3 \|g\|\cM_1^2(\mF_{\bX})  l^{2(\alpha+1)}\eta^2 \right\} \\
	&\le 8 \exp\Big\{ - \tilde c_2 n\min(\eta^2, \eta)\Big\} + 4\exp\Big\{ - \tilde c_2  \cM_1^{-2}(\mF_{\bX}) n l^{-2(\alpha + 1)}\Big\},
	\end{split}
	\end{equation}
where $g(\cdot)=\sum_{m=1}^{\infty}g_{jm}\phi_{jm}(\cdot),$ $\|g^{-jl}\|_{\lambda} = \big({\sum}_{m: m\neq l}\lambda_{jm} g_{jm}^2 \big)^{1/2},$
$\rho_{2} = 2 c_0^{-1}  
$
and
$\rho_{3} = 4(6+2\sqrt{2}) c_0^{-2} \lambda_0^2$ with $c_0 \leq 4\cM_1(\mF_{\bX})\lambda_0l^{\alpha+1}.$
\end{lemma}
{\it Proof.} It follows from the expansion $g(\cdot) = \sum_{m = 1}^\infty g_{jm} \phi_{jm}(\cdot)$ and (\ref{express.phi}) in Lemma~\ref{lemma.lambda.phi} in Appendix~\ref{sec.pf.lemma.lambda.phi} that
\begin{eqnarray*}
\big \langle \widehat \phi_{jl} - \phi_{jl}, g\big \rangle &=&
 \sum_{m: m \neq l}(\widehat \lambda_{jl} - \lambda_{jm})^{-1} g_{jm}
\big \langle \widehat \phi_{jl}, \langle \widehat \Delta_{jj}, \phi_{jm} \rangle \big \rangle +
g_{jl} \big \langle \widehat\phi_{jl} - \phi_{jl}, \phi_{jl} \big\rangle\\
& = &\sum_{m: m \neq l}\Big \{(\widehat \lambda_{jl} - \lambda_{jm})^{-1} -(\lambda_{jl} - \lambda_{jm})^{-1} \Big\} g_{jm}
\big \langle \widehat \phi_{jl}, \langle \widehat \Delta_{jj}, \phi_{jm} \rangle \big \rangle \\
 && + \sum_{m: m \neq l}(\lambda_{jl} - \lambda_{jm})^{-1} g_{jm}
\big \langle \widehat \phi_{jl} - \phi_{jl}, \langle \widehat \Delta_{jj}, \phi_{jm} \rangle \big \rangle \\
&& +  \sum_{m: m \neq l}(\lambda_{jl} - \lambda_{jm})^{-1} g_{jm}
\big \langle \phi_{jl}, \langle \widehat \Delta_{jj}, \phi_{jm} \rangle \big \rangle
+ g_{jl} \big \langle \widehat\phi_{jl} - \phi_{jl}, \phi_{jl} \big\rangle\\
& = & I_{1} +I_{2} + I_3 + I_4.
\end{eqnarray*}

Let $\Omega_{d_j} = \{ 2 \|\widehat \Delta_{jj}\|_{\cS} \le \delta_{jd_j}\}.$ It follows from Condition~\ref{cond_eigen} and (\ref{thm_Sigma_compt}) in Theorem~\ref{thm_max_bound} with the choice of $\eta=\{4\cM_1(\mF_{\bX})\lambda_0l^{\alpha+1}\}^{-1}c_0\leq 1$ that
\begin{equation}
\label{bd.Delta}
    P\big(\Omega_{d_j}^C\big) \leq P\big( \|\widehat \Delta_{jj}\|_{\cS}\geq 2^{-1}c_0 l^{-\alpha-1}\big) \leq 4 \exp\left\{-16^{-1}\tilde c\cM_1^{-2}(\mF_{\bX})\lambda_0^{-2}c_0^2 l^{-2(\alpha + 1)} n\right\}.
\end{equation}
On the event $\Omega_{d_j},$ we can see that ${\sup}_{l \leq d_j}|\widehat\lambda_{jl}-\lambda_{jl}|\leq \lambda_{jd_j}/2,$ which implies that $2^{-1}\lambda_{jl} \leq \widehat\lambda_{jl} \leq 2\lambda_{jl}.$ Moreover, $|\widehat \lambda_{jl}-\lambda_{jl}| \leq 2^{-1}|\lambda_{jl}-\lambda_{jm}|$ for $1\leq l \neq m \leq d_j$ and hence $
|\widehat \lambda_{jl} - \lambda_{jm}|\ge 2^{-1} | \lambda_{jl} - \lambda_{jm}|$ for $j=1,\ldots,p.$ By Condition~\ref{cond_eigen}, $|\lambda_{jl} - \lambda_{jm}| \ge c_0 l^{- \alpha -1}$ for $1\leq m \neq l \leq d_j.$ Using the above results, we have
\begin{eqnarray*}
|I_1|^2 &\le& \big(\widehat \lambda_{jl} - \lambda_{jl}\big)^{2} \sum_{m: m \neq l}(\widehat \lambda_{jl} - \lambda_{jm})^{-2}(\lambda_{jl} - \lambda_{jm})^{-2} g_{jm}^2 \| \widehat \Delta_{jj}\|_{\cS}^2\\
&\le& 4\big(\widehat \lambda_{jl} - \lambda_{jl}\big)^{2}\| \widehat \Delta_{jj}\|_{\cS}^2\sum_{m: m \neq l}(\lambda_{jl} - \lambda_{jm})^{-4} g_{jm}^2\\
&\le& 4 c_0^{-4}\|g^{-jl}\|^2 l^{4(\alpha+1)} \big(\widehat \lambda_{jl} - \lambda_{jl}\big)^{2}  \|\widehat \Delta_{jj}\|_{\cS}^2,
\end{eqnarray*}
where $||g^{-jl}||=(\sum_{m:m \neq l} g_{jm}^2)^{1/2}.$
This together with (\ref{eigen.bd}) implies that, on the event $\Omega_{d_j},$
$$ |I_1| \le 2 c_0^{-2} \|g^{-jl}\| l^{2(\alpha+1)}\big \| \widehat \Delta_{jj}\big\|_{\cS}^2.$$
Similarly, we can show that
$$
|I_2| \le c_0^{-1} \|g^{-jl}\|  l^{\alpha+1} \|\widehat \phi_{jl} - \phi_{jl}\| \|\widehat \Delta_{jj}\|_{\cS} \le 2\sqrt{2}c_0^{-2} \|g^{-jl}\| l^{2(\alpha+1)}  \|\widehat \Delta_{jj}\|_{\cS}^2.
$$
Moreover, by the result $||\widehat\phi_{jl}-\phi_{jl}||^2=\langle\widehat\phi_{jl}-\phi_{jl},-2\phi_{jl}\rangle$ and (\ref{eigen.bd}) we have
$$|I_4| = 2^{-1}|g_{jl}| \|\widehat \phi_{jl} - \phi_{jl}\|^2 \le 4 c_0^{-2} |g_{jl}| l^{2(\alpha + 1)} \|\widehat \Delta_{jj}\|_{\cS}^2.
$$
Combing the above upper bound results, we have
$$
|I_1| + |I_2| +|I_4| \le (6+2\sqrt{2}) c_0^{-2} \|g\| l^{2(\alpha + 1)} \|\widehat \Delta_{jj}\|_{\cS}^2.
$$

Let $\widetilde{\lambda}_{g} = {\sum}_{m: m \neq l} \lambda_{jm} (\lambda_{jl} - \lambda_{jm})^{-2} g_{jm}^2 \leq c_0^{-2}  l^{2(\alpha + 1)}   \|g^{-jl}\|_{\lambda}^2.$ Then it follows from (\ref{thm_Phi_12}) in Theorem~\ref{thm_con_op} that 
\begin{equation}
\label{I3.con.bd}
P\left\{ \Big|\lambda_{jl}^{-1/2}\widetilde{\lambda}_{g}^{-1/2} I_3\Big | \ge 2 \cM_1(\mF_{\bX}) \eta\right \} \le 4 \exp\Big\{ - c n \min(\eta^2, \eta)\Big\}.
\end{equation}
Define $\Omega_{1,\eta} = \Big \{ \|\widehat \Delta_{jj}\|_{\cS} \le 2 \cM_1(\mF_{\bX}) \lambda_0 \eta\Big\}$ and
$
\Omega_{2,\eta} = \Big\{ \big|I_3\big| \le 2 c_0^{-1}  \lambda_{jl}^{1/2} \|g^{-jl}\| _{\lambda}\cM_1(\mF_{\bX}) l^{\alpha + 1} \eta\Big \}.
$ Let
$
\rho_{2} = 2 c_0^{-1}  
$
and
$\rho_{3} = 4(6+2\sqrt{2}) c_0^{-2} \lambda_0^2.$
Under the event $\Omega_{d_j} \cap \Omega_{1,\eta}\cap \Omega_{2,\eta},$ we obtain that
$$
\left|\big \langle \widehat \phi_{jl} - \phi_{jl}, g\big \rangle \right| \le
\rho_2\|g^{-jl}\|_{\lambda}  \cM_1(\mF_{\bX}) \lambda_{jl}^{1/2}l^{\alpha + 1}  \eta
+ \rho_3 \|g\|\cM_1^2(\mF_{\bX})  l^{2(\alpha+1)}\eta^2.
$$

Let $\tilde c_2 = \min(16^{-1}\lambda_0^{-2}c_0^2\tilde c,\tilde c_1)$.
It follows from (\ref{thm_Sigma_compt}) in Theorem~\ref{thm_max_bound} and (\ref{I3.con.bd}) that $$P\big(\Omega_{1,\eta}^C\cup \Omega_{2,\eta}^C\big)
\le 8 \exp\Big\{ - \tilde c_2n \min(\eta^2, \eta)\Big\}.$$ This together with (\ref{bd.Delta}) completes the proof of (\ref{thm_con_phi_g}). $\square$

\subsection{Lemma~\ref{lemma_Gamma_max} and its proof}
\label{sec.lemma_Gamma_max}
\begin{lemma}
	\label{lemma_Gamma_max}
	Suppose that Conditions \ref{cond.cov.func}-\ref{cond_eigen} hold.  Then there exist some positive constants $C_{\Gamma}$, $c_5$ and $c_6$ such that
	$$
	\big\|\widehat \bGamma - \bGamma\big\|_{\max} \le C_\Gamma \cM_1(\mF_{\bX}) q^{\alpha +1}\sqrt{\frac{\log (pq)}{n}}
	$$
	with probability greater than $1 - c_5(pq)^{-c_6}$.
\end{lemma}
{\it Proof.}
Note that
$$
\|\widehat \bGamma - \bGamma\|_{\max} = \max_{1 \le j,k \le p, 1 \le l,m \le q} \Big| \widehat{\lambda}_{jl}^{-1/2} \widehat{\lambda}_{km}^{-1/2} \widehat \sigma_{jklm} - \lambda_{jl}^{-1/2} \lambda_{km}^{-1/2}\sigma_{jklm}\Big|.
$$
Let $\widehat{s}_{jklm} = \frac{\widehat \lambda_{jl} \widehat \lambda_{km}}{\lambda_{jl} \lambda_{km}}$ for each $(j,k,l,m)$. Then we have
\begin{eqnarray*}
	\widehat{\lambda}_{jl}^{-1/2} \widehat{\lambda}_{km}^{-1/2} \widehat \sigma_{jklm} - \lambda_{jl}^{-1/2} \lambda_{km}^{-1/2}\sigma_{jklm}
	= \widehat{s}_{jklm}^{-1/2}\left(\frac{\widehat \sigma_{jklm} - \sigma_{jklm}}{\lambda_{jl}^{1/2} \lambda_{km}^{1/2}} \right)
	+ \Big(\widehat{s}_{jklm}^{-1/2} - 1\Big)\frac{\sigma_{jklm}}{\lambda_{jl}^{1/2} \lambda_{km}^{1/2}}.
\end{eqnarray*}
Let
$
\Omega_{\lambda} = \left\{\sup_{1\le j \le p, 1\le l \le q}\left|\frac{\widehat \lambda_{jl} - \lambda_{jl}}{\lambda_{jl}}\right| \le 1/5\right\}.
$
Observe that
$$
\widehat s_{jklm} - 1 = \left(\frac{\widehat \lambda_{jl} - \lambda_{jl}}{\lambda_{jl}} + 1\right)\left( \frac{\widehat \lambda_{km} - \lambda_{km} }{\lambda_{km}}\right) + \frac{\widehat \lambda_{jl} - \lambda_{jl}}{ \lambda_{jl}}.
$$
Then under the event $\Omega_{\lambda},$ we have $|\widehat s_{jklm} - 1| \le 1/2,$ and thus $\widehat{s}_{jklm}^{-1/2} \le \sqrt{2}.$ Moreover, provided that fact that $|(1+x)^{-1/2} - 1| \le x$ if $|x| \le 1/2$, we have
$$
\Big|\widehat{s}_{jklm}^{-1/2} - 1\Big| \le \frac{6}{5} \left(\left|  \frac{\widehat \lambda_{km} - \lambda_{km} }{\lambda_{km}}\right| + \left| \frac{\widehat \lambda_{jl} - \lambda_{jl} }{\lambda_{jl}}\right|\right).
$$

Under the event $\Omega_{\lambda},$ the above results together with the fact of $\sigma_{jklm}\leq \lambda_{jl}^{1/2}\lambda_{km}^{1/2}$ imply that
$$
\|\widehat \bGamma - \bGamma\|_{\max} \le \sqrt{2}\max_{1 \le j,k \le p, 1 \le l,m \le q}\left|\frac{\widehat \sigma_{jklm} - \sigma_{jklm}}{\lambda_{jl}^{1/2} \lambda_{km}^{1/2}} \right| + \frac{12}{5}\max_{1\le j \le p, 1\le l \le q}\left|\frac{\widehat \lambda_{jl} - \lambda_{jl}}{\lambda_{jl}}\right|.
$$
Then it follows from Theorems~\ref{thm_lam_phi} and \ref{thm_cov} that there exist some positive constants $C_\Gamma$, $c_5$ and $c_6$ such that
$$
\|\widehat \bGamma - \bGamma\|_{\max} \le C_\Gamma \cM_1(\mF_{\bX}) q^{\alpha +1}\sqrt{\frac{\log (pq)}{n}}
$$
with probability greater than $1 - c_5(pq)^{-c_6}$. The proof is complete. $\square$

\subsection{Lemma~\ref{lemma.moment} and its proof}
The following lemma  shows how to derive the tail probability through moment conditions.

\begin{lemma}
	\label{lemma.moment}
	Let $\tX$ be a random variable. If for some constants $ c_1,c_2 > 0$
	$$
	P\left(|\tX| > t \right ) \le c_1 \exp\{- c_2^{-1}\min(t^2, t)\} ~~ \mbox{for any } t >0,
	$$
	then for any integer $q \ge 1$,
	$$
	\tE(\tX^{2q}) \le q!c_1(4c_2)^q + (2q)!c_1(4c_2)^{2q}.
	$$
	Conversely, if for some positive constants $a_1,a_2$,
	$
	\tE(\tX^{2q}) \le q!a_1a_2^q + (2q)!a_1a_2^{2q}, ~ q \ge 1,
	$ then by letting $c^*_2 = 8\max\{4(a_2 + a^2_2), a_2\}$ and $c_1^* = a_1$, we have that
	$$
	P\left(|\tX| > t \right ) \le c_1^* \exp\{- c^{*-1}_2\min(t^2, t)\} ~~ \mbox{for any } t >0.
	$$
\end{lemma}
{\it Proof.} This lemma can be proved in a similar way to Theorem 2.3 of {\color{blue}Boucheron et al. (2014)} and hence the proof is omitted here. In the proof, the following two
inequalities are used, i.e. for any $c >0$ and $t > 0$,
$$
\frac{1}{2}\min(t^2,t) \le \frac{t^2}{1 + t} \le  \min(t^2,t),
$$
and
$$
\sqrt{\frac{ct}{2}} + \frac{ct}{2} \le \frac{c(t + \sqrt{t^2 + 4t/c})}{2} \le
\sqrt{ct} + ct.~~\square
$$
\subsection{Lemma~\ref{lemma.lambda.phi} and its proof}
\label{sec.pf.lemma.lambda.phi}
\begin{lemma}
	\label{lemma.lambda.phi}
	For each $j = 1,\ldots,p$ and $l=1,\dots,$ the term of $\widehat \lambda_{jl} - \lambda_{jl}$ can be expressed as
	\begin{eqnarray}
	\label{express.lambda}
	\widehat \lambda_{jl} - \lambda_{jl} = \big\langle \phi_{jl}, \widehat \Delta_{jj}(\phi_{jl}) \big \rangle  + R_{jl},
	\end{eqnarray}
	where $|R_{jl}| \le 2 \|\widehat \phi_{jl} - \phi_{jl}\|\|\widehat \Delta_{jj} \|_{\cS}$.
	Furthermore, if $\inf_{m: m \neq l} |\widehat \lambda_{jl} - \lambda_{jm}| > 0$, then
	\begin{eqnarray}
	\label{express.phi}
	\widehat \phi_{jl} - \phi_{jl} = \sum_{m: m \neq l}(\widehat \lambda_{jl} - \lambda_{jm})^{-1} \phi_{jm}
	\big \langle \widehat \phi_{jl}, \widehat \Delta_{jj}(\phi_{jm}) \big \rangle
	+\phi_{jl}\big \langle \widehat\phi_{jl} - \phi_{jl}, \phi_{jl} \big\rangle.
	\end{eqnarray}
\end{lemma}
{\it Proof}. This lemma follows directly from Lemma 5.1 of {\color{blue}Hall and Horowitz (2007)} and hence the proof is omitted here. $\square$

\subsection{Lemma~\ref{lemma.sub.ieq} and its proof}
\label{sec.lemma.sub.ieq}
\begin{lemma}
	\label{lemma.sub.ieq}
	For a $p \times p$ lag-$h$ autocovariance function, $\bSigma_h=\big(\Sigma_{jk}^{(h)}\big)_{1 \leq j,k \leq p},$ $h=0, 1, \dots,$ we have
	$$\|\Sigma_{jk}^{(h)}\|_{\cS} \leq \lambda_0~~\mbox{and}~~\|\Sigma_{jk}^{(h)}(\phi_{km})\| \leq \lambda_{km}^{1/2}\lambda_0^{1/2} ~\text{ for } m \geq 1.$$
\end{lemma}
{\it Proof}. By the expansion $\Sigma_{jk}^{(h)}(u,v)=\sum_{l,m=1}^{\infty} E\{\xi_{tjl}\xi_{(t+h)km}\}\phi_{jl}(u)\phi_{km}(v),$ the orthonormality of $\{\phi_{jl}(\cdot)\}$ for $\{\phi_{km}(\cdot)\}$ for each $(j,k)$ and the Cauchy-Schwarz inequality, we have 
$$
\|\Sigma_{jk}^{(h)}\|_{\cS}  \leq  \Big[\sum_{l,m=1}^{\infty} E(\xi_{tjl}^2)E\{\xi_{(t+h)km}^2\}\Big]^{1/2} \leq \Big\{\sum_{l=1}^{\infty} \lambda_{jl} \sum_{m=1}^{\infty} \lambda_{km}\Big\}^{1/2} \leq \lambda_0. $$

Moreover, applying similar techniques, we have
\begin{eqnarray*}
\|\Sigma_{jk}^{(h)}(\phi_{km})\|^2&=& \int_{\cU}\Big(\sum_{l=1}^{\infty}E(\xi_{tjl}\xi_{tkm})\phi_{jl}(u)\Big)^2du \\
&\leq& \sum_{l=1}^{\infty}E(\xi_{tjl}^2)E(\xi_{tkm}^2)= \int_{\cU}\Sigma_{jj}(u,u)du\lambda_{km}\leq \lambda_0\lambda_{km},
\end{eqnarray*}
which completes the proof. $\square.$

\subsection{Lemma~\ref{lemma.err.A} and its proof}
\label{sec.pf.lemma.err.A}
\begin{lemma}
	\label{lemma.err.A}
	For each $j,k=1,\dots,p,$ let $\{\phi_{jl}(\cdot)\}_{1\leq l \leq q}$ and $\{\widehat\phi_{jl}(\cdot)\}_{1 \leq l \leq q}$ correspond to true and estimated eigenfunctions, respectively, and $\widehat\psi_{jklm}$ be the estimate of $\psi_{jklm}$ for $l,m=1,\dots,q.$ Then we have
	\begin{eqnarray*}
		\left\|\sum_{l=1}^q\sum_{m=1}^q\widehat \phi_{km}(\cdot)\widehat \psi_{jklm} \left\{\widehat \phi_{jl}(\cdot) -\phi_{jl}(\cdot) \right\}\right\|_{\cS}^2 & \leq&
		\sum_{l=1}^q\|\widehat \phi_{jl} - \phi_{jl}\|^2 \sum_{l=1}^q\sum_{m=1}^q\widehat \psi_{jklm}^2,\\
		\left \| \sum_{l=1}^q\sum_{m=1}^q\phi_{km}(\cdot)(\widehat \psi_{jklm} - \psi_{jklm})\phi_{jl}(\cdot)\right\|_{\cS}^2&=&  \sum_{l=1}^q\sum_{m=1}^q(\widehat \psi_{jklm} - \psi_{jklm})^2.
	\end{eqnarray*}
\end{lemma}
{\it Proof}. We prove the first result
\begin{eqnarray*}
	&&  \left\|\sum_{l =1}^q \sum_{m=1}^q\widehat \phi_{km}\widehat \psi_{jklm} (\widehat \bphi_{jl} -\bphi_{jl}) \right\|_\cS^2\\
	& &=\sum_{m =1 }^q \left \|\sum_{l=1}^q \widehat \psi_{jklm} (\widehat \bphi_{jl} -\bphi_{jl}) \right\|^2
	\leq \sum_{m =1 }^q \sum_{l=1}^q \widehat \psi_{jklm}^2 \sum_{l = 1}^q \|\widehat \bphi_{jl} -\bphi_{jl}\|^2,
\end{eqnarray*}
where the first equality from the orthonormality of $\{\widehat \phi_{km}(\cdot)\}_{1 \leq m \leq q}$ and the second inequality comes from Cauchy-Schwarz inequality.
By the orthonormality of of $\{\phi_{km}(\cdot)\}_{1 \leq m \leq q}$  and $\{\phi_{jl}(\cdot)\}_{1 \leq l \leq q},$ we can prove the second result
\begin{eqnarray*}
	&&  \left\|\sum_{l =1}^q \sum_{m=1}^q \phi_{km}(\widehat \psi_{jklm}-\psi_{jklm})\bphi_{jl} \right\|_\cS^2\\
	& = &\sum_{m =1 }^q \left \|\sum_{l=1}^q (\widehat \psi_{jklm} - \psi_{jklm})\bphi_{jl} \right\|^2
	= \sum_{m =1 }^q \sum_{l=1}^q (\widehat \psi_{jklm}-\psi_{jklm})^2,
\end{eqnarray*}
which completes the proof. $\square$

\subsection{Lemma~\ref{lemma.norm.ieq} and its proof}
\label{sec.lemma.norm.ieq}
\begin{lemma}
	\label{lemma.norm.ieq}
	Let $\bA,\bB \in \eR^{pq\times q}$ with $j$-th blocks given by $\bA_j,\bB_j \in \eR^{q\times q},$ respectively. We have
	\begin{equation}
	\label{norm.eq1}
	\llangle\bA,\bB\rrangle \leq ||\bB||^{(q)}_{\max} ||\bA||_1^{(q)}.
	\end{equation}
\end{lemma}
{\it Proof}. By the definition and Cauchy-Schwarz inequality
\begin{eqnarray*}
	\llangle\bA,\bB\rrangle &=& \sum_{j=1}^p \llangle\bA_j,\bB_j\rrangle\\
	&\leq& \sum_{j=1}^p\llangle\bA_j,\bA_j\rrangle^{1/2} \llangle\bB_j,\bB_j\rrangle^{1/2}\\
	&\leq& \underset{j}{\max}||\bB_j||_F\sum_{j=1}^p||\bA_j||_F=||\bB||^{(q)}_{\max}||\bA||_1^{(q)},
\end{eqnarray*}
which completes the proof. $\square$

\renewcommand{\theequation}{D.\arabic{equation}}
\section{An illustrative example}
\label{ap.ill.ex}
In the following, for any $\bA=(A_{jk})_{1 \leq j,k \leq p},\bB=(B_{jk})_{1 \leq j,k \leq p}$ with their $(j,k)$-th components $A_{jk},B_{jk} \in \eS$ and $\bx \in \cH,$
write $  \bA \bB,$ $\bA\bx$ and $\bx^{\T}\bA$ for
$$
\int_{\cU} \bA(u,v') \bB(v',v)dv' ~\mbox{,}~ \int_{\cU} \bA(u,v) \bx(v)dv  ~\mbox{and}~ \int_{\cU} \bx(u)^{\T}\bA(u,v) du,
$$ respectively. For a $p \times p$ matrix, $\bC,$ we denote its maximum eigenvalue, spectral radius and operator norm by $\lambda_{\max}(\bC),$ $\rho(\bC)=|\lambda_{\max}(\bC)|$ and $\|\bC\|=\sqrt{\lambda_{\max}(\bC^{\T}\bC)},$ respectively.


Let $p=2$, $\bx_t=\big(x_{t1},x_{t2}\big)^{\T},$
$\bpsi=\text{diag}\big(\psi_1,\psi_2\big),$
$\bC=\left(\begin{array}{cc}a& b\\0 & a\end{array}\right)$
and $\be_t=\big(e_{t1},e_{t2}\big)^{\T},$ then the VFAR(1) model in
(\ref{vfar1.ex}) and (\ref{vfar1.ex.A}) can be rewritten as
$$
\bpsi(u)\bx_t = \int_{\cU} \bpsi(u)\bC\bpsi(v)\bpsi(v)\bx_{t-1}dv+ \bpsi(u)\be_t,
$$ which leads to a VAR(1) model
\begin{equation}
\label{var1.ex}
\bx_t = \bC \bx_{t-1} + \be_t.
\end{equation}

Provided that $\bA(u,v)=\bpsi(u)\bC\bpsi(v)$ and $||\bC||=\sqrt{\lambda_{\max}(\bC^{\T}\bC)}=\lambda_1$ with
$\bC^{\T}\bC\by=\lambda_1^2\by$ for $||\by||=1,$
it is easy to see that
$$\bA^{\T}\bA=\int \bA^{\T}(u,v')\bA(v',v)dv'=
\int\bpsi(u)\bC^{\T}\bpsi(v')\bpsi(v')\bC\bpsi(v)dv'=\bpsi(u)\bC^{\T}\bC\bpsi(v)$$
and $$\int(\bA^{\T}\bA)(u,v)(\bpsi(v)\by)dv=\int\bpsi(u)\bC^{\T}\bC\bpsi(v)\bpsi(v)\by dv=\bpsi(u)\bC^{\T}\bC\by=\lambda_1^2\bpsi(u)\by.$$ Hence $||\bA||_{\cL}=\sqrt{\lambda_{\max}(\bA^{\T}\bA)}=||\bC||=\lambda_1.$ The left side of Figure~\ref{ill.example} plots $||\bA||_{\cL}$ vs $b$ for different values of $a \in (0,1).$

Let $(\omega_j,\bv_j), j=1,2,$ be the eigen-pairs of $\bC$ satisfying $\bC\bv_j=\omega_j\bv_j.$ Then $$\int_{\cU} \bpsi(u)\bC\bpsi(v)\bpsi(v)\bv_j dv=\int_{\cU}\bA(u,v)\bpsi(v)\bv_j dv=\omega_j\bpsi(u)\bv_j.$$ Hence $\bA$ and $\bC$ share the same eigenvalues, which are $\omega_1=\omega_2=a.$ When $\rho(\bA)=\rho(\bC)=|a|<1,$ (\ref{vfar1.ex}) and (\ref{var1.ex}) correspond to stationary VFAR(1) and VAR(1) models, respectively.

For the VFAR(1) model in (\ref{vfar1.ex}), the spectral density function and the covariance function of $\{\bX_t\}_{t \in \eZ}$ are
\begin{equation}
\label{density.vfar1}
f_{\bX,\theta} = \frac{1}{2\pi}\left(\bSigma_0  + \sum_{h = 1}^\infty \Big\{ \bSigma_0 (\bA^{\T})^{h} \exp(-i h \theta)  +  \bA^{h}\bSigma_0  \exp(i h \theta) \Big\}\right)
\end{equation}
and
\begin{equation}
\label{Sigma0.vfar1}
\bSigma_0=\sigma^2\sum_{h=1}^\infty \bA^{h}(\bA^{h})^{\T},
\end{equation}
respectively.
For the VAR(1) model in (\ref{var1.ex}), the spectral density matrix and the covariance matrix of $\{\bx_t\}_{t \in \eZ}$ are
\begin{equation}
\label{density.var1}
f_{\bx,\btheta}=
\frac{1}{2\pi}\left(\bS_0  + \sum_{h = 1}^\infty \Big\{ \bS_0 (\bC^{\T})^{h} \exp(-i h \theta)  +  \bC^{h}\bS_0  \exp(i h \theta) \Big\}\right)
\end{equation}
and
\begin{equation}
\label{Sigma0.var1}
\bS_0=\sigma^2\sum_{h=1}^\infty \bC^{h}(\bC^{h})^{\T},
\end{equation}
respectively. Noting that $\bA^{h}\bSigma_0=\int_{\cU}\bpsi(u)\bC^{h}\bpsi(v')\bpsi(v')\bS_0\bpsi(v)dv'=\bpsi(u)\bC^{h}\bS_0\bpsi(v)$ and applying similar techniques, we can obtain that $f_{\bX,\theta}=\bpsi(u)f_{\bx,\theta}\bpsi(v)$ and $\bSigma_0=\bpsi(u)\bS_0\bpsi(v).$ The functional stability measure of $\{\bX_t\}_{t \in \eZ}$ under (\ref{vfar1.ex}) is
$$	2 \pi \cdot \underset{\theta\in [-\pi, \pi], \bPhi \in \cH_0}{\text{ess}\sup} \frac{\bPhi^{\T}f_{\bX,\theta}\bPhi}{\bPhi^{\T}\bSigma_0\bPhi}= 2 \pi \cdot \underset{\theta\in [-\pi, \pi], \bPhi \in \cH_0}{\text{ess}\sup} \frac{(\bpsi\bPhi)^{\T}f_{\bx,\theta}\bpsi\bPhi}{(\bpsi\bPhi)^{\T}\bS_0\bpsi\bPhi},$$
where $\bpsi\bPhi \in {\eR}^2$ with $(\bpsi\bPhi)_j=\langle\bphi_j,\bPhi\rangle, j=1,2.$ Hence the functional stability measure of $\{\bX_t\}_{t \in \eZ}$ under (\ref{vfar1.ex}) is the same as that of $\{\bx_t\}_{t \in \eZ}$ under (\ref{var1.ex}), i.e. the essential supremum of the maximal eigenvalue of $2\pi\bS_0^{-1/2}f_{\bx,\theta}\bS_0^{-1/2}$ over $\theta \in [-\pi,\pi].$ Some calculations yield $f_{\bx,\theta}$ and $\bS_0$ as follows.

By (\ref{density.var1}), we have
\begin{eqnarray*}
	f_{\bx,\btheta}&=&
	\frac{1}{2\pi}\left[\bS_0+ \bS_0\left(\begin{array}{cc} \sum_{h=1}^{\infty}a^{h}\exp(-ih\theta)& 0\\
		\sum_{h=1}^{\infty} h a^{h-1}\exp(-ih\theta)b &\sum_{h=1}^{\infty}a^{h}\exp(-ih\theta)
	\end{array}\right)\right.\\
	&& \left.+ \left(\begin{array}{cc} \sum_{h=1}^{\infty}a^{h}\exp(ih\theta)& \sum_{h=1}^{\infty} h a^{h-1}\exp(ih\theta)b\\
		0 &\sum_{h=1}^{\infty}a^{h}\exp(ih\theta)
	\end{array}\right)\bS_0\right]\\
	&=&\frac{1}{2\pi}\left[\bS_0+\bS_0\left(\begin{array}{cc} \frac{\alpha\exp(-i\theta)}{1-a\exp(-i\theta)}& 0\\
		\frac{b\exp(-i\theta)}{(1-a\exp(-i\theta))^2}&\frac{a\exp(-i\theta)}{1-a\exp(-i\theta)}
	\end{array}\right)
	+
	\left(\begin{array}{cc} \frac{a\exp(i\theta)}{1-a\exp(i\theta)}& \frac{b\exp(i\theta)}{(1-a\exp(i\theta))^2}\\
		0&\frac{a\exp(i\theta)}{1-a\exp(i\theta)}
	\end{array}\right)\bS_0
	\right].
\end{eqnarray*}

By (\ref{Sigma0.var1}), we have
\begin{eqnarray*}
	\bS_0&=&\left(\begin{array}{cc}
		\sum_{h=0}^{\infty} a^{2h} + \sum_{h=0}^{\infty}h^2a^{2h-2}b^2& \sum_{h=0}^{\infty} ha^{2h-1}b\\
		\sum_{h=0}^{\infty} ha^{2h-1}b & \sum_{h=0}^{\infty} a^{2h}\end{array}\right)\\
	&=&\left(\begin{array}{cc}\frac{1}{1-a^2} + \frac{(a^2+1)b^2}{(1-a^2)^3}& \frac{ab}{(1-a^2)^2}\\\frac{ab}{(1-a^2)^2} & \frac{1}{1-a^2}\end{array}\right).
\end{eqnarray*}
The right side of Figure~\ref{ill.example} plots functional stability measures of $\{\bX_t\}_{t \in \eZ}$ vs $b$ for different values of $a \in (0,1).$

\renewcommand{\theequation}{E.\arabic{equation}}
\section{Derivations for VFAR models}
\label{ap.deriv}
\subsection{Matrix representation of a VFAR($L$) model in (\ref{vfar3})}
\label{ap.vfar.deriv}
Note that the VFAR(L) model in (\ref{vfar1}) can be equivalently represented as
\begin{equation}
\label{vfar2}
X_{tj}(u) = \sum_{h = 1}^L \sum_{k=1}^p
\langle A_{jk}^{(h)}(u,\cdot), X_{(t-h)k}(\cdot)\rangle + \varepsilon_{tj}(u), ~~t = L+1,\ldots,n,j=1\dots,p.
\end{equation}
It then follows from the Karhunen-Lo\'eve expansion that (\ref{vfar2}) can be rewritten as
\begin{eqnarray*}
	\sum_{l=1}^{\infty} \xi_{tjl}\phi_{jl}(u) &=& \sum_{h=1}^{L} \sum_{k=1}^p \sum_{m=1}^{\infty} \langle A_{jk}^{(h)}(u,\cdot), \phi_{km}(\cdot)\rangle \xi_{(t-h)km} + \varepsilon_{tj}(u).
\end{eqnarray*}
This, together with orthonormality of $\{\phi_{jm}(\cdot)\}_{m \geq 1},$ implies that
$$	\xi_{tjl} = \sum_{h=1}^{L} \sum_{k=1}^p \sum_{m=1}^{q_k} \langle\phi_{jl},A_{jk}^{(h)}(\phi_{km}) \rangle\ \xi_{(t-h)km}+ r_{tjl} + \epsilon_{tjl},
$$
where $r_{tjl}=\sum_{h=1}^{L} \sum_{k=1}^p \sum_{m=q_k+1}^{\infty} \langle\phi_{jl},A_{jk}^{(h)}(\phi_{km}) \rangle  \xi_{(t-h)km}$ and
$\epsilon_{tjl}=\langle \phi_{jl},\varepsilon_{tj}\rangle$ for $l=1,\dots, q_j,$ represent the approximation and random errors, respectively.
Let $\br_{tj}=(r_{tj1}, \dots, r_{tjq_j})^{\T}$ and  $\beps_{tj}=(\epsilon_{tj1}, \dots, \epsilon_{tjq_j})^{\T}.$ Let $\bR_j, \bE_j$ be $(n-L) \times q_j$ matrices whose row vectors are formed by $\{\br_{tj}, t=L+1, \dots, n\}$ and $\{\beps_{tj}, t=L+1, \dots, n\}$ respectively. Then (\ref{vfar2}) can be represented in the matrix form of (\ref{vfar3})

\subsection{VFAR(1) representation of a VFAR($L$) model}
\label{ap.vfar1.rep}
We can represent a $p$-dimensional VFAR($L$) model in (\ref{vfar1}) as a $pL$-dimensional VFAR(1) model in the form of
\begin{equation}
    \label{vfar.eqv}
    \widetilde\bX_t(u)=\int_{\cU}\widetilde \bA_{1}(u,v)\widetilde\bX_{t-1}(v)dv+ \widetilde\bvarepsilon_{t-1}(u), ~u \in \cU,
\end{equation}
where $\widetilde \bX_{t}=\left(\begin{array}{c}\bX_{t}\\\bX_{t-1}\\\vdots\\\bX_{t-L+1}\end{array}\right),$ $\widetilde \bA_1=\left(\begin{array}{ccccc}\bA_{1} & \bA_{2} & \cdots & \bA_{L-1} & \bA_{L}\\ \bI_{p} & {\bf 0} & \cdots & {\bf 0} & {\bf 0}\\{\bf 0} & \bI_{p} & \cdots & {\bf 0} & {\bf 0}\\\vdots & \vdots & \ddots & \vdots & \vdots\\{\bf 0} & {\bf 0} & \cdots & \bI_{p} & {\bf 0}\end{array}\right),$ $\widetilde \bvarepsilon_{t}=\left(\begin{array}{c}\bvarepsilon_{t}\\
\bvarepsilon_{t-1}\\\vdots\\ \bvarepsilon_{t-L+1}\end{array}\right)$ and
$\bI_p$  denotes the identity operator.
In the non-functional setting, a similar VAR(1) representation of a VAR($L$) model can be found in {{\color{blue} Basu and Michailidis (2015)}}.

\renewcommand{\theequation}{F.\arabic{equation}}
\subsection{VFAR(1) representation of the simulation example}
\label{ap.vfar1.sim}
Noting that $\btheta_t=\bB\btheta_{t-1}+\bfeta_t,$ we have $\btheta_{tj}=\sum_{k=1}^p\bB_{jk}\btheta_{(t-1)k}+\bfeta_{tj}$ for $j=1,\dots,p.$ Multiplying both sides by $\bs(u)^{\T}$ and applying $\int_{\cU}\bs(v)\bs(v)^{\T}dv=\bI_5,$ we obtain that $\bs(u)^{\T}\btheta_{tj}=\int_{\cU}\sum_{k=1}^p\bs(u)^{\T}\bB_{jk}\bs(v)\bs(v)^{\T}\btheta_{(t-1)k}dv+ \bs(u)^{\T}\bfeta_{tj}.$ Letting $A_{jk}(u,v)=\bs(u)^{\T}\bB_{jk}\bs(v),$ $X_{tj}(u)=\bs(u)^{\T}\btheta_{tj}$ and $\varepsilon_{tj}(u)=\bs(u)^{\T}\bfeta_{tj},$ we have $X_{tj}(u)=\sum_{k=1}^p\langle A_{jk}(u,\cdot)X_{(t-1)k}(\cdot)\rangle + \varepsilon_{tj}(u).$

\section{Algorithms in fitting VFAR models}
\label{ap.alg}

\subsection{Selection of tuning parameters}
\label{ap.select.tune}
To fit the proposed sparse VFAR model, we need choose values for three tuning parameters, $q_j$ (the number of selected principal components for $j=1,\dots,p$), $\eta_j$ (the smoothing parameter when performing regularized FPCA, as described in Section~\ref{ap.rfpca} below) and $\gamma_{nj}$ (the regularization parameter in (\ref{vfar.crit}) to control the block sparsity level in $\{\widehat\bPsi_{jk}^{(h)}: h=1,\dots,L, k=1,\dots,p\}$).

We adopt a $K$-fold cross-validated method to choose $(q_j,\eta_j)$ for each $j.$ Specifically, let $W_{tjs}$ be observed values of $X_{tj}(u_s)$ at $u_1, \dots, u_T.$ We randomly divide the set $\{1, \dots, n\}$ into $K$ groups, $\cD_1, \dots, \cD_K$ of approximately equal size, with the first group treated as a validation set. Implementing regularized FPCA on the remaining $K-1$ groups, we obtain estimated mean function $\widehat\mu_{jl}^{(-1)}(u),$ FPC scores $\widehat\xi_{tjl,\eta_j}^{(-1)}$ and eigenfunctions $\widehat \phi_{jl}^{(-1)}(u;\eta_j)$ for $l=1, \dots, q_j.$ The predicted curve for the $t$-th sample in group one can be computed by $\widehat W_{tjs}^{(1)}=\widehat\mu_{jl}^{(-1)}(u_s)+ \sum_{l=1}^{q_j}\widehat\xi_{tjl,\eta_j}^{(-1)}\widehat \phi_{jl}^{(-1)}(u_s;\eta_j)$. This procedure is repeated $K$ times. Finally, we choose $q_j$ and $\eta_j$ as the values that minimize the mean cross-validated error, $$\text{CV}(q_j,\eta_j)=(KT)^{-1}\sum_{k=1}^K \sum_{s=1}^T\sum_{t \in \cD_k} (W_{tjs}-\widehat W_{tjs}^{(k)})^2.$$

The optimal $\gamma_{nj}$'s are selected by minimizing AICs or BICs in (\ref{ic_vfar}), where details can be found in Section~\ref{sec.sim}. 

\subsection{Regularized FPCA}
\label{ap.rfpca}
In this section, we drop subscripts $j$ for simplicity of notation. Suppose we observe $\bX(\cdot)=(X_1(\cdot),\dots, X_n(\cdot))^{\T}$ on $\cU,$ our goal is to find the first $q$ regularized principal component functions $\{ \phi_{l}(\cdot), l=1, \dots,q\}.$ We obtain the $l$-th leading principal component $\phi_l(\cdot)$ through a smoothing approach, which maximizes the following penalized sample variance [(9.1) in {\color{blue} Ramsay and Silverman (2005)}]
\begin{equation}
\label{pen.fpca}
\text{PEN}_{\eta}(\phi_l)=\frac{\var(\langle\phi_l, X_i\rangle)}{||\phi_l||^2 + \eta ||\phi''_l||^2},
\end{equation}
subject to $||\phi_l||=1$ and $\langle\phi_l,\phi_{l'}\rangle+\eta\langle\phi_l'',\phi_{l'}''\rangle=0, l'=1,\dots,l-1,$ where $\eta \geq 0$ is a smoothing parameter to control the roughness of $\phi_l(\cdot).$

Suppose that $\bX(u)=\bdelta^{\T}\bbb(u)$ and $\phi_l(u)=\bzeta_l^{\T}\bbb(u)$ where $\bbb(\cdot)$ is a $G$-dimensional B-spline basis function, $\bdelta \in {\eR}^{n\times G}$
and $\bzeta_l \in {\eR}^G$ are the basis coefficients for $\bX(\cdot)$ and $\phi_l(\cdot),$ respectively. Let $\bJ=\int\bbb(u)\bbb(u)^{\T}du,$ $\bU=\bJ\bdelta^{\T}\bdelta\bJ$ and $\bQ=\int\bbb''(u)\bbb''(u)^{\T}du,$ (\ref{pen.fpca}) is equivalent to maximizing
\begin{equation}
\label{pen.fpca1}
\text{PEN}_{\eta}(\phi_l)=\frac{\bzeta_l^{\T}\bU\bzeta_l}{\bzeta_l^{\T}(\bJ+\eta\bQ)\bzeta_l},
\end{equation}
subject to $\bzeta_l^{\T}\bJ\bzeta_l=1$ and $\bzeta_l^{\T}(\bJ+\eta \bQ)\bzeta_{l'}=0, l'=1,\dots, l-1.$ By {singular value decomposition} (SVD), we obtain eigen-pairs, $(\bS_1,\bP_1)$ and $(\bS_2,\bP_2)$ such that $\bJ+\eta\bQ=\bP_1\bS_1^{-2}\bP_1^{\T}$ and $\bS_1\bP_1^{\T}\bU\bP_1\bS_1=\bP_2\bS_2^{-2}\bP_2^{\T}.$ Then (\ref{pen.fpca1}) becomes
$\text{PEN}_{\eta}(\phi_l)=\frac{\bx_l^{\T}\bP_2^{\T}\bS_2\bP_2\bx_l}{\bx_l^{\T}\bx_l},$ where $\bx_l=\bS_1^{-1}\bP_1^{\T}\bzeta_l.$
This suggests us to perform SVD on $\bP_2^{\T}\bS_2\bP_2,$ where we can obtain $\widehat\bx_l,$
$\widehat\bzeta_l=\bP_1\bS_1\widehat\bx_l$ and $\widehat\phi_l(u)=\widehat\bzeta_l^{\T}\bbb(u)/\big(\widehat\bzeta_l^{\T}\bJ\widehat\bzeta_l\big)^{1/2},l=1,\dots, q.$ In practice, we can set $G$ to a pre-specified large enough value, and implement the cross-validation procedure as described in Section~\ref{ap.select.tune} to select $q$ and $\eta.$

\subsection{Block FISTA algorithm to solve (\ref{vfar.crit})}
\label{ap.alg.vfar}

The optimization problem in (\ref{vfar.crit}) can be reformulated as follows.
\begin{eqnarray}
\label{generalproblem}
\min_{\bX \in {\eR}^{r \times q_j}} g(\bX), ~~g(\bX) = f(\bX) + \gamma_{nj} \sum_{k = 1}^{pL} \|\bX_k\|_F,
\end{eqnarray}
where $f(\bX) = 2^{-1}\text{trace}\left\{(\bY - \bB \bX)^{\T} (\bY - \bB \bX)\right\},$ $r = \sum_{h=1}^L\sum_{k=1}^p q_k,$ $\bY \in \eR^{(n-L)\times q_j},$ $\bB \in \eR^{(n-L) \times r},$ and $\bX = (\bX_1^{\T},\ldots,\bX_{pL}^{\T})^{\T} \in \eR^{r \times q_j}$ with $\bX_k \in \eR^{q_k \times q_j}$ for $k=1,\dots,p.$
(\ref{generalproblem}) is a convex problem including the smooth part for $\bX_k,$ i.e. $f(\bX)$ and the non-smooth part for $\bX_k$, i.e.$\gamma_{nj} \sum_{k = 1}^{pL} \|\bX_k\|_F$. To solve the minimization problem in (\ref{generalproblem}), we adopt a block version of {fast iterative shrinkage-thresholding algorithm} (FISTA) ({\color{blue}Beck and Teboulle, 2009}) combined with a restarting technique ({\color{blue}O'Donoghue and Candes, 2015}), namely block FISTA.

The basic idea behind our proposed block FISTA is summarized as follows. Let $\nabla f(\bX)$ be the gradient of $f(\bX)$ at $\bX.$ We start with an initial value $\bX^{(0)}.$ At the $(m+1)$-th iteration we first try to solve a regularized sub-problem
\begin{eqnarray}
\label{subprob-1}
\min_{\bX \in \eR^{r \times q_j}} \text{trace}\left\{(\nabla f(\bX^{(m)})^{\T} (\bX -  \bX^{(m)})\right\}
+ (2C)^{-1}\left \|\bX - \bX^{(m)} \right \|_F^2
+ \gamma_{nj} \sum_{k = 1}^{pL} \|\bX_k\|_F,
\end{eqnarray}
where $\bX^{(m)}$ is the $m$-th iterate and $C>0$ is a small constant controlling the stepsize at $(m+1)$-th step. The second term in (\ref{subprob-1}) can be interpreted as a quadratic regularization, which restricts the updated iterate not to be very far from $\bX^{(m)}$. 
The analytical solution to (\ref{subprob-1}) takes the form of
\begin{eqnarray}
\label{thresholdsolution}
\widetilde{\bX}^{(m+1)} = \big(\widetilde{\bX}_k^{(m+1)}\big)~~ \mbox{with}~~\widetilde{\bX}^{(m+1)}_k =  \left( 1 - \gamma_{nj} C{\|\bZ_k^{(m)}\|_F^{-1}}\right)_{+} \bZ_k^{(m)},k=1, \dots, pL,
\end{eqnarray}
where $\bZ^{(m)} = \bX^{(m)} - C\nabla f(\bX^{(m)})=\big(\big(\bZ_1^{(m)}\big)^{\T}, \dots, \big(\bZ_{pL}^{(m)}\big)^{\T}\big)^{\T}$ and $x_+ = \max(0,x).$ (See also (3.a) and (3.b) of Algorithm~\ref{alg1}).


We then take block FISTA ({\color{blue}Beck and Teboulle, 2009}) by adding  an extrapolation step in the algorithm (see also (3.c) and (3.d) of Algorithm~\ref{alg1}):
$$
\bX^{(m+1)} = \widetilde{\bX}^{(m+1)} + \omega^{(m+1)}(\widetilde{\bX}^{(m+1)} - \widetilde{\bX}^{(m)}),
$$
where the weight $\omega^{(m+1)}$ is specified in Algorithm~\ref{alg1}.  Finally, at the end of each iteration, we evaluate the generalized gradient at $\bX^{(m+1)}$ by computing the sign of
$$
\text{trace}\left\{(\bX^{(m)} - \widetilde{\bX}^{(m+1)})^{\T}(\widetilde{\bX}^{(m+1)} - \widetilde{\bX}^{(m)})\right\},
$$
which can be thought of a proxy of $\text{trace}\left\{(\nabla g(\bX^{(m)}))^{\T} (\widetilde{\bX}^{(m+1)} -  \widetilde{\bX}^{(m)})\right\}$.
For a positive sign, i.e. the objective function is increasing at $\widetilde{\bX}^{(m+1)}$, we then restart our accelerated algorithm by setting $\bX^{(m+1)} = \bX^{(m)}$ and $ \omega^{(m+1)} = \omega^{(1)}$ ({\color{blue}O'Donoghue and Candes, 2015}). This step can guarantee that the objective function $g$ decreases over each iteration. We iterative the above steps until convergence. We summarize the restarting-based block FISTA in Algorithm~\ref{alg1}. In practice, one issue is how to choose the stepsize parameter $C$. In general, the proposed scheme is guaranteed to converge when $C < \big(\lambda_{\max}(\bB^{\T}\bB)\big)^{-1}.$ 
Here we choose $C = 0.9 \big(\lambda_{\max}(\bB^{\T}\bB)\big)^{-1}$, which turns to work well in empirical studies. Alternatively, $C$ can be selected through a line search and one simple backtracking rule.

\begin{center}
	\begin{algorithm}[!t] \caption{\label{alg1}{Block FISTA for solving (\ref{generalproblem})}}
		\begin{enumerate}
			\item[1.] Input: $C = 0.9 \big(\lambda_{\max}(\bB^{\T}\bB)\big)^{-1}$, $\theta_0 = 1$, $\bX^{(0)} = (\bX_1^{(0)\T},\ldots,\bX_{pL}^{(0)\T})^{\T}= {\bf 0},$
			$\bZ^{(0)} = (\bZ_1^{(0)\T},\ldots,\bZ_{pL}^{(0)\T})^{\T}= {\bf 0},$ $\widetilde{\bX}^{(0)} = (\widetilde{\bX}_1^{(0)\T},\ldots,\widetilde{\bX}_{pL}^{(0)\T})^{\T}= {\bf 0}.$
			\item[2.] For $m = 0, 1, \ldots$ do 	
			\begin{itemize}		
				\item[] (3.a)~ $\bZ^{(m)} = \bX^{(m)} - C\nabla f(\bX^{(m)}),$
				\item[] (3.b)~ $\widetilde{\bX}^{(m+1)}_k =  \left( 1 - \gamma_{nj} C{\|\bZ_k^{(m)}\|_F^{-1}}\right)_+ \bZ_k^{(m)},k=1, \dots, pL,$
				\item[] (3.c)~ $\theta_{m+1} = \big(1+ \sqrt{1+ 4 \theta_m^2}\big)/2,$
				\item[] (3.d)~ $\bX^{(m+1)} = \widetilde{\bX}^{(m+1)} + \frac{\theta_m - 1}{\theta_{m+1}}\big(\widetilde{\bX}^{(m+1)} - \widetilde{\bX}^{(m)}\big),$		
				\item[] (3.e)~ If $\text{trace}\left\{(\bX^{(m)} - \widetilde{\bX}^{(m+1)})^{\T}(\widetilde{\bX}^{(m+1)} - \widetilde{\bX}^{(m)})\right\} > 0$, set
				$$
				\bX^{(m+1)} = \bX^{(m)}, \theta_{m+1} = 1.
				$$
			\end{itemize}
			\item[] end do until convergence.
			\item[3.] Output: the final estimator $\bX^{(m+1)}.$
		\end{enumerate}
	\end{algorithm}		
\end{center}

\section{Additional empirical results}
\label{ap.emp}

\subsection{Simulation studies}
\label{ap.sim}
\begin{figure*}
	\centering
	\includegraphics[width=0.75\textwidth]{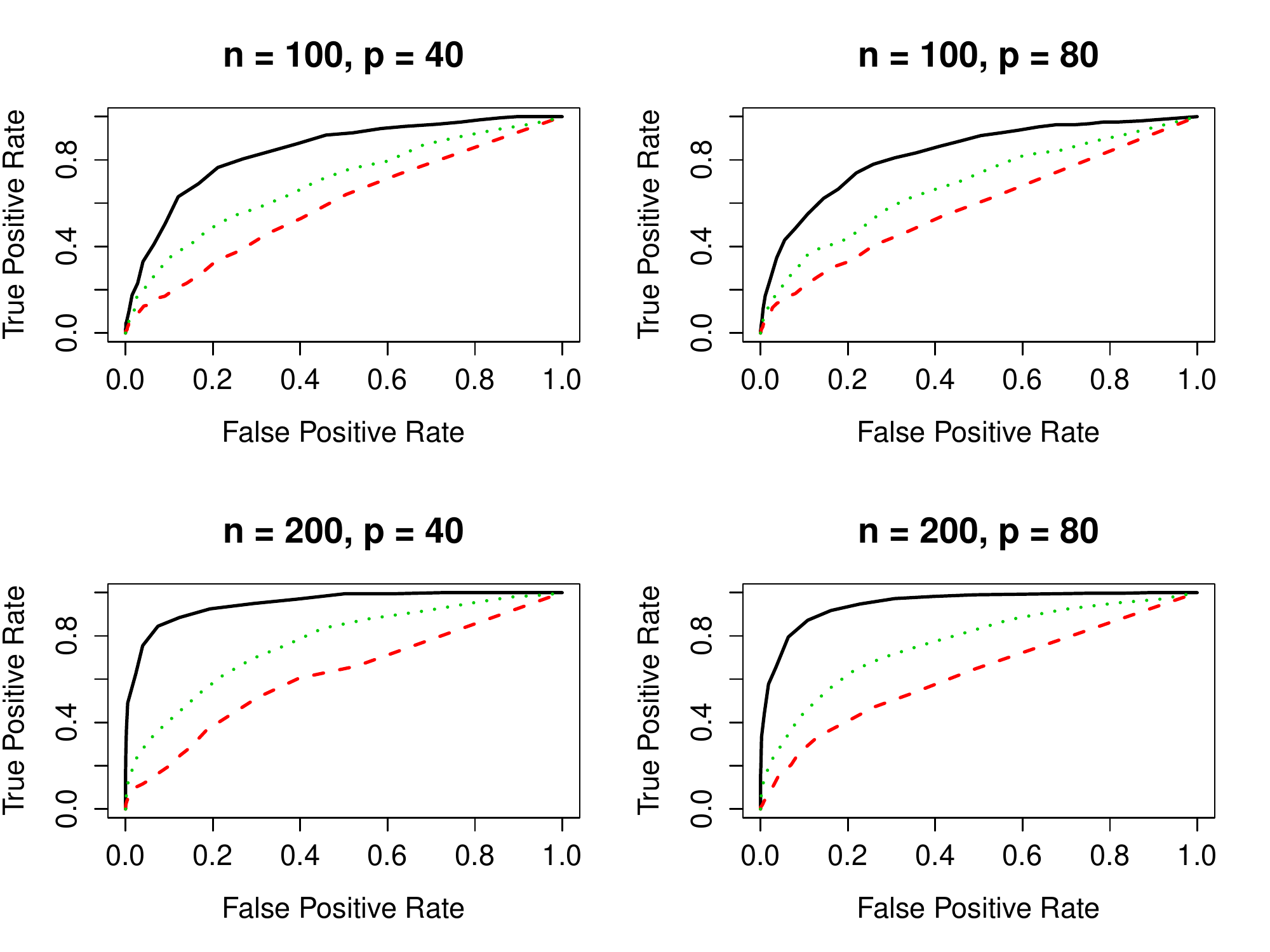}
	\caption{\label{roc.curve.sparse}{Comparisons of median estimated ROC curves over 100 simulation runs. $\ell_1/\ell_2$-$\LS_{\text{a}}$  (black solid), $\ell_1/\ell_2$-$\LS_2$  (green dotted) and $\ell_1$-$\LS_1$  (red dashed) for Model~(i).}}
\end{figure*}

\begin{figure*}
	\centering
	\includegraphics[width=0.75\textwidth]{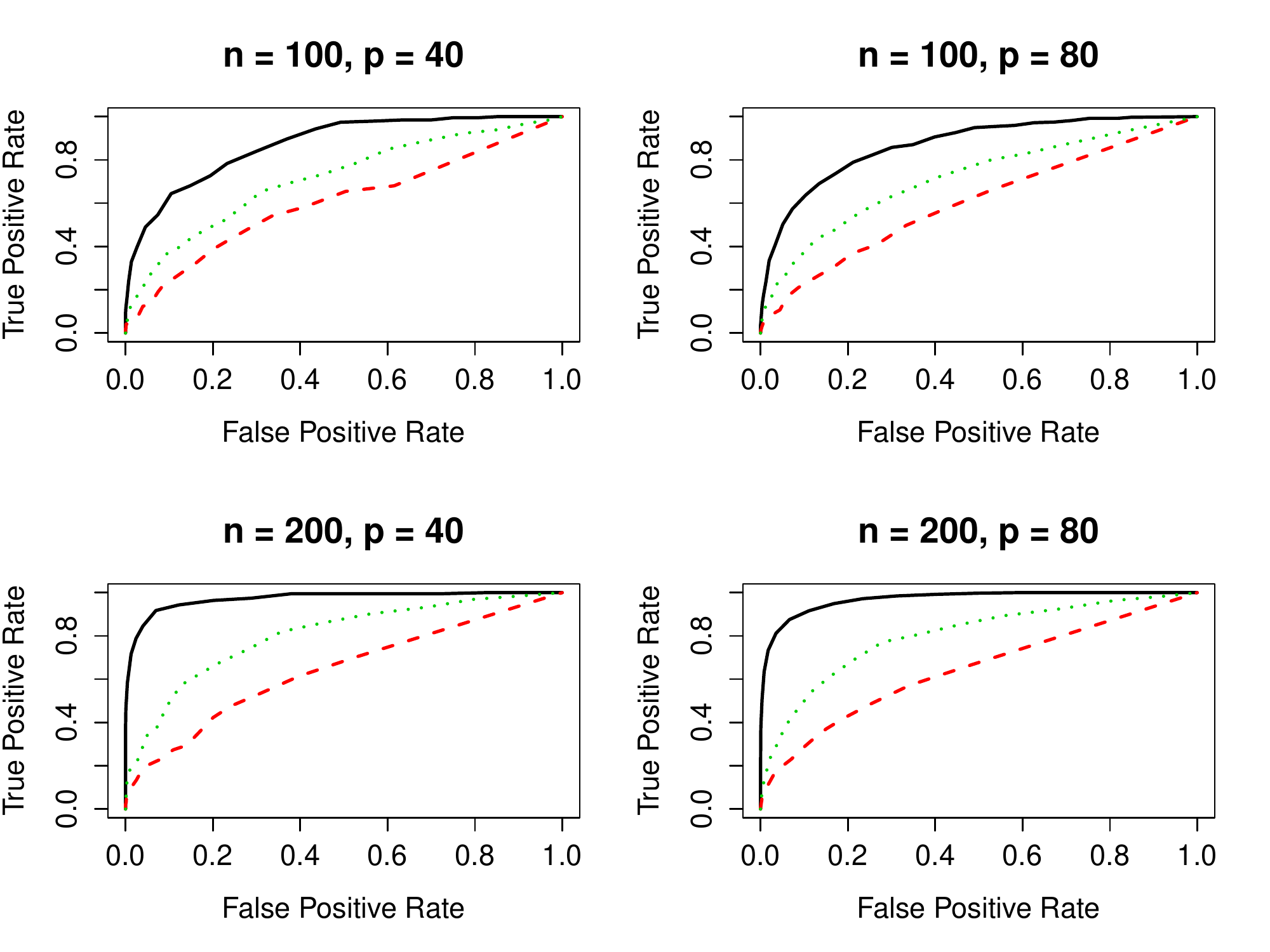}
	\caption{\label{roc.curve.band}{Comparisons of median estimated ROC curves over 100 simulation runs. $\ell_1/\ell_2$-$\LS_{\text{a}}$  (black solid), $\ell_1/\ell_2$-$\LS_2$  (green dotted) and $\ell_1$-$\LS_1$  (red dashed) for Model~(ii).}}
\end{figure*}

Figures~\ref{roc.curve.sparse} and \ref{roc.curve.band} plot the median best ROC curves (we rank ROC curves by the corresponding AUROCs) over the 100 stimulation runs in Models~(i) and (ii), respectively. Again we see that $\ell_1/\ell_2$-$\text{LS}_2,$ which explains the partial curve information, although performing better than $\ell_1$-$\text{LS}_1$ is substantially outperformed by $\ell_1/\ell_2$-$\text{LS}_{\text{a}}$ in terms of model selection consistency.

\subsection{Real data analysis}
\label{supp_sec_real}
\begin{figure*}
	\centering
	\includegraphics[width=7.00cm,height=7.00cm]{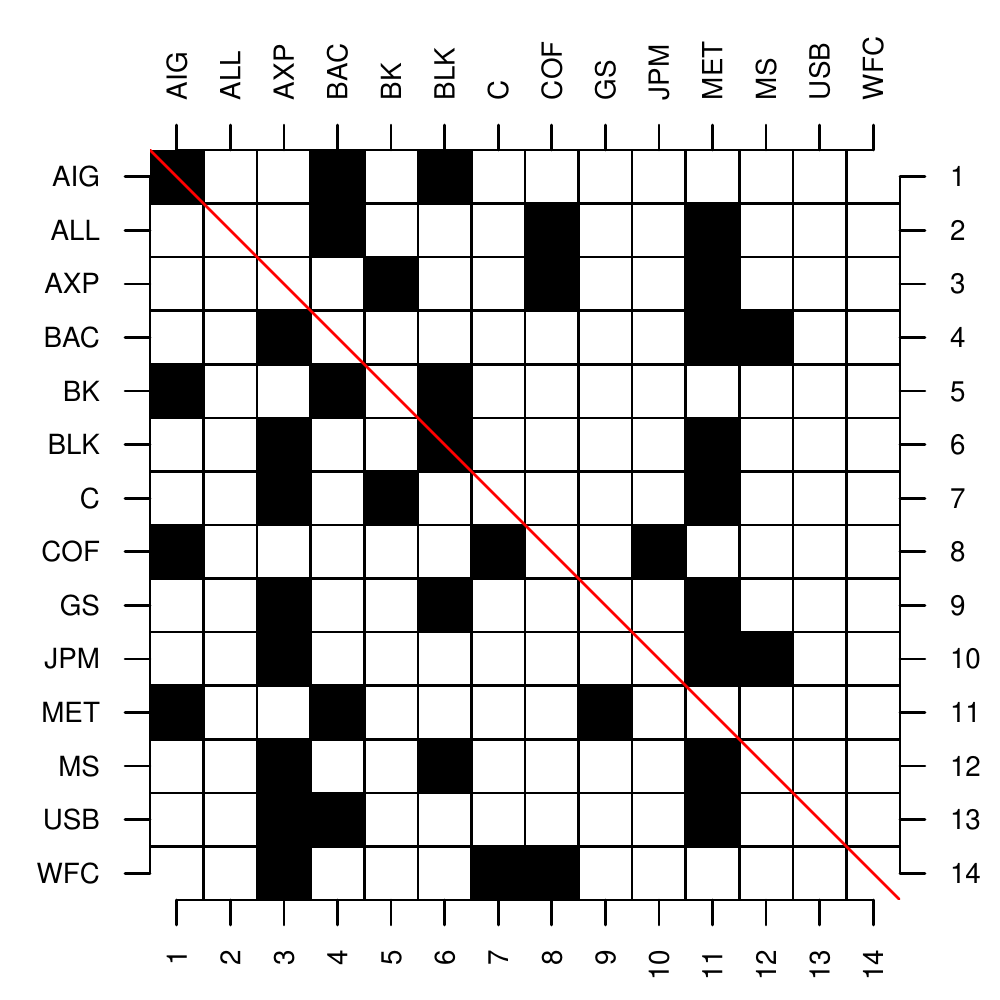}\includegraphics[width=7.00cm,height=7.00cm]{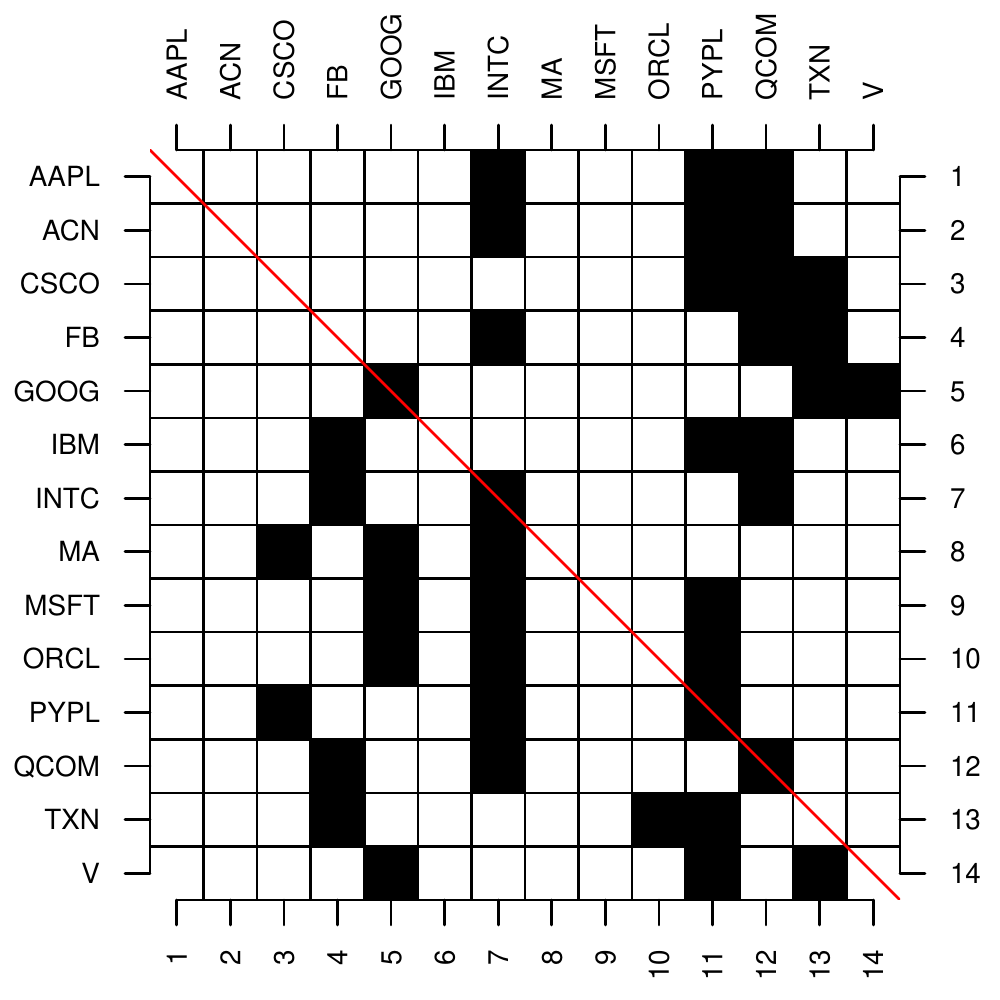}
	\includegraphics[width=9.00cm,height=9.00cm]{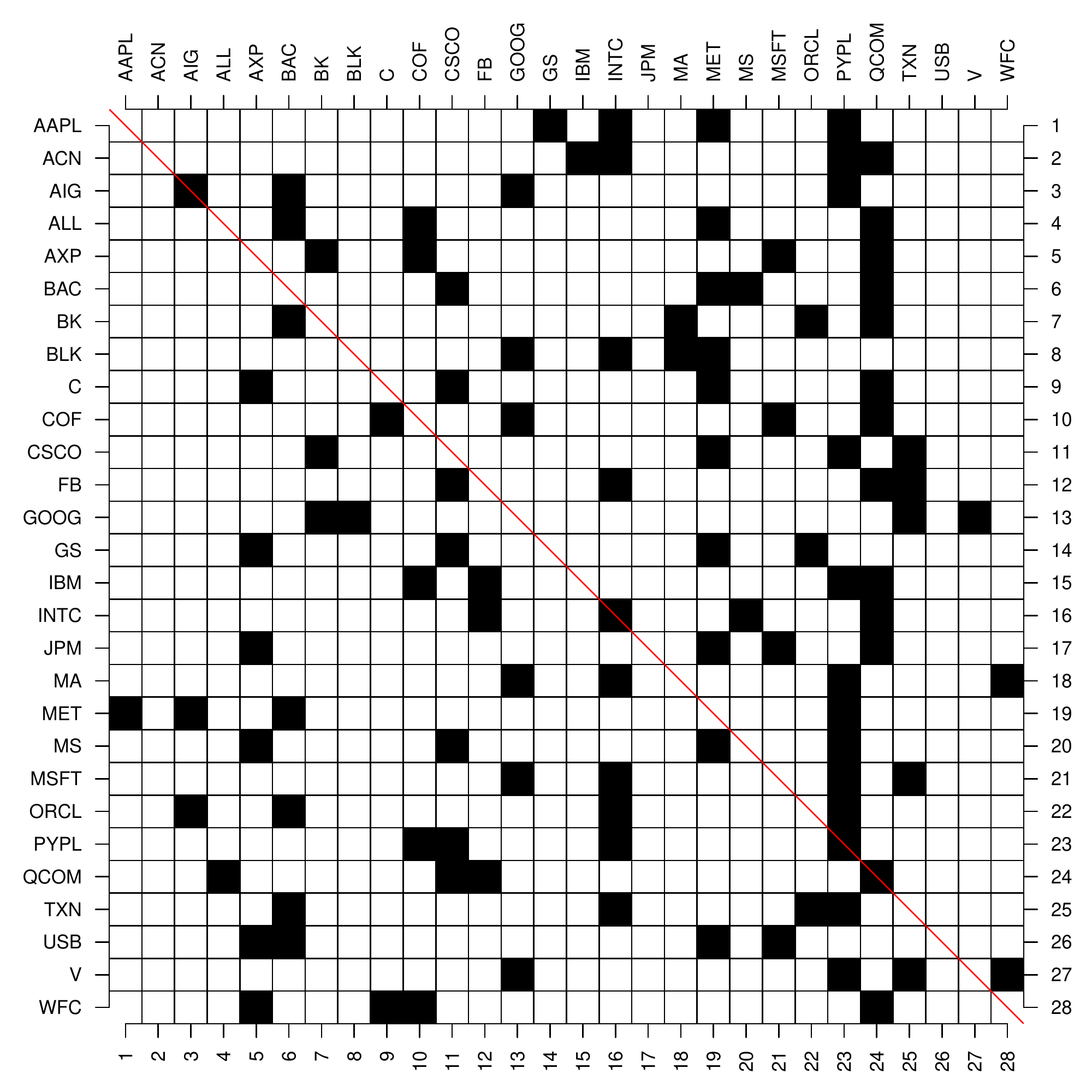}
	\caption{\label{A.plot}{\it \small Topleft, topright and bottom graphs plot the sparsity structures in {$\widehat\bA$} for stocks in only the financial sector, only the IT sector and both financial and IT sectors, respectively. Black and white correspond to non-zero and zero functional components in $\widehat\bA,$ respectively.}}
\end{figure*}

\begin{figure*}
	\centering
	\includegraphics[width=12.50cm,height=12.50cm]{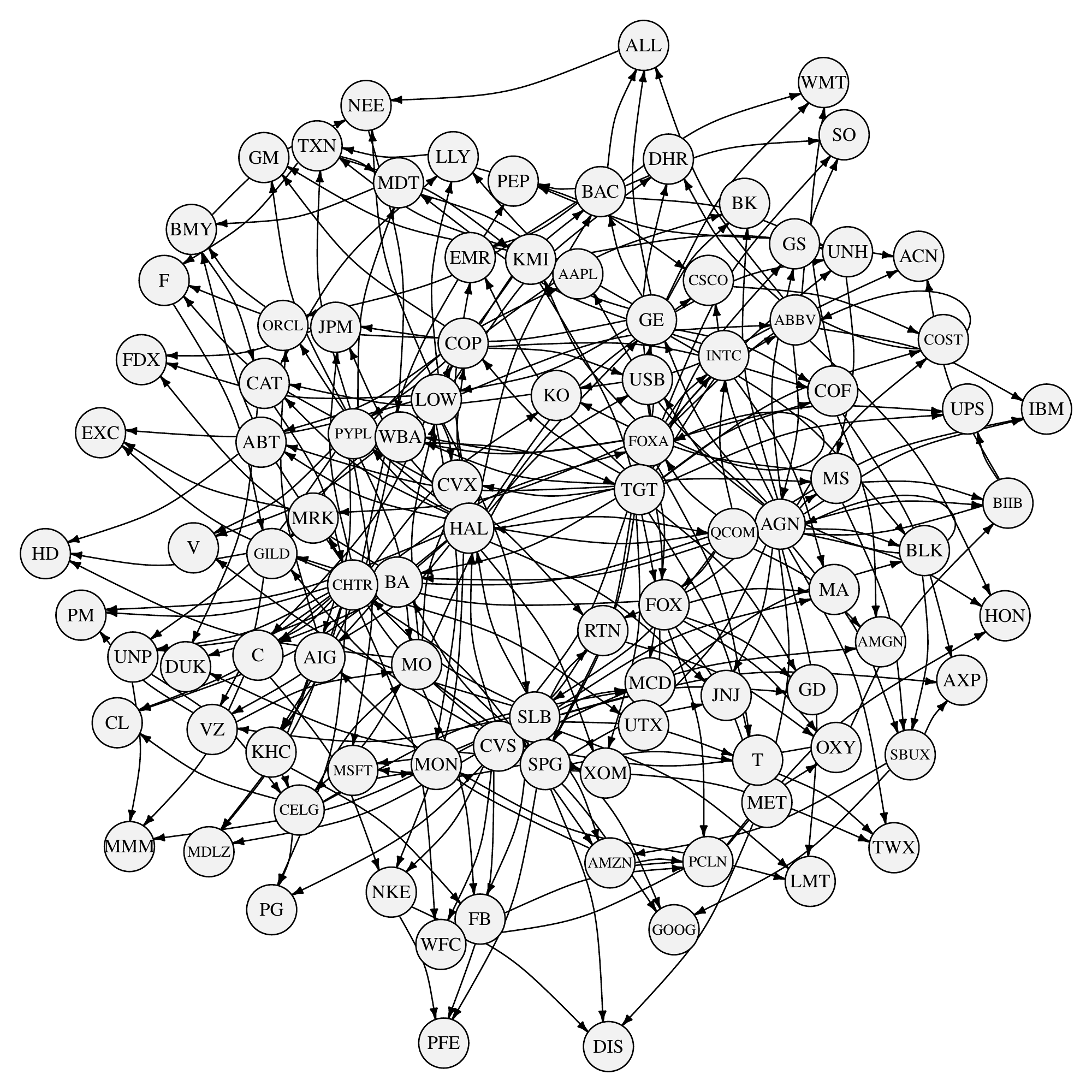}
	\caption{\label{net98.plot}{\it \small The directed graph with indegree=3 for $p=98$ stocks.}}
\end{figure*}

Table~\ref{name.tb} provides tickers, company names and classified sectors of 98 stocks under our study. Figure~\ref{A.plot} plots the sparsity patterns in $\widehat\bA$ (estimated transition function) for either 18 or 36 stocks. For a large collection of $p=98$ S\&P100 stocks, to better visualize a large causal network, we set the row-wise sparsity to $3/98.$ We then plot a large and sparse directed graph in Figure~\ref{net98.plot}. 
We observe that companies, e.g. Allergan, Halliburton, Target Corp., have relatively higher causal impacts on all 98 stocks in terms of their CIDR curves.


\begin{table}
	\caption{\label{name.tb} List of S\&P~100 stocks under study.}
	\vspace{0.5cm}
\linespread{1.25}
\begin{adjustbox}{angle=90}
\tiny
\begin{tabular}{lll|lll}
\hline
Ticker & Company name & Sector & Ticker & Company name & Sector \\ \hline
AAPL & APPLE INC & Information Technology & JPM & JPMORGAN CHASE \& CO & Financials \\
ABBV & ABBVIE INC & Health Care & KHC & KRAFT HEINZ & Consumer Staples \\
ABT & ABBOTT LABORATORIES & Health Care & KMI & KINDER MORGAN INC & Energy \\
ACN & ACCENTURE PLC CLASS A & Information Technology & KO & COCA-COLA & Consumer Staples \\
AGN & ALLERGAN & Health Care & LLY & ELI LILLY & Health Care \\
AIG & AMERICAN INTERNATIONAL GROUP INC & Financials & LMT & LOCKHEED MARTIN CORP & Industrials \\
ALL & ALLSTATE CORP & Financials & LOW & LOWES COMPANIES INC & Consumer Discretionary \\
AMGN & AMGEN INC & Health Care & MA & MASTERCARD INC CLASS A & Information Technology \\
AMZN & AMAZON COM INC & Consumer Discretionary & MCD & MCDONALDS CORP & Consumer Discretionary \\
AXP & AMERICAN EXPRESS & Financials & MDLZ & MONDELEZ INTERNATIONAL INC CLASS A & Consumer Staples \\
BA & BOEING & Industrials & MDT & MEDTRONIC PLC & Health Care \\
BAC & BANK OF AMERICA CORP & Financials & MET & METLIFE INC & Financials \\
BIIB & BIOGEN INC INC & Health Care & MMM & 3M & Industrials \\
BK & BANK OF NEW YORK MELLON CORP & Financials & MO & ALTRIA GROUP INC & Consumer Staples \\
BLK & BLACKROCK INC & Financials & MON & MONSANTO & Materials \\
BMY & BRISTOL MYERS SQUIBB & Health Care & MRK & MERCK \& CO INC & Health Care \\
C & CITIGROUP INC & Financials & MS & MORGAN STANLEY & Financials \\
CAT & CATERPILLAR INC & Industrials & MSFT & MICROSOFT CORP & Information Technology \\
CELG & CELGENE CORP & Health Care & NEE & NEXTERA ENERGY INC & Utilities \\
CHTR & CHARTER COMMUNICATIONS INC CLASS A & Consumer Discretionary & NKE & NIKE INC CLASS B & Consumer Discretionary \\
CL & COLGATE-PALMOLIVE & Consumer Staples & ORCL & ORACLE CORP & Information Technology \\
COF & CAPITAL ONE FINANCIAL CORP & Financials & OXY & OCCIDENTAL PETROLEUM CORP & Energy \\
COP & CONOCOPHILLIPS & Energy & PCLN & THE PRICELINE GROUP INC & Consumer Discretionary \\
COST & COSTCO WHOLESALE CORP & Consumer Staples & PEP & PEPSICO INC & Consumer Staples \\
CSCO & CISCO SYSTEMS INC & Information Technology & PFE & PFIZER INC & Health Care \\
CVS & CVS HEALTH CORP & Consumer Staples & PG & PROCTER \& GAMBLE & Consumer Staples \\
CVX & CHEVRON CORP & Energy & PM & PHILIP MORRIS INTERNATIONAL INC & Consumer Staples \\
DHR & DANAHER CORP & Health Care & PYPL & PAYPAL HOLDINGS INC & Information Technology \\
DIS & WALT DISNEY & Consumer Discretionary & QCOM & QUALCOMM INC & Information Technology \\
DUK & DUKE ENERGY CORP & Utilities & RTN & RAYTHEON & Industrials \\
EMR & EMERSON ELECTRIC & Industrials & SBUX & STARBUCKS CORP & Consumer Discretionary \\
EXC & EXELON CORP & Utilities & SLB & SCHLUMBERGER NV & Energy \\
F & F MOTOR & Consumer Discretionary & SO & SOUTHERN & Utilities \\
FB & FACEBOOK CLASS A INC & Information Technology & SPG & SIMON PROPERTY GROUP REIT INC & Real Estate \\
FDX & FEDEX CORP & Industrials & T & AT\&T INC & Telecommunications \\
FOX & TWENTY-FIRST CENTURY FOX INC CLASS & Consumer Discretionary & TGT & TARGET CORP & Consumer Discretionary \\
FOXA & TWENTY FIRST CENTURY FOX INC CLASS & Consumer Discretionary & TWX & TIME WARNER INC & Consumer Discretionary \\
GD & GENERAL DYNAMICS CORP & Industrials & TXN & TEXAS INSTRUMENT INC & Information Technology \\
GE & GENERAL ELECTRIC & Industrials & UNH & UNITEDHEALTH GROUP INC & Health Care \\
GILD & GILEAD SCIENCES INC & Health Care & UNP & UNION PACIFIC CORP & Industrials \\
GM & GENERAL MOTORS & Consumer Discretionary & UPS & UNITED PARCEL SERVICE INC CLASS B & Industrials \\
GOOG & ALPHABET INC CLASS C & Information Technology & USB & US BANCORP & Financials \\
GS & GOLDMAN SACHS GROUP INC & Financials & UTX & UNITED TECHNOLOGIES CORP & Industrials \\
HAL & HALLIBURTON & Energy & V & VISA INC CLASS A & Information Technology \\
HD & HOME DEPOT INC & Consumer Discretionary & VZ & VERIZON COMMUNICATIONS INC & Telecommunications \\
HON & HONEYWELL INTERNATIONAL INC & Industrials & WBA & WALGREEN BOOTS ALLIANCE INC & Consumer Staples \\
IBM & INTERNATIONAL BUSINESS MACHINES CO & Information Technology & WFC & WELLS FARGO & Financials \\
INTC & INTEL CORPORATION CORP & Information Technology & WMT & WALMART STORES INC & Consumer Staples \\
JNJ & JOHNSON \& JOHNSON & Health Care & XOM & EXXON MOBIL CORP & Energy \\ \hline
\end{tabular}
\end{adjustbox}
\end{table}
\linespread{1.5}

\newpage
\linespread{1.15}\selectfont
\section*{References}
\begin{description}
	\item
	Basu, S.  and Michailidis, G. (2015). Regularized estimation in sparse high-dimensional time series models. {\it The Annals of Statistics}, {\bf 43}, 1535-1567.
	\item
	Beck, A.  and Teboulle, M. (2009). A fast iterative shrinkage-thresholding algorithm for linear inverse problems. {\it SIAM J. Imaging Sciences}, {\bf 2}, 183-202.
	\item
	Boucheron, S., Lugosi, G. and Massart, P. (2014). {\it Concentration Inequalities: A Nonasymptotic Theory of Independence}. Oxford University Press.
	\item
	Hall, P.  and Horowitz, J. L. (2007). Methodology and convergence rates for functional linear regression. {\it The Annals of Statistics}, {\bf 35}, 70-91.
	\item O'Donoghue, B. and Candes, E. (2015). Adaptive restart for accelerated gradient schemes. {\it Foundations of Computational Mathematics}, {\bf 15}, 715-732.
	\item
	Qiao, X., Guo, S. and James, G. M. (2019). Functional graphical models. {\it Journal of the American Statistical Association}, {\bf 114}, 211-222.
	\item 
	Ramsay, J. O. and Silverman, B. W. (2005). {\it Functional data analysis (2nd et.)}. Springer, New York.
	\item 
	Rudelson, M. and Vershynin, R. (2014). Hanson-wright inequality and sub-Gaussian concentration. {\it Electronic Communications in Probability}, {\bf 18}, 1-9.
\end{description}
\end{document}